# ISSUES IN ANTENNA OPTIMIZATION - A MONOPOLE CASE STUDY

**Richard A. Formato**

Registered Patent Attorney & Consulting Engineer
Cataldo & Fisher, LLC, 400 TradeCenter, Suite 5900
Woburn, MA 01801 USA
rf2@ieee.org, rformato@cataldofisher.com

**Abstract —** A typical antenna optimization design problem is presented, and various issues involved in the design process are discussed. Defining a suitable objective function is a central question, as is the type of optimization algorithm that should be used, stochastic versus deterministic. These questions are addressed by way of an example. A single-resistor loaded broadband HF monopole design is considered in detail, and the resulting antenna compared to published results for similar continuously loaded and discrete resistor loaded designs.

***Ver.2 April 8, 2011*** [Formatting of Table I corrected. References updated. Source code added in Appendix.]
***Brewster, Massachusetts***

# ISSUES IN ANTENNA OPTIMIZATION - A MONOPOLE CASE STUDY

**R. A. Formato**

Registered Patent Attorney & Consulting Engineer
Cataldo & Fisher, LLC, 400 TradeCenter, Suite 5900
Woburn, MA 01801 USA
rf2@ieee.org, rformato@cataldofisher.com

**Abstract —** A typical antenna optimization design problem is presented, and various issues involved in the design process are discussed. Defining a suitable objective function is a central question, as is the type of optimization algorithm that should be used, stochastic versus deterministic. These questions are addressed by way of an example. A single-resistor loaded broadband HF monopole design is considered in detail, and the resulting antenna compared to published results for similar continuously loaded and discrete resistor loaded designs.

**Keywords —** Antenna optimization, optimization, optimization algorithm, wire antenna, monopole, HF, impedance loading, impedance profile, resistive loading, broadband, bandwidth, numerical methods, CFO, Central Force Optimization, metaheuristic

## 1. INTRODUCTION

"*Good against remotes is one thing. Good against the living, that's something else.*" Han Solo thus cautioned Luke Skywalker as he practiced lightsaber skills in the classic 1977 film **Star Wars**. This note echoes a similar sentiment when it comes to designing antennas with optimization algorithms, "*Good against benchmarks is one thing. Good against 'real world' antennas, that's something else.*" This admonition is examined by way of an example, designing an optimized resistively-loaded broadband high-frequency (HF) monopole antenna. It highlights the importance of, and the difficulties in, choosing an appropriate objective function and the advantages of using a deterministic optimizer in doing so.

Optimization algorithms typically are evaluated against benchmarks with known extrema (fitnesses and locations). How well an algorithm works is measured by its accuracy and efficiency, referring respectively to how close it gets to the extrema and how much computational effort is expended in the process (usually the number of function evaluations).

An algorithm's performance often depends on user-specified setup parameters, and it may change dramatically with different values. Additional complications are introduced by inherently stochastic optimizers, such as Particle Swarm or Ant Colony Optimization, because this type of algorithm returns a different answer on successive runs, relying as they do on true random variables computed from a probability distribution. Even before an antenna problem has been precisely stated, the designer must choose suitable run parameters for a stochastic optimizer and somehow guess how well it will work on the antenna problem, neither of which is a simple matter.

The picture is further complicated because real world antenna problems introduce yet another level of complexity, defining a suitable objective function. If an optimization algorithm's performance is sensitive to the setup parameters, and its results vary from one run to the next, then the added problem of having to define a "good" objective function can be daunting. Of course, this question does not come up in benchmark testing because the benchmark *is* the objective function. But, as the results reported here show, this question is central in optimizing even a simple antenna.

## 2. DESIGN GOALS

The first step in antenna optimization is defining a clear set of performance goals. There are many measures to consider, such as directivity, radiation pattern, bandwidth, efficiency, and physical size, among others. Goals for all parameters must be articulated in order to define an objective function that effectively measures how well they are met. There are two main objectives for the monopole example described here: (1) as flat as possible an impedance bandwidth from 5 to 30 MHz, and (2) maximum gain. The metallic monopole element is 10.7 meters high with 0.005 meter radius (dimensions chosen for comparison to other designs) fed against a perfect electrically conductive (PEC) ground plane.

## 3. IMPORTANCE OF THE OBJECTIVE FUNCTION

The next step is defining an objective function that measures how well a particular antenna design meets the performance goals. For purposes of illustration, the monopole's decision space, $\Pi$, is chosen to be two-dimensional (2-D) because this topology ("landscape") can be visualized. The decision variables are (1) the value of the loading resistor, $R$ ($\Omega$), and (2) its placement along the monopole, $H$ (m), as shown in Figure 1. $\Pi : \left\{ (R, H) \,\middle|\, 0 \leq R \leq 1000\,\Omega \,,\, 0.05 \leq H \leq 10.65\,\text{m} \right\}$ therefore defines the decision space for this example. Of course, real world antenna problems usually contain many more than two variables, often far more, which considerably complicates the definition of a good objective function because then the landscape cannot be visualized.

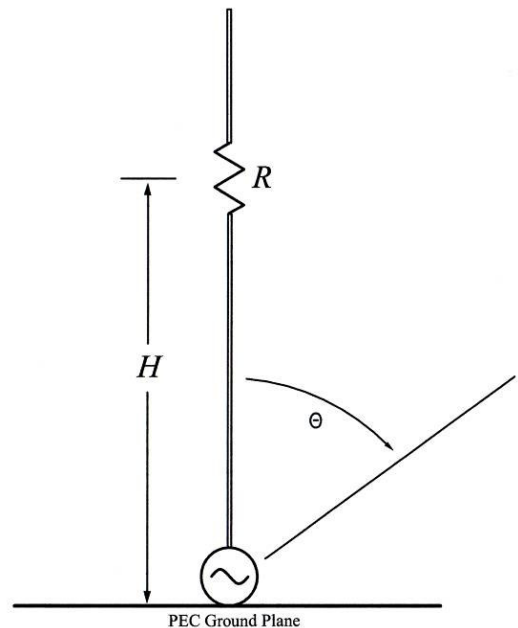

**Figure 1.** Base-Fed Monopole

The antenna parameters considered for inclusion in the objective function in this example are: minimum radiation efficiency, $Min(\varepsilon)$; minimum maximum gain, $Min(G_{\max})$; voltage standing wave ratio $VSWR // Z_0$ computed relative to a purely resistive feed system characteristic impedance $Z_0$; and monopole input impedance $Z_{in} = R_{in} + jX_{in}$. Each parameter is evaluated as a function of frequency, and $G_{\max}$ at a given frequency is the maximum gain over the polar angle $\Theta$ in Figure 1.

Needless to say, there are limitless ways these parameters can be combined, and the question is, Which combination is best? Unfortunately there is no answer to this query, other than trying different possibilities and evaluating each one's performance. A deterministic optimizer can make a big difference in this regard. Because stochastic algorithms return different results for every run, there is no (good) way of determining whether or not better designs are the result of a more suitable objective function or the inherent variability of the algorithm itself. By contrast, a deterministic optimizer always returns the same results for given setup parameters, so that any improvement in antenna performance is attributable to changes in the objective function. These considerations are illustrated below using three different objective functions for the loaded monopole.

## 4. FUNCTION $f_1(R, H, Z_0)$

The first monopole objective function (to be maximized) combines the minimum radiation efficiency and maximum gain with the maximum VSWR excursion in a simple formula, $f_1(R, H, Z_0) = \dfrac{Min(\varepsilon) + Min(G_{max})}{\Delta VSWR(Z_0)}$, where $\Delta VSWR(Z_0)$ is the difference between maximum and minimum standing wave ratios over the 5-30 MHz HF band relative to $Z_0$. The fitness increases with increasing efficiency and minimum gain and decreasing VSWR difference. $f_1$'s landscape with $Z_0 = 50\,\Omega$ appears in Figures 2(a)-(d), respectively, which show perspective and plan views and projections onto the $R - Z$ and $H - Z$ planes. $f_1$ is smoothly varying and unimodal with a maximum fitness of 2.3764… at the point $(R, H) = (5.025126\,\Omega, 1.621357\,\text{m})$. The global maximum was located by computing $f_1(R, H, 50)$ over the decision space $\Pi$ using a grid of 200x200 points and searching for the maximum. This procedure is used for each of the 2-D landscapes discussed here.

The monopole was modeled using Numerical Electromagnetics Code, version 2, double precision (NEC-2D), which is freely available online [1]. The 10.7 meter tall element was divided into 107 segments with the resistor placed at the segment's midpoint using "LD" loading cards (see [2,3] for NEC modeling guidelines and data input formats). Because the antenna is loaded by segment number, not distance above the ground plane, the height coordinate $H$ was converted to the loading segment number as $n = \left\lfloor 0.5 + \dfrac{H}{\Delta} \right\rfloor$ where $\Delta = 0.1\,\text{m}$ is the segment length. A typical NEC input file appears in Figure 3.

While $f_1$'s functional form may seem quite reasonable for measuring the monopole's performance, an examination of its topology reveals two potential concerns: maximum fitness occurs close to $\Pi$'s lower resistance boundary, and it varies very little with height. The first characteristic may impede an optimization algorithm's ability to search $\Pi$ while the second may impede convergence (exploration vs. exploitation). In a higher dimensionality decision space these characteristics cannot be ascertained by inspection, which is a further complication in defining a useful 'real world' objective function.

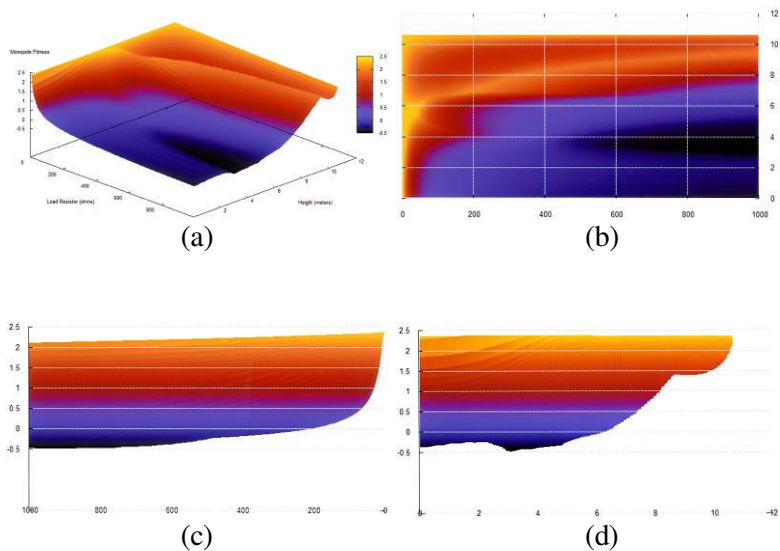

**Figure 2.** Landscape for Objective Function $f_1(R, H, 50)$

```
CM File: DES1.NEC
CM NEC2D run using R,Z values from DESIGN #1 DS plot
CM R=5.025126 ohms, Z=1.621357 m
CM seg # = INT(0.5+Z/SegLen) = 16
CM Zo=50 ohms
CE
GW1,107,0.,0.,0.,0.,0.,10.7,.005
GE1
LD0,1,16,16,5.025126,0.,0.
GN1
FR 0,26,0,0,5.,1.
EX 0,1,1,1,1.,0.
RP 0,10,1,1001,0.,0.,10.,0.,100000.
EN
```

**Figure 3.** Typical Monopole NEC Input File.

Apart from the question of how "searchable" $\Pi$ is for $f_1$'s maxima, perhaps the more important question is how well $f_1$ actually reflects a good monopole design, that is, one that performs well against the stated performance objectives. In this example, because $f_1$'s maximum can be visualized and located, the resulting "best" monopole design can be evaluated by computing its performance using the known maximum's coordinates. A feed system characteristic impedance of $Z_0 = 50\,\Omega$ is assumed because typical HF transmitters are designed for $50\,\Omega$ systems, and the results appear in Figure 4. In the plots, calculated data points are shown as symbols (5-30 MHz every 1 MHz), and the solid curves are interpolated using a natural cubic spline. Total power gain was computed every $10°$ for $0° \le \Theta \le 90°$ where $\Theta$ is the polar angle in NEC's standard right-handed spherical polar coordinate system (see Figure 1).

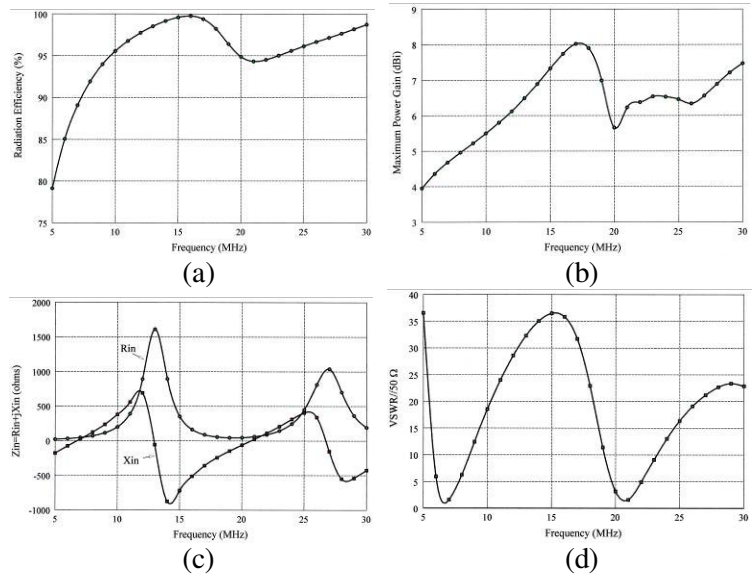

**Figure 4.** Performance of Monopole Design #1, $f_1(5.025, 1.62, 50)$.

The radiation efficiency and maximum gain numbers are quite good. Minimum efficiency is just below 80%, and for the most part the efficiency exceeds 95% above 10 MHz. Of course, this result is not altogether unexpected in view of the very light loading, $R \approx 5\,\Omega$. The maximum power gain figures also are quite good, with minimum $G_{\max} \approx 4\,\text{dBi}$. This result also is expected in view of the light loading.

But in stark contrast, the VSWR performance is very poor. The goal of flattening VSWR as much as possible was missed completely. VSWR varies from 1.61 to nearly 37 with pronounced fluctuations. The impedance bandwidth of this design, typically specified as $\text{VSWR} \le 2:1$ (return loss $\le -10\,\text{dB}$), is extremely small. Perhaps somewhat surprisingly, even though $f_1(R, H, 50)$'s functional form appeared to be a reasonable measure of how well the loaded monopole meets the design goals, the fact is, it is not. The best that a perfectly accurate optimization algorithm could do is discover the design in Figures 3 and 4, and that design happens to be quite poor. This example shows how important it is to choose an appropriate objective function.

## 5. FUNCTION $f_2(R, H, Z_0)$

$f_1(R, H, Z_0)$'s disappointing results make it clear that another, hopefully better, objective function must be defined. $f_1$'s major failing was its inability to flatten the VSWR curve, which suggests that a more aggressive approach is required. Dealing directly with $Z_{in} = R_{in} + jX_{in}$, for example, might work better than trying to minimize VSWR variability. To that end the objective function $f_2(R, H, Z_0) = \dfrac{Min(\varepsilon)}{\left|Z_0 - Max(R_{in})\right| \cdot \left|Max(X_{in})\right|}$ will be considered next. As before, this functional form is simple and ostensibly serves to achieve the design goals. As will be seen below, it does perform better than $f_1(R, H)$ with respect to VSWR, but its topology is such that many optimization algorithms will have considerable difficulty locating maxima.

$f_2(R,H,Z_0)$'s global maximum value is 0.11117… at the point $(R,H) = (819.095477\,,\,2.953015)$. The loading resistance $R \approx 819\,\Omega$ is much heavier than before, which will reduce efficiency and maximum gain but hopefully will tame the VSWR. This reflects the inevitable trade-off in using impedance loading for improving antenna bandwidth, which increases with heavier loading at the expense of radiation efficiency and gain.

Figure 5 shows $f_2$'s landscape. It comprises a series of spikes along a bullet-shaped curve in the $(R,H)$-plane, and the peaks are quite sharp. For example, changing $H$ slightly from 2.8731… to 2.8198… with $R = 819.0954…$ results in nearly three orders of magnitude decrease in fitness. Topologies like this usually are described as "pathological" because many optimization algorithms have difficulty dealing with them. Thus, even though objective function $f_2$ may be better than $f_1$ for achieving the design goals, its pathological landscape may impede an optimization algorithm to such a degree that better designs are not discovered.

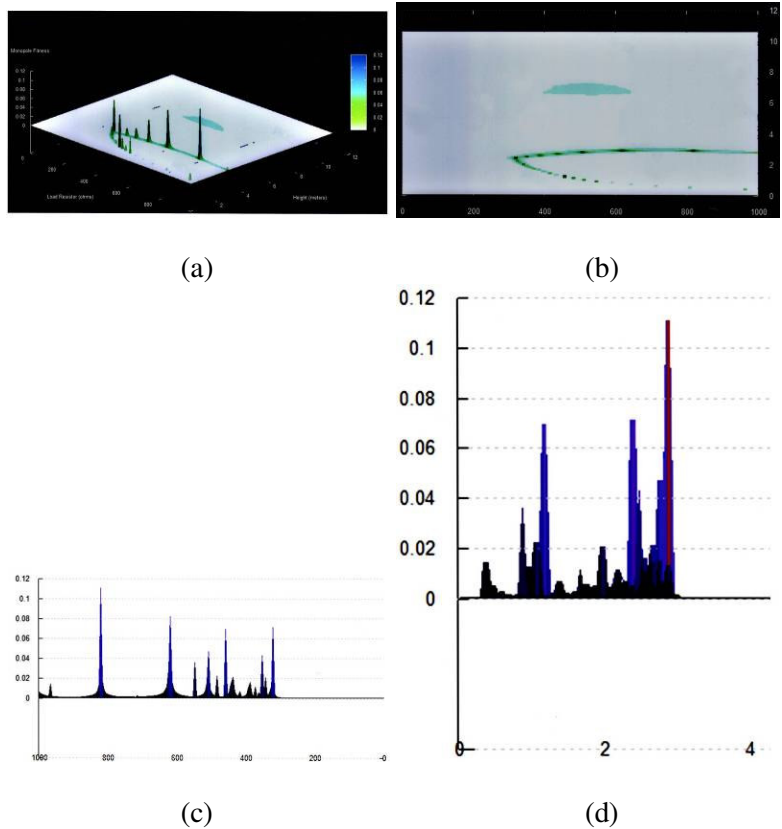

(a)    (b)

(c)    (d)

**Figure 5.**  Landscape for Objective Function $f_2(R,H,50)$

As in §4, NEC-2D was used to model the monopole with $f_2(R,H,50)$'s best fitness, and the results appear in Figure 6. As expected, the radiation efficiency is much lower, especially at low frequencies. Below 15 MHz it ranges from about 5% to 25%. The efficiency does increase substantially mid-band, reaching a peak near 80% at 19 MHz and falling thereafter. The power gain more or less tracks the efficiency, but it is quite low at low frequencies. Above 15 MHz, however, the gain is moderate to good. The heavier loading in this case

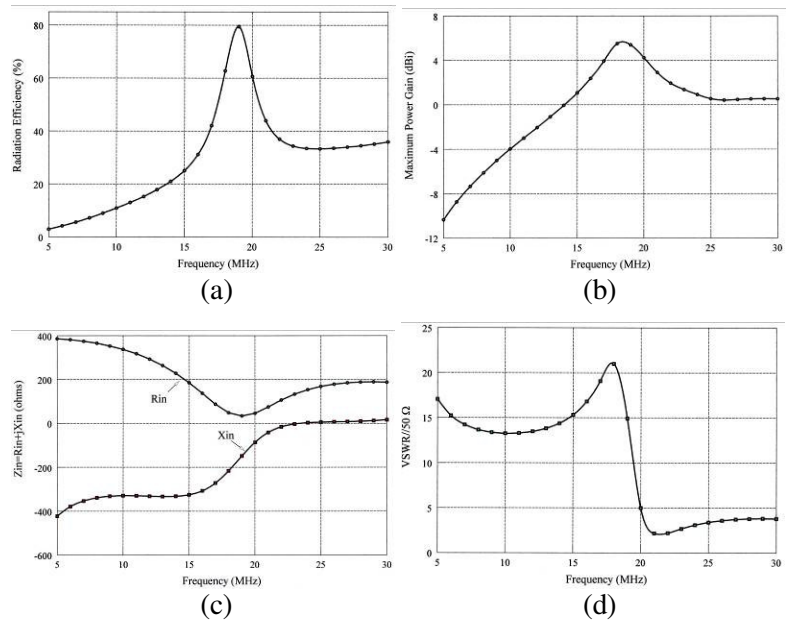

**Figure 6.** Performance of Monopole Design #2, $f_2(819.1, 2.95, 50)$.

considerably reduced the input impedance variation resulting in a fairly smooth variation in $R_{in}$ and, to a lesser degree, in $X_{in}$ as well. As a result VSWR variability is less than in the previous design, but still quite substantial. The VSWR is well-behaved and moderate, $\leq 5:1$, above 20

MHz, but it is very high at lower frequencies with a peak $\approx 21:1$ at 18 MHz. Thus, while $f_2$ is an improvement over $f_1$ in terms of meeting the design goals, it still falls far short of yielding a good monopole design. In addition, its pathological landscape may defeat the effectiveness of many optimization algorithms. Further refinement of the objective function is required.

## 6.  FUNCTION $f_3(R, H, Z_0)$

Because VSWR variability is the biggest problem with the first two objective functions, an even more aggressive approach will be taken with

$$f_3(R, H, Z_0) = \frac{Min(\varepsilon)}{\left|Z_0 - Max(R_{in})\right| \cdot \Delta VSWR(Z_0) \cdot [Max(Xin) - Min(Xin)]}.$$

The gain does not appear in the numerator because it tracks fairly well with efficiency. The denominator comprises three factors that minimize VSWR in different ways. The first drives the real part of the input impedance toward the feed system characteristic impedance. The second minimizes the VSWR variability across the band, while the third attempts to flatten the input reactance.

Because this functional form is determined empirically, other forms probably merit consideration as well. For example, $f_3(R, H, Z_0)$ could be written

$$f_3(R, H, Z_0) = \frac{Min(\varepsilon)^{\eta_1}}{\left|Z_0 - Max(R_{in})\right|^{\eta_2} \cdot \Delta VSWR(Z_0)^{\eta_3} \cdot [Max(Xin) - Min(Xin)]^{\eta_4}}$$

where the exponents $\eta_i$ are constants or functions of frequency. The terms in $f_3$ could be combined differently, say, by addition with weighting coefficients. Other functions, such as logarithms or trigonometric functions, might be useful in combining the antenna's performance measures. And, of course, other performance measures might be included as well. All of these considerations are involved in defining suitable objective functions. As the results for $f_1$ and $f_2$ show, presumably good ones can turn out to be quite poor. It consequently is imperative to investigate how the objective function's form influences the resulting antenna design.

$f_3(R,H,50)$ is unimodal with a smoothly varying topology and a maximum fitness of $1.4624...\text{x}10^{-5}$ at $(R,H) = (502.512563, 7.2143215)$. Its landscape is plotted in Figure 7. This objective function results in the design whose performance is shown in Figure 8. The radiation efficiency increases more or less monotonically from just over 15% at 5 MHz to nearly 40% at 27 MHz and about 38% at 30 MHz. Maximum power gain ranges from a low near -3.1 dBi to a maximum of 2 dBi. The input impedance is well behaved across the HF band, and the resulting VSWR is much flatter than in the previous cases. Maximum VSWR is just below 13:1 at 5 MHz, and it falls very quickly to just above 3:1 at 7.5 MHz. The VSWR increases to ~8:1 at 13 MHz and remains fairly flat thereafter. A comparison of these data to the curves in Figures 4(d) and 6(d) clearly shows that $f_3(R, H, 50)$ is the best objective function of the three. Its monopole design is superior to the others, and its topology lends itself well to being searched by an optimization algorithm.

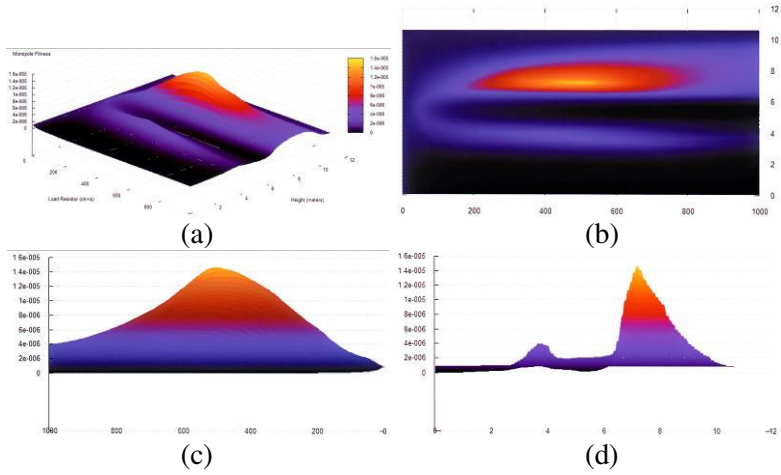

(a)                                                    (b)

(c)                                                    (d)

**Figure 7.**   Landscape for Objective Function $f_3(R,H,50)$

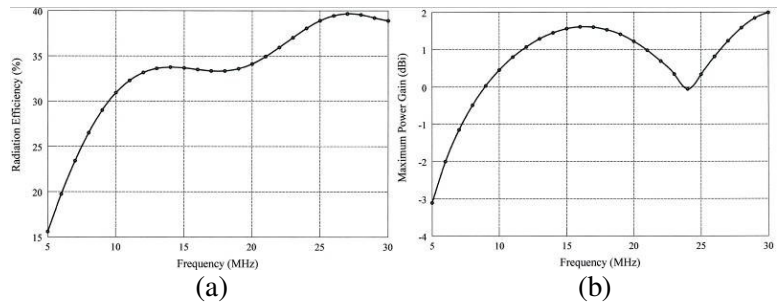

(a)                                                    (b)

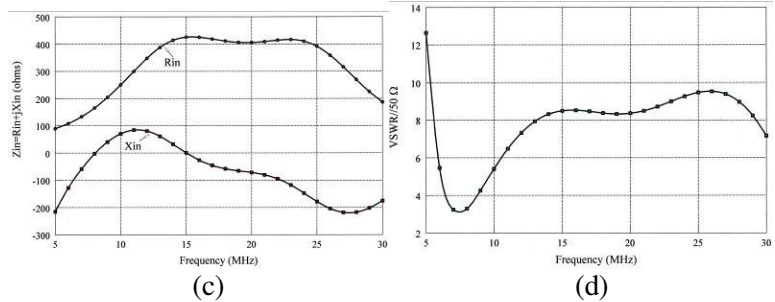

**Figure 8.** Performance of Monopole Design #3, $f_3(502.5, 7.214, 50)$.

## 7. OPTIMIZATION ALGORITHMS

The previous sections discussed some of the issues in defining suitable objective functions for the broadband HF monopole design problem. Three functions were considered, and the results varied considerably from one to the next, with the last one being the best. Each of these objective functions has a known global maximum that can be visualized because the monopole decision space is 2-D. Unfortunately this cannot be done in higher dimensionality spaces, so that their topologies are unknown. The problem faced by the antenna designer therefore is defining an effective objective function that can be searched accurately and efficiently in the *n*-D decision space. The type of optimization algorithm can be an important factor in aiding or inhibiting the process of defining a suitable objective function.

Because stochastic optimization algorithms return different results on successive runs, it is difficult to assess the effects of changing the objective function on their accuracy and efficiency. For example, if Particle Swarm Optimization is applied to the monopole problem, the antenna designer cannot know why successive runs using, say, $f_2(R,H,Z_0)$ and $f_3(R,H,Z_0)$, yield different results. It may be a consequence of the different objective functions (for example, pathological vs. well-behaved), or it may be the algorithm's inherent randomness. Which of these it is can be ascertained only by doing a statistical analysis that probably requires tens or hundreds, possibly thousands, of runs. This dilemma is avoided by using a *deterministic optimizer*, one that yields the same answer for every run with the same setup.

Central Force Optimization [4-14] is a deterministic search and optimization metaheuristic that has performed well against various problems in applied electromagnetics [15-20]. CFO therefore was used to search the three decision spaces for each of the monopole objective functions. A "parameter free" implementation was employed as described in [4] (without directional information in errant probe repositioning, [12]). CFO pseudocode appears in Figure 9. Hardwired parameter values were $N_d = 2$, $F_{rep}^{init} = 0.5$, $\Delta F_{rep} = 0.1$, $F_{rep}^{min} = 0.05$, $\gamma_{start} = 0$, $\gamma_{stop} = 1$, $\Delta\gamma = 0.3333$, $(N_p/N_d)_{max} = 6$, $N_t = 200$ with an early termination criterion of fitness variation $\leq 10^{-6}$ for 25 consecutive steps starting at step #35.

Table 1 summarizes the CFO results. It shows that the algorithm performed well against objective functions $f_1$ and $f_3$ by discovering maxima close to the known values. The results for $f_1$ are consistent with its topology in which the global maximum is near $\Pi$'s lower boundary in $R$ and not particularly sensitive to variations in $H$. CFO essentially recovered $f_3$'s global maximum, but clearly it had a problem with $f_2$'s pathological landscape. NEC-2D's computed performance for the antenna design using CFO's $(R,H)$ coordinates for the best objective function, $f_3$, appears in Figure 10, which, as expected, is essentially the same as Figure 8.

CFO's performance probably would be better still if a more stringent early termination criterion were employed, or if none were used in a much longer run. The purpose of this note, however, is to discuss real-world design issues, and one of them is having to make the engineering decision of when the design is "good enough" relative to the resources expended. In this case, CFO achieved an acceptable design very close to the known best design using $f_3$ and a total of 4,636 function evaluations. This meets the "good enough" test.

This simple monopole example demonstrates that how well an "optimized" antenna performs can be highly dependent upon both the objective function against which it is optimized and how accurately and efficiently the optimization program performs against that function's landscape. Because CFO is deterministic, it allows the antenna designer to investigate the effects of changing the objective function's form and parameters that determine its landscape. It is evident that defining an effective function is much easier when the optimizer returns the same results every time instead of different ones. In the author's opinion this is an important consideration in addressing real-world antenna problems, or, for that matter, any problem in which definition of the objective function is an issue.

**Table 1.** CFO Optimization Results

| Objective Function | Known Max Fitness / Coords | CFO Max Fitness / Coords |
|---|---|---|
| $f_1(R,H,50)$ | 2.376 / (5.025,1.621) | 2.371 / (8.179,5.361) |
| $f_2(R,H,50)$ | 0.1112 / (819.1,2.953) | 0.0684 / (322.2,2.366) |
| $f_3(R,H,50)$ | $1.462 \times 10^{-5}$ / (502.5,7.241) | $1.401 \times 10^{-5}$ / (499.6,7.302) |

```
Procedure  CFO [ f(x⃗), N_d, Π ]
```

Internals:  $N_t$, $F_{rep}^{init}$, $\Delta F_{rep}$, $F_{rep}^{min}$, $\left(\dfrac{N_p}{N_d}\right)_{MAX}$, $\gamma_{start}$, $\gamma_{stop}$, $\Delta\gamma$.

Initialize $f_{max}^{global}(\vec{x})$ = very large negative number, say, $-10^{+4200}$.

For $N_p/N_d = 2$ to $\left(\dfrac{N_p}{N_d}\right)_{MAX}$ by $2$:

(a.0)     Total number of probes:  $N_p = N_d \cdot \left(\dfrac{N_p}{N_d}\right)$

For $\gamma = \gamma_{start}$ to $\gamma_{stop}$ by $\Delta\gamma$:

(a.1)  Re-initialize data structures for position/ acceleration vectors & fitness matrix.

(a.2)  Compute IPD (see [4]).

(a.3)  Compute initial fitness matrix, $M_0^p, 1 \le p \le N_p$.

(a.4)  Initialize $F_{rep} = F_{rep}^{init}$.

For $j = 0$ to $N_t$ (or earlier termination – see [4]):

(b)     Compute position vectors, $\vec{R}_j^p, 1 \le p \le N_p$ (eq.(2) in [4])

(c)     Retrieve errant probes $(1 \le p \le N_p)$:

If $\vec{R}_j^p \cdot \hat{e}_i < x_i^{min} \therefore$  $\vec{R}_j^p \cdot \hat{e}_i = \max\{x_i^{min} + F_{rep}(\vec{R}_{j-1}^p \cdot \hat{e}_i - x_i^{min}), x_i^{min}\}$.

If $\vec{R}_j^p \cdot \hat{e}_i > x_i^{max} \therefore$  $\vec{R}_j^p \cdot \hat{e}_i = \min\{x_i^{max} - F_{rep}(x_i^{max} - \vec{R}_{j-1}^p \cdot \hat{e}_i), x_i^{max}\}$.

(d)     Compute fitness matrix for current probe distribution, $M_j^p, 1 \le p \le N_p$.

(e)     Compute accelerations using current probe distribution and fitnesses (eq.(1) in [4]).

(f)     Increment $F_{rep}$: $F_{rep} = F_{rep} + \Delta F_{rep}$; If $F_{rep} > 1 \therefore F_{rep} = F_{rep}^{min}$.

(g)     If $j \ge 20$ and $j\,MOD\,10 = 0 \therefore$

(i)  Shrink $\Omega$ around $\vec{R}_{best}$ (see [4].

(ii)  Retrieve errant probes [procedure Step (c)].

Next $j$

(h)     Reset $\Omega$ boundaries [values before shrinking].

(i)     If $f_{max}(\vec{x}) \ge f_{max}^{global}(\vec{x}) \therefore$  $f_{max}^{global}(\vec{x}) = f_{max}(\vec{x})$.

Next $\gamma$

Next $N_p/N_d$

**Figure 9.** CFO Pseudocode

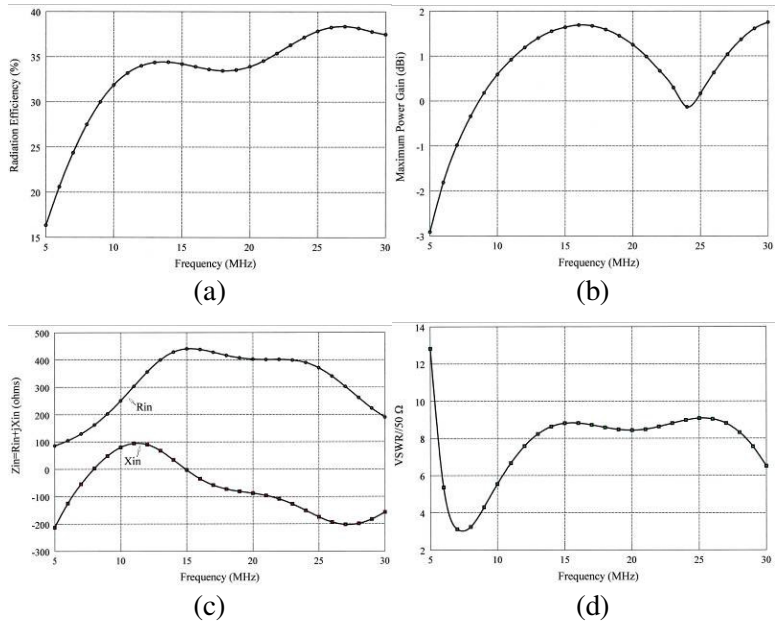

(a)

(b)

(c)

(d)

**Figure 10.** CFO-Optimized Monopole, $f_3(499.6, 7.302, 50)$.

## 8. OTHER LOADED MONOPOLE DESIGNS

The example antenna structure in this paper was inspired by the HF monopole designed and tested by Rama Rao and Debroux [21,22]. Their antenna employed an analytically computed continuous resistive loading profile that achieved VSWR $\leq 2:1$ from 5-30 MHz with the use of a matching network. Radiation efficiency, gain and pattern were reasonable and well-behaved for supporting moderate to long range HF links. A CFO-designed discrete-resistance loading profile for the same antenna geometry provided even better performance [23], and it is instructive to compare that CFO design, which utilized fourteen discrete resistors, to the $f_3(499.6, 7.302, 50)$ design above that uses only a single resistor.

Figure 11 shows the NEC input file for the 14-segment, CFO-optimized monopole. Note that this antenna was optimized with $Z_0 = 300\,\Omega$ instead of $Z_0 = 50\,\Omega$ because the $300\,\Omega$ reference was used in [21,22]. Each monopole segment is loaded at its center with a discrete resistor whose value ranges from about $7.2\,\Omega$ to $82.7\,\Omega$ (NEC "LD" cards). This antenna's performance is shown in Figure 12. Radiation efficiency ranges from a low of about 8% to a high of 45%, while the maximum gain increases from ~-5.5 dBi at 5 MHz to ~3 dBi at 28 MHz with a pronounced mid-band dip. The most important performance measure, VSWR//300Ω, is quite good, being less than 2.25:1 at all frequencies, and less than 2:1 above about 5.5 MHz. This antenna meets the VSWR $\leq 2:1$ goal essentially across the entire HF band with no matching network. The only additional element needed to feed this antenna from a $50\,\Omega$ system is a low-loss, broadband 6:1 unun, which is readily available.

The results for the 14-segment monopole and the VSWR data in Figure 10(d) suggest that a better feed system impedance for the $f_3(499.6, 7.302, 50)$ design actually might be $300\,\Omega$. This conjecture was investigated by recalculating this design's VSWR parametrically in $Z_0$ from $225\,\Omega$ to $350\,\Omega$. Figure 13 plots the results. The best overall performance indeed does occur with $Z_0 \approx 300\,\Omega$. It is $\leq 2:1$ from ~7.5 MHz through ~28 MHz. VSWR above 28 MHz remains fairly low, below 2.5:1; but at the low end of the band it does increase quickly as the frequency drops. VSWR maximum is ~5.5:1 at 5 MHz, but even this value is quite acceptable because it is high only in a fairly narrow band, and values less than 10:1 are readily matched with simple networks.

```
CM File: LD_MONO.NEC
CM Run ID 02-16-2011 10:43:46
CM Nd= 14, p= 2, j= 120
CM Zo=300 ohms
CE
GW1,14,0.,0.,0.,0.,0.,10.668,.0254
GE1
LD0,1,1,1,82.7045,0.,0.
LD0,1,2,2,29.31145,0.,0.
LD0,1,3,3,9.2825,0.,0.
LD0,1,4,4,7.154042,0.,0.
LD0,1,5,5,7.397769,0.,0.
LD0,1,6,6,7.310225,0.,0.
LD0,1,7,7,27.58697,0.,0.
LD0,1,8,8,26.55749,0.,0.
LD0,1,9,9,24.70102,0.,0.
LD0,1,10,10,22.80148,0.,0.
LD0,1,11,11,20.82445,0.,0.
LD0,1,12,12,16.44918,0.,0.
LD0,1,13,13,11.4537,0.,0.
LD0,1,14,14,9.471994,0.,0.
GN1
FR 0,26,0,0,5.,.1.
EX 0,1,1,1,1.,0.
RP 0,10,1,1001,0.,0.,10.,0.,100000.
EN
```

**Figure 11.** NEC Input File CFO-Optimized 14-Seg Monopole.

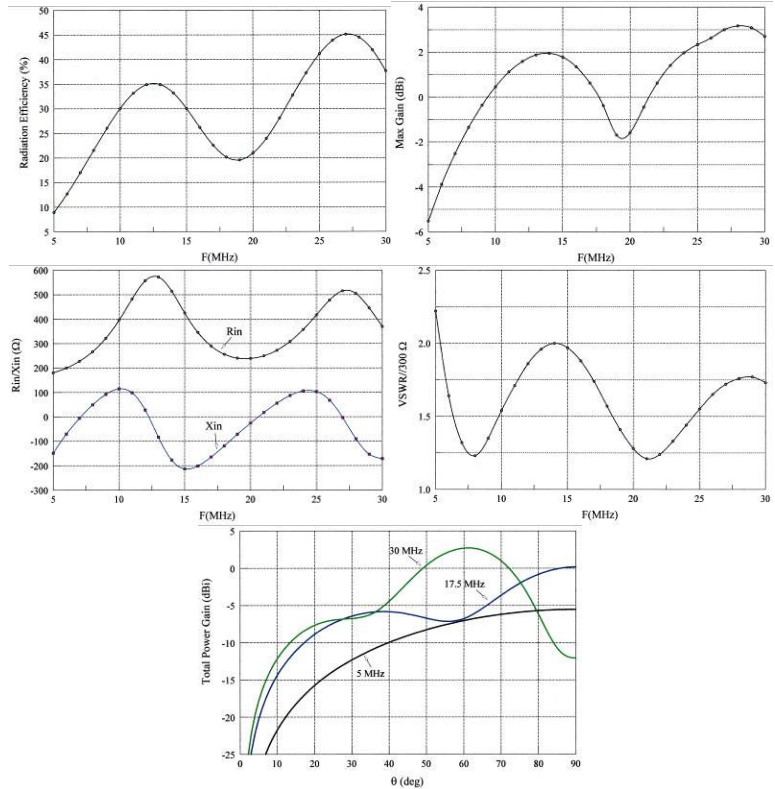

**Figure 12.** Performance, CFO-optimized 14-Seg Monopole, $Z_0 = 300\,\Omega$ .

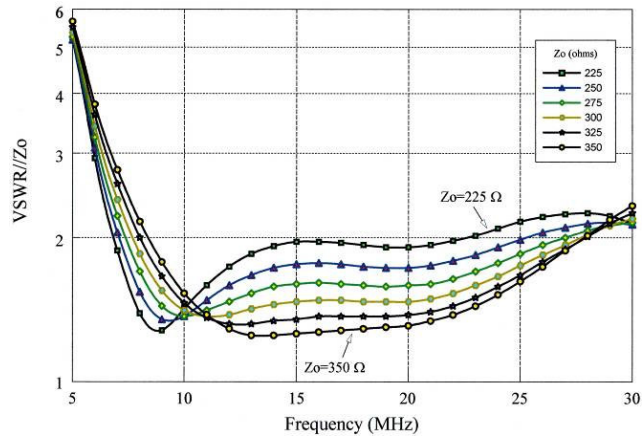

**Figure 13.** VSWR for CFO Design $f_3(499.6, 7.302, 50)$ Parametric in $Z_0$ .

At this point in the design process, it is apparent that a metallic monopole loaded by a single correctly placed resistor may perform nearly as well as one employing fourteen resistors, and probably better than the designs in [21,22] employing continuous loading. Because the objective function's landscape changes with $Z_0$, even if only slightly, and because the CFO-optimized antenna with $Z_0 = 50\,\Omega$ exhibits better VSWR when $Z_0$ is increased to $300\,\Omega$, it is instructive to tweak the previous design by making another CFO run against $f_3$'s landscape with $Z_0 = 300\,\Omega$, that is, with $f_3(R,H,300)$ as the objective function. CFO's best fitness in this case is $9.66379 \times 10^{-5}$ at $(R,H) = (501.78982, 7.107012)$ using 5,016 function evaluations. The tweaked design's performance appears in Figure 14.

The single-resistor monopole outperforms the 14-segment monopole by every measure except VSWR. The minimum radiation efficiency is 15%, and it increases nearly monotonically. By contrast, radiation efficiency for the 14-segment antenna starts off near 8% and exhibits considerable fluctuation with increasing frequency. The tweaked design's gain increases from ~ -3.25 dBi at 5 MHz to ~2.25 dBi at 30 MHz with a dip to 0 dBi at 24 MHz, whereas the 14-segment design has a similar behavior but a lower gain at 5 MHz (~ -5.6 dBi) and a very slightly higher gain at 28 MHz (~3.1 dBi). The VSWR performance of the tweaked antenna is somewhat worse, but nonetheless quite good. It is below 2:1 from 7.5-26 MHz and only slightly above that through 30 MHz where it reaches 2.5:1. Below 7.5 MHz VSWR increases quickly with decreasing frequency reaching just over 5:1 at 5 MHz. But this degree of variability can easily be handled by a simple matching network. It is reasonable to expect a VSWR below 2:1 across the entire 5-30 MHz band, possibly well below 2:1.

Perhaps the most important measure of the tweaked monopole's effectiveness is its radiation pattern compared to the 14-segment design. The tweaked design generally exhibits higher power gain at all polar angles, especially in the angle range of interest for moderate to long range HF links, $60° \leq \Theta \leq 80°$. Comparing the 14-segment loaded monopole to the single-resistor $f_3(R,H,300)$ design, it is clear that the single-resistor loading provides better overall performance. In addition, this monopole is much simpler to fabricate and maintain, and it arguably is substantially better than the continuously loaded designs in [21,22].

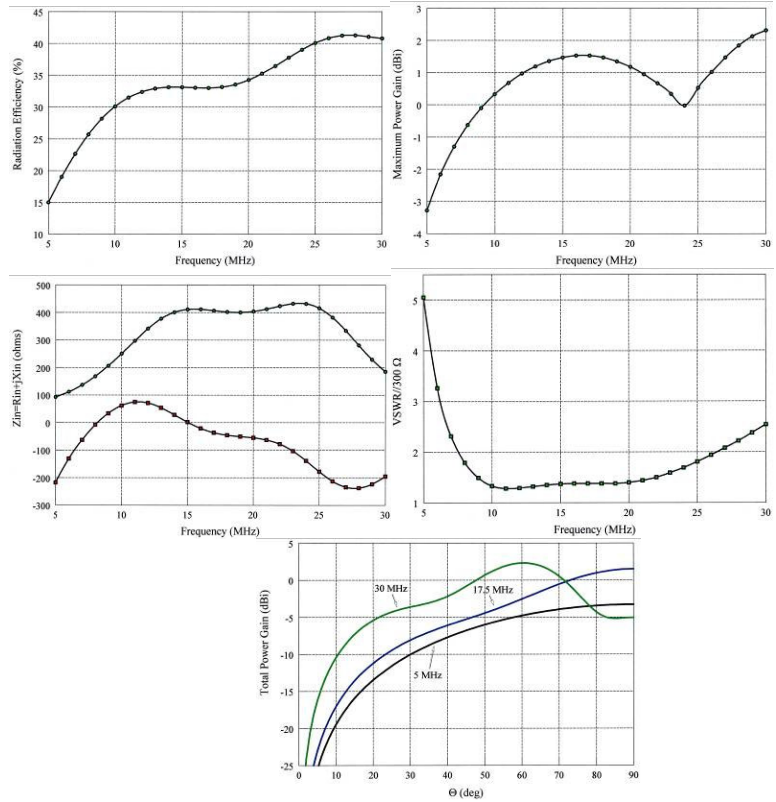

**Figure 14.** Performance of Tweaked Single-Resistor CFO-optimized Monopole, $f_3(501.8, 7.107, 300)$, $Z_0 = 300 \, \Omega$.

## 9. CONCLUSION

This paper discussed an example antenna optimization problem using a single-resistor loaded HF monopole. It addressed issues in defining an objective function that effectively measures antenna performance, and the suitability of stochastic and deterministic optimization algorithms. Assessing how well a particular objective function will achieve design goals is a difficult question because the landscape of functions beyond 2-D cannot be visualized (note that CFO's tendency to distribute probes may be useful in this regard; see §9 in [16]).

Three different objective functions and their landscapes were considered. The deterministic Central Force Optimization algorithm was applied each objective function's topology and the results compared. The best objective function then was used to develop a final tweaked design that was compared to a similar, previously CFO-optimized design utilizing fourteen discrete loading resistors. The single-resistor monopole performs as well or better than the fourteen resistor version and better than other continuously loaded antennas reported in the literature. CFO source code appears in the Appendix, and any papers not freely accessible online are available upon request to the author (rf2@ieee.org).

***Ver.2 April 8, 2011*** [Formatting of Table I corrected. References updated. Source code added in Appendix.]

***Brewster, Massachusetts***

# APPENDIX – CFO SOURCE CODE

```
'Program 'CFO_LOADED_MONO_02-11-2011.BAS' compiled with
'Power Basic/Windows Compiler 9.03 (www.PowerBasic.com).

'PROGRAM USED TO GENERATE DATA FOR SINGLE RESISTOR LOADED MONOPOLE PAPER

'LAST MOD 03-04-2011 ~1641 HRS EST

'(c) 2006-2011 Richard A. Formato

'ALL RIGHTS RESERVED WORLDWIDE

'THIS PROGRAM IS FREEWARE.  IT MAY BE COPIED AND
'DISTRIBUTED WITHOUT LIMITATION AS LONG AS THIS
'COPYRIGHT NOTICE AND THE GNUPLOT AND REFERENCE
'INFORMATION BELOW ARE INCLUDED WITHOUT MODIFICATION,
'AND AS LONG AS NO FEE OR COMPENSATION IS CHARGED,
'INCLUDING "TIE-IN" OR "BUNDLING" FEES CHARGED FOR
'OTHER PRODUCTS.

'===========================================================

'THIS PROGRAM REQUIRES wgnuplot.exe TO DISPLAY PLOTS.
'Gnuplot is a copyrighted freeware plotting program
'available at http://www.gnuplot.info/index.html.

'IT ALSO REQUIRES A VERSION OF THE Numerical Electromagnetics
'Code (NEC) in order to run the PBM benchmarks.  If this file
'is not present, a runtime error occurs.  Remove the code
'that checks for the NEC EXE is there is no interest in the
'PBM functions.

'===========================================================

#COMPILE EXE

#DIM ALL

%USEMACROS = 1

#INCLUDE "Win32API.inc"

DEFEXT A-Z

'------ EQUATES -----

%IDC_FRAME1       = 101
%IDC_FRAME2       = 102

%IDC_Function_Number1  = 121
%IDC_Function_Number2  = 122
%IDC_Function_Number3  = 123
%IDC_Function_Number4  = 124
%IDC_Function_Number5  = 125
%IDC_Function_Number6  = 126
%IDC_Function_Number7  = 127
%IDC_Function_Number8  = 128
%IDC_Function_Number9  = 129
%IDC_Function_Number10 = 130
%IDC_Function_Number11 = 131
%IDC_Function_Number12 = 132
%IDC_Function_Number13 = 133
%IDC_Function_Number14 = 134
%IDC_Function_Number15 = 135
%IDC_Function_Number16 = 136
%IDC_Function_Number17 = 137
%IDC_Function_Number18 = 138
%IDC_Function_Number19 = 139
%IDC_Function_Number20 = 140
%IDC_Function_Number21 = 141
%IDC_Function_Number22 = 142
%IDC_Function_Number23 = 143
%IDC_Function_Number24 = 144
%IDC_Function_Number25 = 145
%IDC_Function_Number26 = 146
%IDC_Function_Number27 = 147
%IDC_Function_Number28 = 148
%IDC_Function_Number29 = 149
%IDC_Function_Number30 = 150
%IDC_Function_Number31 = 151
%IDC_Function_Number32 = 152
%IDC_Function_Number33 = 153
%IDC_Function_Number34 = 154
%IDC_Function_Number35 = 155
%IDC_Function_Number36 = 156
%IDC_Function_Number37 = 157
%IDC_Function_Number38 = 158
%IDC_Function_Number39 = 159
%IDC_Function_Number40 = 160
%IDC_Function_Number41 = 161
%IDC_Function_Number42 = 162
%IDC_Function_Number43 = 163
%IDC_Function_Number44 = 164
%IDC_Function_Number45 = 165
%IDC_Function_Number46 = 166
%IDC_Function_Number47 = 167
%IDC_Function_Number48 = 168
%IDC_Function_Number49 = 169
%IDC_Function_Number50 = 170

'------------------------------ GLOBAL CONSTANTS & SYMBOLS ---------------------------

GLOBAL MonopoleSecsPerRun AS EXT

GLOBAL XiOffset() AS EXT 'offset array for Rosenbrock F6 function

GLOBAL XiMin(), XiMax(), DiagLength, StartingXiMin(), StartingXiMax() AS EXT 'decision space boundaries, length of diagonal

GLOBAL Aij() AS EXT 'array for Shekel's Foxholes function

GLOBAL EulerConst, Pi, Pi2, Pi4, TwoPi, FourPi, e, Root2 AS EXT 'mathematical constants

GLOBAL Alphabet$, Digits$, RunID$  'upper/lower case alphabet, digits 0-9 & Run ID
```

```
GLOBAL Quote$, SpecialCharacters$   'quotation mark & special symbols

GLOBAL Mu0, Eps0, c, eta0 AS EXT   'E&M constants

GLOBAL Rad2Deg, Deg2Rad, Feet2Meters, Meters2Feet, Inches2Meters, Meters2Inches AS EXT 'conversion factors

GLOBAL Miles2Meters, Meters2Miles, NautMi2Meters, Meters2NautMi AS EXT          'conversion factors

GLOBAL ScreenWidth&, ScreenHeight& 'screen width & height

GLOBAL xOffset&, yOffset&           'offsets for probe plot windows

GLOBAL FunctionNumber%

GLOBAL AddNoiseToPBM2$

'------------------------------ TEST FUNCTION DECLARATIONS ----------------------------
DECLARE FUNCTION F1(R(),Nd%,p%,j&)              'F1 (n-D)

DECLARE FUNCTION F2(R(),Nd%,p%,j&)              'F2(n-D)

DECLARE FUNCTION F3(R(),Nd%,p%,j&)              'F3 (n-D)

DECLARE FUNCTION F4(R(),Nd%,p%,j&)              'F4 (n-D)

DECLARE FUNCTION F5(R(),Nd%,p%,j&)              'F5 (n-D)

DECLARE FUNCTION F6(R(),Nd%,p%,j&)              'F6 (n-D)

DECLARE FUNCTION F7(R(),Nd%,p%,j&)              'F7 (n-D)

DECLARE FUNCTION F8(R(),Nd%,p%,j&)              'F8 (n-D)

DECLARE FUNCTION F9(R(),Nd%,p%,j&)              'F9 (n-D)

DECLARE FUNCTION F10(R(),Nd%,p%,j&)             'F10 (n-D)

DECLARE FUNCTION F11(R(),Nd%,p%,j&)             'F11 (n-D)

DECLARE FUNCTION F12(R(),Nd%,p%,j&)             'F12 (n-D)

DECLARE FUNCTION u(xi,a,k,m)                    'Auxiliary function for F12 & F13

DECLARE FUNCTION F13(R(),Nd%,p%,j&)             'F13 (n-D)

DECLARE FUNCTION F14(R(),Nd%,p%,j&)             'F14 (n-D)

DECLARE FUNCTION F15(R(),Nd%,p%,j&)             'F15 (n-D)

DECLARE FUNCTION F16(R(),Nd%,p%,j&)             'F16 (n-D)

DECLARE FUNCTION F17(R(),Nd%,p%,j&)             'F17 (n-D)

DECLARE FUNCTION F18(R(),Nd%,p%,j&)             'F18 (n-D)

DECLARE FUNCTION F19(R(),Nd%,p%,j&)             'F19 (n-D)

DECLARE FUNCTION F20(R(),Nd%,p%,j&)             'F20 (n-D)

DECLARE FUNCTION F21(R(),Nd%,p%,j&)             'F21 (n-D)

DECLARE FUNCTION F22(R(),Nd%,p%,j&)             'F22 (n-D)

DECLARE FUNCTION F23(R(),Nd%,p%,j&)             'F23 (n-D)

DECLARE FUNCTION F24(R(),Nd%,p%,j&)             'F24 (n-D)

DECLARE FUNCTION F25(R(),Nd%,p%,j&)             'F25 (n-D)

DECLARE FUNCTION F26(R(),Nd%,p%,j&)             'F26 (n-D)

DECLARE FUNCTION F27(R(),Nd%,p%,j&)             'F27 (n-D)

DECLARE FUNCTION ParrottF4(R(),Nd%,p%,j&)       'Parrott F4 (1-D)

DECLARE FUNCTION SGO(R(),Nd%,p%,j&)             'SGO Function (2-D)

DECLARE FUNCTION GoldsteinPrice(R(),Nd%,p%,j&)  'Goldstein-Price Function (2-D)

DECLARE FUNCTION StepFunction(R(),Nd%,p%,j&)    'Step Function (n-D)

DECLARE FUNCTION Schwefel226(R(),Nd%,p%,j&)     'Schwefel Prob. 2.26 (n-D)

DECLARE FUNCTION Colville(R(),Nd%,p%,j&)        'Colville Function (4-D)

DECLARE FUNCTION Griewank(R(),Nd%,p%,j&)        'Griewank (n-D)

DECLARE FUNCTION Himmelblau(R(),Nd%,p%,j&)      'Himmelblau (2-D)

DECLARE FUNCTION Rosenbrock(R(),Nd%,p%,j&)      'Rosenbrock (n-D)

DECLARE FUNCTION Sphere(R(),Nd%,p%,j&)          'Sphere (n-D)

DECLARE FUNCTION HimmelblauNLO(R(),Nd%,p%,j&)   'Himmelblau NLO (5-D)

DECLARE FUNCTION Tripod(R(),Nd%,p%,j&)          'Tripod (2-D)

DECLARE FUNCTION Sign(X)                        'Auxiliary function for Tripod

DECLARE FUNCTION RosenbrockF6(R(),Nd%,p%,j&)    'Rosenbrock F6 (10-D)

DECLARE FUNCTION CompressionSpring(R(),Nd%,p%,j&)'Compression Spring (3-D)

DECLARE FUNCTION GearTrain(R(),Nd%,p%,j&)       'Gear Train (4-D)

DECLARE FUNCTION PBM_1(R(),Nd%,p%,j&)           'PBM Benchmark #1

DECLARE FUNCTION PBM_2(R(),Nd%,p%,j&)           'PBM Benchmark #2

DECLARE FUNCTION PBM_3(R(),Nd%,p%,j&)           'PBM Benchmark #3

DECLARE FUNCTION PBM_4(R(),Nd%,p%,j&)           'PBM Benchmark #4

DECLARE FUNCTION PBM_5(R(),Nd%,p%,j&)           'PBM Benchmark #5

DECLARE FUNCTION LOADED_MONO(R(),Nd%,p%,j&)     'Loaded Monopole
```

```
'------------------------------- SUB DECLARATIONS -------------------------------

DECLARE SUB CopyBestMatrices(Np%,Nd%,Nt&,R(),M(),Rbest(),Mbest())

DECLARE SUB CheckNECFiles(NECfileError$)

DECLARE SUB GetTestFunctionNumber(FunctionName$)

DECLARE SUB FillArrayAij

DECLARE SUB Plot3DbestProbeTrajectories(NumTrajectories%,M(),R(),Np%,Nd%,LastStep&,FunctionName$)

DECLARE SUB Plot2DbestProbeTrajectories(NumTrajectories%,M(),R(),Np%,Nd%,LastStep&,FunctionName$)

DECLARE SUB Plot2DindividualProbeTrajectories(NumTrajectories%,M(),R(),Np%,Nd%,LastStep&,FunctionName$)

DECLARE SUB Show2Dprobes(R(),Np%,Nt&,j&,Frep,BestFitness,BestProbeNumber%,BestTimeStep&,FunctionName$,RepositionFactor$,Gamma)

DECLARE SUB Show3Dprobes(R(),Np%,Nd%,Nt&,j&,Frep,BestFitness,BestProbeNumber%,BestTimeStep&,FunctionName$,RepositionFactor$,Gamma)

DECLARE SUB StatusWindow(FunctionName$,StatusWindowHandle???)

DECLARE SUB
PlotResults(FunctionName$,Nd%,Np%,BestFitnessOverall,BestNpNd%,BestGamma,Neval&&,Rbest(),Mbest(),BestProbeNumberOverall%,BestTimeStepOverall&,LastSte
pBestRun&,Alpha,Beta)

DECLARE SUB
DisplayRunParameters(FunctionName$,Nd%,Np%,Nt&,G,DeltaT,Alpha,Beta,Frep,R(),A(),M(),PlaceInitialProbes$,InitialAcceleration$,RepositionFactor$,RunCFO
$,ShrinkDS$,CheckForEarlyTermination$)

DECLARE SUB GetBestFitness(M(),Np%,StepNumber&,BestFitness,BestProbeNumber%,BestTimeStep&)

DECLARE SUB
TabulateIDprobeCoordinates(Max1DprobesPlotted%,Nd%,Np%,LastStep&,G,DeltaT,Alpha,Beta,Frep,R(),M(),PlaceInitialProbes$,InitialAcceleration$,Reposition
Factor$,FunctionName$,Gamma)

DECLARE SUB
GetPlotAnnotation(PlotAnnotation$,Nd%,Np%,Nt&,G,DeltaT,Alpha,Beta,Frep,M(),PlaceInitialProbes$,InitialAcceleration$,RepositionFactor$,FunctionName$,G
amma)

DECLARE SUB
ChangeRunParameters(NumProbesPerDimension%,Np%,Nd%,Nt&,G,Alpha,Beta,DeltaT,Frep,PlaceInitialProbes$,InitialAcceleration$,RepositionFactor$,FunctionNa
me$)

DECLARE SUB CLEANUP

DECLARE SUB
Plot1DProbePositions(Max1DprobesPlotted%,Nd%,Np%,LastStep&,G,DeltaT,Alpha,Beta,Frep,R(),M(),PlaceInitialProbes$,InitialAcceleration$,RepositionFactor
$,FunctionName$,Gamma)

DECLARE SUB DisplayMmatrix(Np%,Nt&,M())

DECLARE SUB DisplayMbestMatrix(Np%,Nt&,Mbest())

DECLARE SUB DisplayMmatrixThisTimeStep(Np%,j&,M())

DECLARE SUB DisplayAmatrix(Np%,Nd%,Nt&,A())

DECLARE SUB DisplayAmatrixThisTimeStep(Np%,Nd%,j&,A())

DECLARE SUB DisplayRmatrix(Np%,Nd%,Nt&,R())

DECLARE SUB DisplayRmatrixThisTimeStep(Np%,Nd%,j&,R(),Gamma)

DECLARE SUB DisplayXiminMax(Nd%,xiMin(),XiMax())

DECLARE SUB DisplayRunParameters2(FunctionName$,Nd%,Np%,Nt&,G,DeltaT,Alpha,Beta,Frep,PlaceInitialProbes$,InitialAcceleration$,RepositionFactor$)

DECLARE SUB
PlotBestProbeVSTimeStep(Nd%,Np%,LastStep&,G,DeltaT,Alpha,Beta,Frep,M(),PlaceInitialProbes$,InitialAcceleration$,RepositionFactor$,FunctionName$,Gamma
)

DECLARE SUB
PlotBestFitnessEvolution(Nd%,Np%,LastStep&,G,DeltaT,Alpha,Beta,Frep,M(),PlaceInitialProbes$,InitialAcceleration$,RepositionFactor$,FunctionName$,Gamm
a)

DECLARE SUB
PlotAverageDistance(Nd%,Np%,LastStep&,G,DeltaT,Alpha,Beta,Frep,M(),PlaceInitialProbes$,InitialAcceleration$,RepositionFactor$,FunctionName$,R(),DiagL
ength,Gamma)

DECLARE SUB Plot2Dfunction(FunctionName$,R())

DECLARE SUB Plot1Dfunction(FunctionName$,R())

DECLARE SUB
GetFunctionRunParameters(FunctionName$,Nd%,Np%,Nt&,G,DeltaT,Alpha,Beta,Frep,R(),A(),M(),DiagLength,PlaceInitialProbes$,InitialAcceleration$,Repositio
nFactor$)

DECLARE SUB InitialProbeDistribution(Np%,Nd%,Nt&,R(),PlaceInitialProbes$,Gamma)

DECLARE SUB RetrieveErrantProbes(Np%,Nd%,j&,R(),Frep)

DECLARE SUB Retrieveerrantprobes2(Np%,Nd%,j&,R(),A(),Frep)

DECLARE SUB CFO(FunctionName$,Nd%,Nt&,R(),A(),M(),DiagLength,BestFitnessOverall,BestNpNd%,BestGamma,Neval&&,Rbest(),_
              Mbest(),BestProbeNumberOverall%,BestTimeStepOverall&,LastStepBestRun&,Alpha,Beta)  'Self-contained CFO routine -> NO USER-SPECIFIED
PARAMETERS

DECLARE SUB IPD(Np%,Nd%,Nt&,R(),Gamma)

DECLARE SUB ResetDecisionSpaceBoundaries(Nd%)

DECLARE SUB ThreeDplot(PlotFileName$,PlotTitle$,Annotation$,xCoord$,yCoord$,zCoord$,_
              XaxisLabel$,YaxisLabel$,ZaxisLabel$,zMin$,zMax$,GnuPlotEXE$,A$)

DECLARE SUB ThreeDplot2(PlotFileName$,PlotTitle$,Annotation$,xCoord$,yCoord$,zCoord$,XaxisLabel$,_
              YaxisLabel$,ZaxisLabel$,zMin$,zMax$,GnuPlotEXE$,A$,xStart$,xStop$,yStart$,yStop$)

DECLARE SUB ThreeDplot3(PlotFileName$,PlotTitle$,Annotation$,xCoord$,yCoord$,zCoord$,_
              XaxisLabel$,YaxisLabel$,ZaxisLabel$,zMin$,zMax$,GnuPlotEXE$,xStart$,xStop$,yStart$,yStop$)

DECLARE SUB TwoDplot(PlotFileName$,PlotTitle$,xCoord$,yCoord$,XaxisLabel$,YaxisLabel$,_
              LogXaxis$,LogYaxis$,xMin$,xMax$,yMin$,yMax$,xTics$,yTics$,GnuPlotEXE$,LineType$,Annotation$)

DECLARE SUB TwoDplot2Curves(PlotFileName1$,PlotFileName2$,PlotTitle$,Annotation$,xCoord$,yCoord$,XaxisLabel$,YaxisLabel$,_
              LogXaxis$,LogYaxis$,xMin$,xMax$,yMin$,yMax$,xTics$,yTics$,GnuPlotEXE$,LineSize)
```

```
DECLARE SUB TwoDplot3curves(NumCurves%,PlotFileName1$,PlotFileName2$,PlotFileName3$,PlotTitle$,Annotation$,xCoord$,yCoord$,XaxisLabel$,YaxisLabel$, _
                            LogXaxis$,LogYaxis$,xMin$,xMax$,yMin$,yMax$,xTics$,yTics$,GnuPlotEXE$)

DECLARE SUB CreateGNUplotINIfile(PlotWindowULC_X%,PlotWindowULC_Y%,PlotWindowWidth%,PlotWindowHeight%)

DECLARE SUB Delay(NumSecs)

DECLARE SUB MathematicalConstants

DECLARE SUB AlphabetAndDigits

DECLARE SUB SpecialSymbols

DECLARE SUB EMconstants

DECLARE SUB ConversionFactors

DECLARE SUB ShowConstants

DECLARE SUB GetNECdata(NumFreqs%,Zo,FrequencyMHZ(),RadEfficiencyPCT(),MaxGainDBI(),MinGainDBI(),RinOhms(),XinOhms(),VSWR(),FileStatus$)

DECLARE SUB ComplexMultiply(ReA,ImA,ReB,ImB,ReC,ImC)

DECLARE SUB ComplexDivide(ReA,ImA,ReB,ImB,ReC,ImC)

'------ FUNCTION DECLARATIONS -------

DECLARE FUNCTION Probeweight2(Nd%,Np%,R(),M(),p%,j&)

DECLARE FUNCTION Probeweight(Nd%,R(),p%,j&) 'computes a 'weighting factor' based on probe's position (greater weight if closer to decision space
boundary)

DECLARE FUNCTION SlopeRatio(M(),Np%,StepNumber&)

DECLARE CALLBACK FUNCTION DlgProc

DECLARE FUNCTION HasFITNESSsaturated$(Nsteps&,j&,Np%,Nd%,M(),R(),DiagLength)

DECLARE FUNCTION HasDAVGsaturated$(Nsteps&,j&,Np%,Nd%,M(),R(),DiagLength)

DECLARE FUNCTION OscillationInDavg$(j&,Np%,Nd%,M(),R(),DiagLength)

DECLARE FUNCTION DavgThisStep(j&,Np%,Nd%,M(),R(),DiagLength)

DECLARE FUNCTION NoSpaces$(X,NumDigits%)

DECLARE FUNCTION FormatFP$(X,Ndigits%)

DECLARE FUNCTION FormatInteger$(M%)

DECLARE FUNCTION TerminateNowForSaturation$(j&,Nd%,Np%,Nt&,G,DeltaT,Alpha,Beta,R(),A(),M())

DECLARE FUNCTION MagVector(V(),N%)

DECLARE FUNCTION UniformDeviate(u&&)

DECLARE FUNCTION RandomNum(a,b)

DECLARE FUNCTION GaussianDeviate(Mu,Sigma)

DECLARE FUNCTION UnitStep(X)

DECLARE FUNCTION Fibonacci&&(N%)

DECLARE FUNCTION ObjectiveFunction(R(),Nd%,p%,j&,FunctionName$)

DECLARE FUNCTION UnitStep(X)

DECLARE FUNCTION FP2String2$(X!)

DECLARE FUNCTION FP2String$(X)

DECLARE FUNCTION Int2String$(X%)

DECLARE FUNCTION StandingWaveRatio(Zo,ReZ,ImZ)

'=======================================================================================================
'----- MAIN PROGRAM ------

FUNCTION PBMAIN () AS LONG

'    ------ CFO Parameters -----

    LOCAL Nd%, Np%, Nt&

    LOCAL G, DeltaT, Alpha, Beta, Frep AS EXT

    LOCAL PlaceInitialProbes$, InitialAcceleration$, RepositionFactor$

    LOCAL R(), A(), M(), Rbest(), Mbest() AS EXT    'position, acceleration & fitness matrices

    LOCAL FunctionName$          'name of objective function

'    ----------------------- Miscellaneous Setup Parameters ----------------------

    LOCAL N%, i%, YN&, Neval&&, NevalTotal&&, BestNpNd%, NumTrajectories%, Max1DprobesPlotted%, LastStepBestRun&, Pass%

    LOCAL A$, RunCFO$, CFOversion$, NECfileError$, RunStart$, RunStop$

    LOCAL BestGamma, BestFitnessThisRun, BestFitnessOverall, StartTime, StopTime, OptimumFitness AS EXT

    LOCAL BestProbeNum&, BestTimeStep&, BestProbeNumberOverall%, BestTimeStepOverall&, StatusWindowHandle???

'    ------------------- Global Constants --------------------

    REDIM Aij(1 TO 2, 1 TO 25) '(GLOBAL array for Shekel's Foxholes function)

    CALL FillArrayAij

    CALL MathematicalConstants 'NOTE: Calling order is important!!

    CALL AlphabetAndDigits

    CALL SpecialSymbols
```

```
        CALL EMconstants

        CALL ConversionFactors        ': CALL ShowConstants 'to verify constants have been set

        MonopoleSecsPerRun = 0.55##

'       -------------------------- General Setup -----------------------------

        CFOversion$ = "CFO Ver. 03-04-2011"

        RANDOMIZE TIMER  'seed random number generator with program start time

        DESKTOP GET SIZE TO ScreenWidth&, ScreenHeight&  'get screen size (global variables)

        IF DIR$("wgnuplot.exe") = "" THEN

            MSGBOX("WARNING!  'wgnuplot.exe' not found.  Run terminated.") : EXIT FUNCTION

        END IF

'       -------------------------------------------------------------------- CFO RUN PARAMETERS --------------------------------------------------
------------------

        CALL GetTestFunctionNumber(FunctionName$)' : exit function 'DEBUG

        CALL
GetFunctionRunParameters(FunctionName$,Nd%,Np%,Nt&,G,DeltaT,Alpha,Beta,Frep,R(),A(),M(),DiagLength,PlaceInitialProbes$,InitialAcceleration$,Repositio
nFactor$) 'NOTE: Parameters returned but not used in this version!!

        REDIM R(1 TO Np%, 1 TO Nd%, 0 TO Nt&), A(1 TO Np%, 1 TO Nd%, 0 TO Nt&) 'position, acceleration & fitness matrices

        REDIM Rbest(1 TO Np%, 1 TO Nd%, 0 TO Nt&), Mbest(1 TO Np%, 0 TO Nt&) 'overall best position & fitness matrices

'       -------- PLOT 1D and 2D FUNCTIONS ON-SCREEN FOR VISUALIZATION --------

        IF Nd% = 2 AND INSTR(FunctionName$,"PBM_") > 0 THEN

            CALL CheckNECFiles(NECfileError$)

            IF NECfileError$ = "YES" THEN
                EXIT FUNCTION
            ELSE
                MSGBOX("Begin computing plot of function "+FunctionName$+"?  May take a while - be patient...")
            END IF

        END IF

        SELECT CASE Nd%
            CASE 1 : CALL Plot1Dfunction(FunctionName$,R()) : REDIM R(1 TO Np%, 1 TO Nd%, 0 TO Nt&) 'erases coordinate data in R()used to plot function
            CASE 2 : CALL Plot2Dfunction(FunctionName$,R()) : REDIM R(1 TO Np%, 1 TO Nd%, 0 TO Nt&) 'ditto
        END SELECT

'       --------------------------------------------------------------------- RUN CFO ---------------------------------------------------------------------

        YN& = MSGBOX("RUN CFO ON FUNCTION " + FunctionName$ + "?"+CHR$(13)+CHR$(13)+"Get some coffee & sit back...",%MB_YESNO,"CONFIRM RUN") : IF YN& =
%IDYES THEN RunCFO$ = "YES"

        IF RunCFO$ = "YES" THEN

            StartTime = TIMER : RunStart$ = "Started at "+TIMES$+", "+DATE$

            CALL
CFO(FunctionName$,Nd%,Nt&,R(),A(),M(),DiagLength,BestFitnessOverall,BestNpNd%,BestGamma,Neval&,Rbest(),Mbest(),BestProbeNumberOverall%,BestTimeStepO
verall&,LastStepBestRun&,Alpha,Beta)

            StopTime = TIMER : RunStop$ = "Ended at "+TIMES$+", "+DATE$

            CALL
PlotResults(FunctionName$,Nd%,Np%,BestFitnessOverall,BestNpNd%,BestGamma,Neval&,Rbest(),Mbest(),BestProbeNumberOverall%,BestTimeStepOverall&,LastSte
pBestRun&,Alpha,Beta)

            IF FunctionName$ = "LD_MONO" THEN

                IF DIR$("BestMono.NEC") <> "" THEN KILL "BestMono.NEC"
                IF DIR$("BestMono.OUT") <> "" THEN KILL "BestMono.OUT"
                OptimumFitness        = LOADED_MONO(Rbest(),Nd%,BestProbeNumberOverall%,BestTimeStepOverall&)
                NAME "LD_MONO.NEC"        AS            "BestMono.NEC"
                NAME "LD_MONO.OUT"        AS            "BestMono.OUT"

                N% = FREEFILE : OPEN "BestMono.NEC" FOR APPEND AS #N%
                    PRINT #N%,"" : PRINT #N%,"OPTIMUM FITNESS = "+STR$(ROUND(OptimumFitness,3))+", "+DATE$+", "+TIME$
                    PRINT #N%,"" : PRINT #N%,USING$("Best Gamma: ##.###      BestNp/Nd: ###      Nt: #####      Neval: #######      LastStep:
#####",BestGamma,BestNpNd,Nt&,Neval&,LastStepBestRun&)
                    PRINT #N%,""
                    PRINT #N%,RunStart$
                    PRINT #N%,RunStop$
                CLOSE #N%

                N% = FREEFILE : OPEN "BestMono.OUT" FOR APPEND AS #N%
                    PRINT #N%,"" : PRINT #N%,"OPTIMUM FITNESS = "+STR$(ROUND(OptimumFitness,3))+", "+DATE$+", "+TIME$
                    PRINT #N%,"" : PRINT #N%,USING$("Best Gamma: ##.###      BestNp/Nd: ###      Nt: #####      Neval: #######      LastStep:
#####",BestGamma,BestNpNd,Nt&,Neval&,LastStepBestRun&)
                    PRINT #N%,""
                    PRINT #N%,RunStart$
                    PRINT #N%,RunStop$
                CLOSE #N%
            END IF

            MSGBOX(FunctionName$+CHR$(13)+"Total Function Evaluations = "+STR$(Neval&)+CHR$(13)+"Runtime = "+STR$(ROUND((StopTime-
StartTime)/3600##,2))+" hrs, Avg Time/Run = " +STR$(ROUND((StopTime-StartTime)/Neval&,6)))

        END IF

ExitPBMAIN:

END FUNCTION 'PBMAIN()
'==================================================================================== CFO SUBROUTINE
=================================================================================

SUB
CFO(FunctionName$,Nd%,Nt&,R(),A(),M(),DiagLength,BestFitnessOverall,BestNpNd%,BestGamma,Neval&,Rbest(),Mbest(),BestProbeNumberOverall%,BestTimeStepO
verall&,LastStepBestRun&,Alpha,Beta)

LOCAL p%, i%, j& 'Standard Indices: Probe #, Coordinate #, Time Step #

LOCAL Np% 'Number of Probes
```

```
LOCAL MaxProbesPerDimension% 'Maximum # probes per dimension (depends on Nd%)

LOCAL k%, L% 'Dummy summation indices

LOCAL NumProbesPerDimension%, GammaNumber%, NumGammas% 'Probes/dimension on each probe line; gamm point #; # gamma points

LOCAL SumSQ, Denom, Numerator, Gamma, Frep, DeltaFrep AS EXT

LOCAL BestProbeNumber%, BestTimeStep&, LastStep&, BestProbeNumberThisRun%, BestTimeStepThisRun&

LOCAL BestFitness, BestFitnessThisRun, eta AS EXT

LOCAL FitnessSaturation$
'--------------- Initial Parameter Values ------------------

Nt& = 1000 'set to a large value expecting early termination

Alpha = 1## : Beta = 1## 'CFO EXPONENTS Alpha & Beta

IF FunctionName$ = "F7" THEN Nt& = 100 'to reduce runtime because this function contains random noise

IF FunctionName$ = "LD_MONO" THEN Nt& = 200 'to reduce runtime for Loaded Monopole

LastStep& = Nt&

Neval1&& = 0

DeltaFrep = 0.1##

NumGammas% = 4 '11

SELECT CASE Nd% 'set Np%/Nd% based on Nd% to avoid excessive runtimes
    CASE   1 TO   6 : MaxProbesPerDimension% = 14
    CASE   7 TO  10 : MaxProbesPerDimension% = 12
    CASE  11 TO  15 : MaxProbesPerDimension% = 10
    CASE  16 TO  20 : MaxProbesPerDimension% =  8
    CASE  21 TO  30 : MaxProbesPerDimension% =  6
    CASE ELSE       : MaxProbesPerDimension% =  4

END SELECT

IF FunctionName$ = "LD_MONO" THEN
    MaxProbesPerDimension% = 6 'to reduce runtime for Loaded Monopole
    MSGBOX("Approximate runtime: "+STR$(ROUND(MonopoleSecsPerRun*NumGammas*Nd%*MaxProbesPerDimension%*Nt&/40##,1))+" minutes")
END IF
'  ----------------------- Np/Nd LOOP --------------------

BestFitnessOverall = -1E4200 'very large negative number...

FOR NumProbesPerDimension% = 2 TO MaxProbesPerDimension% STEP 2

'FOR NumProbesPerDimension% = 4 TO MaxProbesPerDimension% STEP 2 'original code

Np% = NumProbesPerDimension%*Nd%

'  ----------------------- GAMMA LOOP --------------------

FOR GammaNumber% = 1 TO NumGammas%

Gamma = (GammaNumber%-1)/(NumGammas%-1)

REDIM R(1 TO Np%, 1 TO Nd%, 0 TO Nt&), A(1 TO Np%, 1 TO Nd%, 0 TO Nt&), M(1 TO Np%, 0 TO Nt&) 're-initializes Position Vector/Acceleration/Fitness
matrices to zero
'STEP (A1) ------------- Compute Initial Probe Distribution (Step 0)-----------------

    CALL IPD(Np%,Nd%,Nt&,R(),Gamma) 'Probe Line IPD intersecting on diagonal at a point determined by Gamma
'STEP (A2) ------------- Compute Initial Fitness Matrix (Step 0) ------------------------

    FOR p% = 1 TO Np% : M(p%,0) = ObjectiveFunction(R(),Nd%,p%,0,FunctionName$) : INCR Neval&& : NEXT p%

'STEP (A3) ------------- Set Initial Probe Accelerations to ZERO (Step 0)---------------

    FOR p% = 1 TO Np% : FOR i% = 1 TO Nd% : A(p%,i%,0) = 0## : NEXT i% : NEXT p%

'STEP (A4) ------------- Initialize Frep ----------------

    Frep = 0.5##

'  ===================================== LOOP ON TIME STEPS STARTING AT STEP #1 =====================================

    BestFitnessThisRun = M(1,0)

    FOR j& = 1 TO Nt&

'STEP (B) ---------- Compute Probe Position Vectors for this Time Step --------

        FOR p% = 1 TO Np% : FOR i% = 1 TO Nd% : R(p%,i%,j&) = R(p%,i%,j&-1) + A(p%,i%,j&-1) : NEXT i% : NEXT p% 'note: factor of 1/2 combined with
G=2 to produce unity coefficient
'STEP (C) ---------- Retrieve Errant Probes ---------------

        CALL RetrieveErrantProbes(Np%,Nd%,j&,R(),Frep)

'          CALL RetrieveErrantprobes2(Np%,Nd%,j&,R(),A(),Frep) 'added 04-01-10

'STEP (D) ---------- Compute Fitness Matrix for Current Probe Distribution ---------

        FOR p% = 1 TO Np% : M(p%,j&) = ObjectiveFunction(R(),Nd%,p%,j&,FunctionName$) : INCR Neval&& : NEXT p%

'STEP (E) ---------- Compute Accelerations Based on Current Probe Distribution & Fitnesses ---------------

        FOR p% = 1 TO Np%

            FOR i% = 1 TO Nd%

                A(p%,i%,j&) = 0

                FOR k% = 1 TO Np%

                    IF k% <> p% THEN
```

```basic
                    SumSQ = 0## : FOR L% = 1 TO Nd%  : SumSQ = SumSQ + (R(k%,L%,j&)-R(p%,L%,j&))^2 : NEXT L% 'dummy index

                    IF SumSQ <> 0## THEN 'to avoid zero denominator (added 03-20-10)

                        Denom = SQR(SumSQ) : Numerator = UnitStep((M(k%,j&)-M(p%,j&)))*(M(k%,j&)-M(p%,j&))

                        A(p%,i%,j&) = A(p%,i%,j&) + (R(k%,i%,j&)-R(p%,i%,j&))*Numerator^Alpha/Denom^Beta 'ORIGINAL VERSION WITH VARIABLE Alpha &
Beta

                    END IF 'added 03-20-10

                END IF

            NEXT k% 'dummy index

        NEXT i% 'coord (dimension) #

    NEXT p% 'probe #
'   --------- Get Best Fitness & Corresponding Probe # and Time Step ---------
    CALL GetBestFitness(M(),Np%,j&,BestFitness,BestProbeNumber%,BestTimeStep&)

    IF BestFitness >= BestFitnessThisRun THEN

        BestFitnessThisRun = BestFitness : BestProbeNumberThisRun% = BestProbeNumber% : BestTimeStepThisRun& = BestTimeStep&

    END IF
'   ----- Increment Frep -----
    Frep = Frep + DeltaFrep

    IF Frep > 1## THEN Frep = 0.05## 'keep Frep in range [0.05,1]
'   --------- Starting at Step #20 Shrink Decision Space Around Best Probe Every 20th Step -----------
'   IF j& MOD 10 = 0 AND j& >= 20 THEN
    IF j& MOD 20 = 0 AND j& >= 20 THEN

        FOR i% = 1 TO Nd% : XiMin(i%) = XiMin(i%)+(R(BestProbeNumber%,i%,BestTimeStep&)-XiMin(i%))/2## : XiMax(i%) = XiMax(i%)-(XiMax(i%)-
R(BestProbeNumber%,i%,BestTimeStep&))/2## : NEXT i% 'shrink DS by 0.5

        CALL RetrieveErrantProbes(Np%,Nd%,j&,R(),Frep) 'TO RETRIEVE PROBES LYING OUTSIDE SHRUNKEN DS 'ADDED 02-07-2010
'           CALL RetrieveErrantProbes2(Np%,Nd%,j&,R(),A(),Frep) 'added 04-01-10

    END IF
'   --------- If SlopeRatio Changes Abruptly, Shrink Decision Space Around Best Probe -----------
'   eta = 0.9##
'   IF SlopeRatio(M(),Np%,j&) >= 3## THEN
'       FOR i% = 1 TO Nd% : XiMin(i%) = XiMin(i%)+eta*(R(BestProbeNumber%,i%,BestTimeStep&)-XiMin(i%)) : XiMax(i%) = XiMax(i%)-eta*(XiMax(i%)-
R(BestProbeNumber%,i%,BestTimeStep&)): NEXT i%
'           CALL RetrieveErrantProbes(Np%,Nd%,j&,R(),Frep) 'TO RETRIEVE PROBES LYING OUTSIDE SHRUNKEN DS 'ADDED 02-07-2010
'           CALL Retrieveerrantprobes2(Np%,Nd%,j&,R(),A(),Frep) 'added 04-01-10
'   END IF
'STEP (F) ---------- Check for Early Run Termination ---------
    IF HasFITNESSsaturated$(25,j&,Np%,Nd%,M(),R(),DiagLength) = "YES" THEN
        LastStep& = j&
        EXIT FOR
    END IF
    NEXT j& 'END TIME STEP LOOP
'--------------- Best Overall Fitness & Corresponding Parameters -------------------
IF BestFitnessThisRun >= BestFitnessOverall THEN
    BestFitnessOverall = BestFitnessThisRun : BestProbeNumberOverall% = BestProbeNumberThisRun% : BestTimeStepOverall& = BestTimeStepThisRun&
    BestNpNd% = NumProbesPerDimension% : BestGamma = Gamma : LastStepBestRun& = LastStep&
    CALL CopyBestMatrices(Np%,Nd%,Nt&,R(),M(),Rbest(),Mbest())
END IF
'STEP (G) ----- Reset Decision Space Boundaries to Initial Values -----
CALL ResetDecisionSpaceBoundaries(Nd%)
NEXT GammaNumber% 'END GAMMA LOOP
NEXT NumProbesPerDimension% 'END Np/Nd LOOP
END SUB 'CFO()
'=================================================================================================================================================
===============
SUB IPD(Np%,Nd%,Nt&,R(),Gamma)

LOCAL DeltaXi, Delx1, Delx2, Di AS EXT

LOCAL NumProbesPerDimension%, p%, i%, k%, NumX1points%, NumX2points%, x1pointNum%, x2pointNum%

            IF Nd% > 1 THEN

                NumProbesPerDimension% = Np%\Nd% 'even #

            ELSE

                NumProbesPerDimension% = Np%

            END IF
```

```
                    FOR i% = 1 TO Nd%

                        FOR p% = 1 TO Np%

                            R(p%,i%,0) = XiMin(i%) + Gamma*(XiMax(i%)-XiMin(i%))

                        NEXT Np%

                    NEXT i%

                    FOR i% = 1 TO Nd% 'place probes probe line-by-probe line (i% is dimension number)

                        Deltaxi = (XiMax(i%)-XiMin(i%))/(NumProbesPerDimension%-1)

                        FOR k% = 1 TO NumProbesPerDimension%

                            p% = k% + NumProbesPerDimension%*(i%-1) 'probe #

                            R(p%,i%,0) = XiMin(i%) + (k%-1)*Deltaxi

                        NEXT k%

                    NEXT i%

END SUB 'IPD()

'----

FUNCTION HasFITNESSsaturated$(Nsteps&,j&,Np%,Nd%,M(),R(),DiagLength)

LOCAL A$, B$

LOCAL k&, p%

LOCAL BestFitness, SumOfBestFitnesses, BestFitnessStepJ, FitnessSatTOL AS EXT

    A$ = "NO" : B$ = "j="+STR$(j&)+CHR$(13)

    FitnessSatTOL = 0.000001## 'tolerance for FITNESS saturation

    IF j& < Nsteps& + 10 THEN GOTO ExitHasFITNESSsaturated 'execute at least 10 steps after averaging interval before performing this check

    SumOfBestFitnesses = 0##

    FOR k& = j&-Nsteps&+1 TO j& 'GET BEST FITNESSES STEP-BY-STEP FOR Nsteps& INCLUDING THIS STEP j& AND COMPUTE AVERAGE VALUE.

'       BestFitness = M(k&,1) 'ORIG CODE 03-23-2010: THIS IS A MISTAKE!

        BestFitness = -1E4200 'THIS LINE CORRECTED 03-23-2010 PER DISCUSSION WITH ROB GREEN.
                              'INITIALIZE BEST FITNESS AT k&-th TIME STEP TO AN EXTREMELY LARGE NEGATIVE NUMBER.

        FOR p% = 1 TO Np% 'PROBE-BY-PROBE GET MAXIMUM FITNESS

            IF M(p%,k&) >= BestFitness THEN BestFitness = M(p%,k&)

        NEXT p%

        IF k& = j& THEN BestFitnessStepJ = BestFitness 'IF AT THE END OF AVERAGING INTERVAL, SAVE BEST FITNESS FOR CURRENT TIME STEP j&

        SumOfBestFitnesses = SumOfBestFitnesses + BestFitness

        B$ = B$ + "k="+STR$(k&)+"  BestFit="+STR$(BestFitness)+"   SumFit="+STR$(SumOfBestFitnesses)+CHR$(13)

    NEXT k&

    IF ABS(SumOfBestFitnesses/Nsteps&-BestFitnessStepJ) =< FitnessSatTOL THEN A$ = "YES" 'saturation if (avg value - last value) are within TOL

ExitHasFITNESSsaturated:

    HasFITNESSsaturated$ = A$

END FUNCTION 'HasFITNESSsaturated$()

'-----------------------------------

SUB Retrieveerrantprobes2(Np%,Nd%,j&,R(),A(),Frep) 'added 04-01-10

LOCAL ErrantProbe$

LOCAL p%, i%, k%

LOCAL Xik, dMax, Eta(), EtaStar, SumSQ, MagAj1 AS EXT

REDIM Eta(1 TO Nd%, 1 TO 2) 'Eta(i%,k%)

    FOR p% = 1 TO Np%

        ErrantProbe$ = "NO" 'presume each probe is inside DS

        FOR i% = 1 TO Nd% 'check to see if probe p lies outside DS (any coordinate exceeding a boundary)

            IF (R(p%,i%,j&) > XiMax(i%) OR R(p%,i%,j&) < XiMin(i%)) AND A(p%,i%,j&-1) <> 0## THEN 'probe lies outside DS

                ErrantProbe$ = "YES" : EXIT FOR

            END IF

        NEXT i%

        IF ErrantProbe$ = "YES" THEN 'reposition probe p inside DS with acceleration direction preserved

            FOR i% = 1 TO Nd% 'compute array of Eta values

                FOR k% = 1 TO 2

                    SELECT CASE k%
                        CASE 1 : Xik = XiMin(i%)
                        CASE 2 : Xik = XiMax(i%)
                    END SELECT

                    Eta(i%,k%) = (Xik-R(p%,i%,j&-1))/A(p%,i%,j&-1)

                NEXT k%

            NEXT i%

            EtaStar = 1E4200 'very large positive number
```

```
                FOR i% = 1 TO Nd% 'get min Eta value > 0
                    FOR k% = 1 TO 2
                        IF Eta(i%,k%) =< EtaStar AND Eta(i%,k%) >= 0## THEN EtaStar = Eta(i%,k%)
                    NEXT k%

                NEXT i%
                SumSQ = 0## : FOR i% = 1 TO Nd% : SumSQ = SumSQ + A(p%,i%,j&-1)^2 : NEXT i% : MagAj1 = SQR(SumSQ) 'magnitude of acceleration vector
A(probe #p) at step j&-1

                dMax = EtaStar*MagAj1 'distance to nearest boundary plane from position of probe p% at step j&-1

                FOR i% = 1 TO Nd% 'change probe p's i-th coordinate in proportion to the ratio of position vector lengths at steps j& and j&-1
                    R(p%,i%,j&) = R(p%,i%,j&-1) + Frep*dMax*A(p%,i%,j&-1)/MagAj1 'preserves acceleration directional information by scaling dMax by Frep,
which is arbitrary but precisely known
                NEXT i%
            END IF 'ErrantProbe$ = "YES"
        NEXT p%
END SUB 'Retrieveerrantprobes2()
'-----------------------------
SUB RetrieveErrantProbes(Np%,Nd%,j&,R(),Frep)
LOCAL p%, i%
    FOR p% = 1 TO Np%
        FOR i% = 1 TO Nd%
            IF R(p%,i%,j&) < XiMin(i%) THEN R(p%,i%,j&) = MAX(XiMin(i%) + Frep*(R(p%,i%,j&-1)-XiMin(i%)),XiMin(i%)) 'CHANGED 02-07-10
            IF R(p%,i%,j&) > XiMax(i%) THEN R(p%,i%,j&) = MIN(XiMax(i%) - Frep*(XiMax(i%)-R(p%,i%,j&-1)),XiMax(i%))
        NEXT i%
    NEXT p%
END SUB 'RetrieveErrantProbes()
'-----------------------------
SUB ResetDecisionSpaceBoundaries(Nd%)
    LOCAL i%
    FOR i% = 1 TO Nd% : XiMin(i%) = StartingXiMin(i%) : XiMax(i%) = StartingXiMax(i%) : NEXT i%
END SUB
'------
SUB CopyBestMatrices(Np%,Nd%,Nt&,R(),M(),Rbest(),Mbest())
LOCAL p%, i%, j&
REDIM Rbest(1 TO Np%, 1 TO Nd%, 0 TO Nt&), Mbest(1 TO Np%, 0 TO Nt&) 're-initializes Best Position Vetor/Fitness matrices to zero
    FOR p% = 1 TO Np%
        FOR i% = 1 TO Nd%
            FOR j& = 0 TO Nt&
                Rbest(p%,i%,j&) = R(p%,i%,j&) : Mbest(p%,j&) = M(p%,j&)
            NEXT j&
        NEXT i%
    NEXT p%
END SUB
'------
'FORGET THIS IDEA !!!!
FUNCTION Probeweight(Nd%,R(),p%,j&) 'computes a 'weighting factor' based on probe's position (greater weight if closer to decision space boundary)
LOCAL MinDistCoordinate%, i% 'Dimension number.  Remember, XiMin(), XiMax()& DiagLength are GLOBAL.
LOCAL MinDistance, dStar, MaxWeight AS EXT
    MinDistance = DiagLength 'largest dimension of decision space
    FOR i% = 1 TO Nd% 'compute distance to closest boundary
        IF ABS(R(p%,i%,j&)-XiMin(i%)) =< MinDistance THEN
            MinDistance = ABS(R(p%,i%,j&)-XiMin(i%)) : MinDistCoordinate% = i%
        END IF
        IF ABS(XiMax(i%)-R(p%,i%,j&)) =< MinDistance THEN
            MinDistance = ABS(XiMax(i%)-R(p%,i%,j&)) : MinDistCoordinate% = i%
        END IF
    NEXT i%
    dStar = MinDistance/(XiMax(MinDistCoordinate%)-XiMin(MinDistCoordinate%)) 'normalized minimum distance, [0-1]
    MaxWeight = 2##
'    Probeweight = 1## + 2##*MaxWeight*abs(dStar-0.5##)
```

```
        Probeweight = 1## + 4##*Maxweight*(dStar-0.5##)^2

END FUNCTION 'Probeweight
'------------------------

'FORGET THIS IDEA !!!!

FUNCTION Probeweight2(Nd%,Np%,R(),M(),p%,j&) 'computes a 'weighting factor' based on probe's position (greater weight if closer to decision space boundary)

LOCAL MinDistCoordinate%, ProbeNum%, BestProbeThisStep%, i% 'Dimension number.  Remember, XiMin(), XiMax()& DiagLength are GLOBAL.

LOCAL Distance, SumSQ, dStar, MaxWeight, BestFitnessThisStep AS EXT

    BestFitnessThisStep = M(1,j&)

    FOR ProbeNum% = 1 TO Np% 'get number of best probe this step

        IF M(ProbeNum%,j&) >= BestFitnessThisStep THEN

            BestFitnessThisStep = M(ProbeNum%,j&) : BestProbeThisStep% = ProbeNum%

        END IF

    NEXT ProbeNum%

    SumSQ = 0##

    FOR i% = 1 TO Nd% 'compute distance from probe #p to thes best probe this step

        SumSQ = (R(p%,i%,j&) - R(BestProbeThisStep%,i%,j&))^2

    NEXT i%

    Distance = SQR(SumSQ)

    dStar = Distance/DiagLength 'range [0-1]
'    ---------------- Compute weight Factor --------------------

    MaxWeight = 0##

'    Probeweight2 = 1## + 2##*MaxWeight*abs(dStar-0.5##)

    Probeweight2 = 1## + 4##*MaxWeight*(dStar-0.5##)^2

END FUNCTION 'Probeweight2
'------------------------

FUNCTION SlopeRatio(M(),Np%,StepNumber&)

LOCAL p% 'probe #

LOCAL NumSteps%

LOCAL BestFitnessAtStepNumber, BestFitnessAtStepNumberMinus1, BestFitnessAtStepNumberMinus2, Z AS EXT

    Z = 1## 'assumes no slope change

    IF StepNumber& < 10 THEN GOTO ExitSlopeRatio 'need at least 10 steps for this test

    NumSteps% = 2

    BestFitnessAtStepNumber       = M(1,StepNumber&)           : FOR p% = 1 TO Np% : IF M(p%,StepNumber&)              >= BestFitnessAtStepNumber
THEN BestFitnessAtStepNumber = M(p%,StepNumber&)           : NEXT p%

    BestFitnessAtStepNumberMinus1 = M(1,StepNumber&-NumSteps%)   : FOR p% = 1 TO Np% : IF M(p%,StepNumber&-NumSteps%)    >=
BestFitnessAtStepNumberMinus1 THEN BestFitnessAtStepNumberMinus1 = M(p%,StepNumber&-NumSteps%)    : NEXT p%

    BestFitnessAtStepNumberMinus2 = M(1,StepNumber&-2*NumSteps%) : FOR p% = 1 TO Np% : IF M(p%,StepNumber&-2*NumSteps%) >=
BestFitnessAtStepNumberMinus2 THEN BestFitnessAtStepNumberMinus2 = M(p%,StepNumber&-2*NumSteps%) : NEXT p%

    Z = (BestFitnessAtStepNumber-BestFitnessAtStepNumberMinus1)/(BestFitnessAtStepNumberMinus1-BestFitnessAtStepNumberMinus2)

ExitSlopeRatio:

    SlopeRatio = Z

END FUNCTION

SUB
PlotResults(FunctionName$,Nd%,Np%,BestFitnessOverall,BestNpNd%,BestGamma,Neval&,Rbest(),Mbest(),BestProbeNumberOverall%,BestTimeStepOverall&,LastSte
pBestRun&,Alpha,Beta)

LOCAL LastStep&, BestFitnessProbeNumber%, BestFitnessTimeStep&, NumTrajectories%, Max1DprobesPlotted%, i%

LOCAL RepositionFactor$, PlaceInitialProbes$, InitialAcceleration$, A$, B$

LOCAL G, DeltaT, Frep AS EXT

    G = 2## : DeltaT = 1## : Frep = 0.5## : RepositionFactor$ = "VARIABLE" : PlaceInitialProbes$ = "UNIFORM ON-AXIS " : InitialAcceleration$ =
"FIXED" 'THESE ARE NOW HARDWIRED IN THE CFO EQUATIONS

    B$ = "" : IF Nd% > 1 THEN B$ = "s"

    A$ = FunctionName$ + CHR$(13) +_
        "Best Fitness = " + REMOVE$(STR$(BestFitnessOverall),ANY" ")    + " returned by" + CHR$(13) +_
        "Probe # "        + REMOVE$(STR$(BestProbeNumberOverall%),ANY" ") +_
        " at Time Step  " + REMOVE$(STR$(BestTimeStepOverall&),ANY" ")    + CHR$(13) + CHR$(13) + "P" + REMOVE$(STR$(BestProbeNumberOverall%),ANY"
") + " coordinate" + B$ + ":" + CHR$(13)

        FOR i% = 1 TO Nd% : A$ = A$ + STR$(i%)+"   "+REMOVE$(STR$(ROUND(Rbest(BestProbeNumberOverall%,i%,BestTimeStepOverall&),8)),ANY" ")+CHR$(13) :
NEXT i%

    MSGBOX(A$)

'   -------------------------------------------- PLOT EVOLUTION OF BEST FITNESS, AVG DISTANCE & BEST PROBE # --------------------------------------
-------------------

    CALL
PlotBestFitnessEvolution(Nd%,Np%,LastStepBestRun&,G,DeltaT,Alpha,Beta,Frep,Mbest(),PlaceInitialProbes$,InitialAcceleration$,RepositionFactor$,Functio
nName$,BestGamma)
```

```
        CALL
PlotAverageDistance(Nd%,Np%,LastStepBestRun&,G,DeltaT,Alpha,Beta,Frep,Mbest(),PlaceInitialProbes$,InitialAcceleration$,RepositionFactor$,FunctionName
$,Rbest(),DiagLength,BestGamma)

        CALL
PlotBestProbeVsTimeStep(Nd%,Np%,LastStepBestRun&,G,DeltaT,Alpha,Beta,Frep,Mbest(),PlaceInitialProbes$,InitialAcceleration$,RepositionFactor$,Function
Name$,BestGamma)

' ----------------------------------------------- PLOT TRAJECTORIES OF BEST PROBES FOR 2/3-D FUNCTIONS -------------------------------------------------
------

    IF Nd% = 2 THEN

        NumTrajectories% = 10 : CALL Plot2DbestProbeTrajectories(NumTrajectories%,Mbest(),Rbest(),Np%,Nd%,LastStepBestRun&,FunctionName$)

        NumTrajectories% = 16 : CALL Plot2DindividualProbeTrajectories(NumTrajectories%,Mbest(),Rbest(),Np%,Nd%,LastStepBestRun&,FunctionName$)

    END IF

    IF Nd% = 3 THEN

        NumTrajectories% = 4 : CALL Plot3DbestProbeTrajectories(NumTrajectories%,Mbest(),Rbest(),Np%,Nd%,LastStepBestRun&,FunctionName$)

    END IF

' ---------- For 1-D Objective Functions, Tabulate Probe Coordinates & if Np% =< Max1DprobesPlotted% Plot Evolution of Probe Positions -----------
-

    IF Nd% = 1 THEN

        Max1DprobesPlotted% = 15

        CALL
Tabulate1DprobeCoordinates(Max1DprobesPlotted%,Np%,LastStepBestRun&,G,DeltaT,Alpha,Beta,Frep,Rbest(),Mbest(),PlaceInitialProbes$,InitialAccelerat
ion$,RepositionFactor$,FunctionName$,BestGamma)

        IF Np% =< Max1DprobesPlotted% THEN _

        CALL
Plot1DprobePositions(Max1DprobesPlotted%,Nd%,Np%,LastStepBestRun&,G,DeltaT,Alpha,Beta,Frep,Rbest(),Mbest(),PlaceInitialProbes$,InitialAcceleration$,R
epositionFactor$,FunctionName$,BestGamma)

        CALL CLEANUP 'delete probe coordinate plot files, if any

    END IF

END SUB 'PlotResults()

'=============================================================================================================================================
===========================================================

SUB CheckNECfiles(NECfileError$)

LOCAL N%

    NECfileError$ = "NO"

' ------------------ NEC Files Required for PBM Antenna Benchmarks ------------------

    IF DIR$("n41_2k1.exe") = "" THEN

        MSGBOX("WARNING!  'n41_2k1.exe' not found.  Run terminated.") : NECfileError$ = "YES" : EXIT SUB

    END IF

    N% = FREEFILE : OPEN "ENDERR.DAT"  FOR OUTPUT AS #N%   : PRINT #N%, "NO"      : CLOSE #N%

    N% = FREEFILE : OPEN "FILE_MSG.DAT" FOR OUTPUT AS #N%   : PRINT #N%, "NO"      : CLOSE #N%

    N% = FREEFILE : OPEN "NHSCALE.DAT"  FOR OUTPUT AS #N% : PRINT #N%, "0.00001" : CLOSE #N%

END SUB

'------

SUB GetBestFitness(M(),Np%,StepNumber&,BestFitness,BestProbeNumber%,BestTimeStep&)

LOCAL p%, i&, A$

    BestFitness = M(1,0)

    FOR i& = 0 TO StepNumber&

        FOR p% = 1 TO Np%

            IF M(p%,i&) >= BestFitness THEN

                BestFitness = M(p%,i&) : BestProbeNumber% = p% : BestTimeStep& = i&

            END IF

        NEXT p%

    NEXT i&

END SUB

'=================================================================== FUNCTION DEFINITIONS ===================================================
=================================================================================================

FUNCTION ObjectiveFunction(R(),Nd%,p%,j&,FunctionName$) 'Objective function to be MAXIMIZED is defined here

    SELECT CASE FunctionName$

        CASE "ParrottF4"    : ObjectiveFunction = ParrottF4(R(),Nd%,p%,j&)        'Parrott F4 (1-D)

        CASE "SGO"          : ObjectiveFunction = SGO(R(),Nd%,p%,j&)              'SGO Function (2-D)

        CASE "GP"           : ObjectiveFunction = GoldsteinPrice(R(),Nd%,p%,j&)   'Goldstein-Price Function (2-D)

        CASE "STEP"         : ObjectiveFunction = StepFunction(R(),Nd%,p%,j&)     'Step Function (n-D)

        CASE "SCHWEFEL_226"  : ObjectiveFunction = Schwefel226(R(),Nd%,p%,j&)      'Schwefel Prob. 2.26 (n-D)

        CASE "COLVILLE"     : ObjectiveFunction = Colville(R(),Nd%,p%,j&)         'Colville Function (4-D)

        CASE "GRIEWANK"     : ObjectiveFunction = Griewank(R(),Nd%,p%,j&)         'Griewank Function (n-D)

        CASE "HIMMELBLAU"   : ObjectiveFunction = Himmelblau(R(),Nd%,p%,j&)       'Himmelblau Function (2-D)
```

```
        CASE "ROSENBROCK"     : ObjectiveFunction = Rosenbrock(R(),Nd%,p%,j&)     'Rosenbrock Function (n-D)

        CASE "SPHERE"         : ObjectiveFunction = Sphere(R(),Nd%,p%,j&)          'Sphere Function (n-D)

        CASE "HIMMELBLAUNLO"  : ObjectiveFunction = HIMMELBLAUNLO(R(),Nd%,p%,j&)   'Himmelblau NLO (5-D)

        CASE "TRIPOD"         : ObjectiveFunction = Tripod(R(),Nd%,p%,j&)          'Tripod (2-D)

        CASE "ROSENBROCKF6"   : ObjectiveFunction = RosenbrockF6(R(),Nd%,p%,j&)    'RosenbrockF6 (10-D)

        CASE "COMPRESSIONSPRING" : ObjectiveFunction = CompressionSpring(R(),Nd%,p%,j&)   'Compression Spring (3-D)

        CASE "GEARTRAIN"      : ObjectiveFunction = GearTrain(R(),Nd%,p%,j&)       'Gear Train (4-D)

'       ------------------------- GSO Paper Benchmark Functions -------------------------

        CASE  "F1"            : ObjectiveFunction = F1(R(),Nd%,p%,j&)              'F1   (n-D)
        CASE  "F2"            : ObjectiveFunction = F2(R(),Nd%,p%,j&)              'F2   (n-D)
        CASE  "F3"            : ObjectiveFunction = F3(R(),Nd%,p%,j&)              'F3   (n-D)
        CASE  "F4"            : ObjectiveFunction = F4(R(),Nd%,p%,j&)              'F4   (n-D)
        CASE  "F5"            : ObjectiveFunction = F5(R(),Nd%,p%,j&)              'F5   (n-D)
        CASE  "F6"            : ObjectiveFunction = F6(R(),Nd%,p%,j&)              'F6   (n-D)
        CASE  "F7"            : ObjectiveFunction = F7(R(),Nd%,p%,j&)              'F7   (n-D)
        CASE  "F8"            : ObjectiveFunction = F8(R(),Nd%,p%,j&)              'F8   (n-D)
        CASE  "F9"            : ObjectiveFunction = F9(R(),Nd%,p%,j&)              'F9   (n-D)
        CASE "F10"            : ObjectiveFunction = F10(R(),Nd%,p%,j&)             'F10  (n-D)
        CASE "F11"            : ObjectiveFunction = F11(R(),Nd%,p%,j&)             'F11  (n-D)
        CASE "F12"            : ObjectiveFunction = F12(R(),Nd%,p%,j&)             'F12  (n-D)
        CASE "F13"            : ObjectiveFunction = F13(R(),Nd%,p%,j&)             'F13  (n-D)
        CASE "F14"            : ObjectiveFunction = F14(R(),Nd%,p%,j&)             'F14  (2-D)
        CASE "F15"            : ObjectiveFunction = F15(R(),Nd%,p%,j&)             'F15  (4-D)
        CASE "F16"            : ObjectiveFunction = F16(R(),Nd%,p%,j&)             'F16  (2-D)
        CASE "F17"            : ObjectiveFunction = F17(R(),Nd%,p%,j&)             'F17  (2-D)
        CASE "F18"            : ObjectiveFunction = F18(R(),Nd%,p%,j&)             'F18  (2-D)
        CASE "F19"            : ObjectiveFunction = F19(R(),Nd%,p%,j&)             'F19  (3-D)
        CASE "F20"            : ObjectiveFunction = F20(R(),Nd%,p%,j&)             'F20  (6-D)
        CASE "F21"            : ObjectiveFunction = F21(R(),Nd%,p%,j&)             'F21  (4-D)
        CASE "F22"            : ObjectiveFunction = F22(R(),Nd%,p%,j&)             'F22  (4-D)
        CASE "F23"            : ObjectiveFunction = F23(R(),Nd%,p%,j&)             'F23  (4-D)

'       ------------------------- PBM Antenna Benchmarks -------------------------

        CASE "PBM_1"          : ObjectiveFunction = PBM_1(R(),Nb%,p%,j&)           'PBM_1 (2-D)
        CASE "PBM_2"          : ObjectiveFunction = PBM_2(R(),Nd%,p%,j&)           'PBM_2 (2-D)
        CASE "PBM_3"          : ObjectiveFunction = PBM_3(R(),Nd%,p%,j&)           'PBM_3 (2-D)
        CASE "PBM_4"          : ObjectiveFunction = PBM_4(R(),Nd%,p%,j&)           'PBM_4 (2-D)
        CASE "PBM_5"          : ObjectiveFunction = PBM_5(R(),Nd%,p%,j&)           'PBM_5 (2-D)

'       ------------------------- Loaded Monopole on PEC Ground Plane -------------------------

        CASE "LD_MONO"        : ObjectiveFunction = LOADED_MONO(R(),Nd%,p%,j&)     'Loaded Monopole

    END SELECT

    END FUNCTION 'ObjectiveFunction()

'------

SUB GetFunctionRunParameters(FunctionName$,Nd%,Np%,Nt&,G,DeltaT,Alpha,Beta,Frep,R(),A(),M(),_
                             DiagLength,PlaceInitialProbes$,InitialAcceleration$,RepositionFactor$)

LOCAL i%, NumProbesPerDimension%, NN%, NumCollinearElements%

    SELECT CASE FunctionName$

        CASE "ParrottF4"

            Nd% = 1 : Np% = 3

            REDIM XiMin(1 TO Nd%), XiMax(1 TO Nd%) : XiMin(1) = 0## : XiMax(1) = 1##

            REDIM StartingXiMin(1 TO Nd%), StartingXiMax(1 TO Nd%) : FOR i% = 1 TO Nd% : StartingXiMin(i%) = XiMin(i%) : StartingXiMax(i%) =
XiMax(i%) : NEXT i%

        CASE "SGO"

            Nd% = 2 : Np% = 8

            REDIM XiMin(1 TO Nd%), XiMax(1 TO Nd%) : FOR i% = 1 TO Nd% : XiMin(i%) = -50## : XiMax(i%) = 50## : NEXT i%

            REDIM StartingXiMin(1 TO Nd%), StartingXiMax(1 TO Nd%) : FOR i% = 1 TO Nd% : StartingXiMin(i%) = XiMin(i%) : StartingXiMax(i%) =
XiMax(i%) : NEXT i%

        CASE "GP"

            Nd% = 2 : Np% = 8

            REDIM XiMin(1 TO Nd%), XiMax(1 TO Nd%) : FOR i% = 1 TO Nd% : XiMin(i%) = -100## : XiMax(i%) = 100## : NEXT i%

            REDIM StartingXiMin(1 TO Nd%), StartingXiMax(1 TO Nd%) : FOR i% = 1 TO Nd% : StartingXiMin(i%) = XiMin(i%) : StartingXiMax(i%) =
XiMax(i%) : NEXT i%

        CASE "STEP"

            Nd% = 2 : Np% = 8

            REDIM XiMin(1 TO Nd%), XiMax(1 TO Nd%) : FOR i% = 1 TO Nd% : XiMin(i%) = -100## : XiMax(i%) = 100## : NEXT i%
'           REDIM XiMin(1 TO Nd%), XiMax(1 TO Nd%) : XiMin(1) = 72## : XiMax(1) = 78## : XiMin(2) = 27## : XiMax(2) = 33## 'use this to plot STEP
detail

            REDIM StartingXiMin(1 TO Nd%), StartingXiMax(1 TO Nd%) : FOR i% = 1 TO Nd% : StartingXiMin(i%) = XiMin(i%) : StartingXiMax(i%) =
XiMax(i%) : NEXT i%

        CASE "SCHWEFEL_226"

            Nd% = 30 : Np% = 120

            REDIM XiMin(1 TO Nd%), XiMax(1 TO Nd%) : FOR i% = 1 TO Nd% : XiMin(i%) = -500## : XiMax(i%) = 500## : NEXT i%

            REDIM StartingXiMin(1 TO Nd%), StartingXiMax(1 TO Nd%) : FOR i% = 1 TO Nd% : StartingXiMin(i%) = XiMin(i%) : StartingXiMax(i%) =
XiMax(i%) : NEXT i%

        CASE "COLVILLE"

            Nd% = 4 : Np% = 16
```

```
        REDIM XiMin(1 TO Nd%), XiMax(1 TO Nd%) : FOR i% = 1 TO Nd% : XiMin(i%) = -10## : XiMax(i%) = 10## : NEXT i%

        REDIM StartingXiMin(1 TO Nd%), StartingXiMax(1 TO Nd%) : FOR i% = 1 TO Nd% : StartingXiMin(i%) = XiMin(i%) : StartingXiMax(i%) =
XiMax(i%) : NEXT i%

      CASE "GRIEWANK"

        Nd% = 2 : Np% = 8

        REDIM XiMin(1 TO Nd%), XiMax(1 TO Nd%) : FOR i% = 1 TO Nd% : XiMin(i%) = -600## : XiMax(i%) = 600## : NEXT i%

        REDIM StartingXiMin(1 TO Nd%), StartingXiMax(1 TO Nd%) : FOR i% = 1 TO Nd% : StartingXiMin(i%) = XiMin(i%) : StartingXiMax(i%) =
XiMax(i%) : NEXT i%

      CASE "HIMMELBLAU"

        Nd% = 2 : Np% = 8

        REDIM XiMin(1 TO Nd%), XiMax(1 TO Nd%) : FOR i% = 1 TO Nd% : XiMin(i%) = -6## : XiMax(i%) = 6## : NEXT i%

        REDIM StartingXiMin(1 TO Nd%), StartingXiMax(1 TO Nd%) : FOR i% = 1 TO Nd% : StartingXiMin(i%) = XiMin(i%) : StartingXiMax(i%) =
XiMax(i%) : NEXT i%

      CASE "ROSENBROCK"  '(n-D)

        Nd% = 2 : Np% = 8

        REDIM XiMin(1 TO Nd%), XiMax(1 TO Nd%) : FOR i% = 1 TO Nd% : XiMin(i%) = -2## : XiMax(i%) = 2## : NEXT i% :'XiMin(i%) = -6## : XiMax(i%)
= 6## : NEXT i%

        REDIM StartingXiMin(1 TO Nd%), StartingXiMax(1 TO Nd%) : FOR i% = 1 TO Nd% : StartingXiMin(i%) = XiMin(i%) : StartingXiMax(i%) =
XiMax(i%) : NEXT i%

      CASE "SPHERE"  '(n-D)

        Nd% = 2 : Np% = 8

        REDIM XiMin(1 TO Nd%), XiMax(1 TO Nd%) : FOR i% = 1 TO Nd% : XiMin(i%) = -100## : XiMax(i%) = 100## : NEXT i%

        REDIM StartingXiMin(1 TO Nd%), StartingXiMax(1 TO Nd%) : FOR i% = 1 TO Nd% : StartingXiMin(i%) = XiMin(i%) : StartingXiMax(i%) =
XiMax(i%) : NEXT i%

      CASE "HIMMELBLAUNLO"  '(5-D)

        Nd% = 5 : Np% = 20

        REDIM XiMin(1 TO Nd%), XiMax(1 TO Nd%)

        XiMin(1) = 78## : XiMax(1) = 102##
        XiMin(2) = 33## : XiMax(2) = 45##
        XiMin(3) = 27## : XiMax(3) = 45##
        XiMin(4) = 27## : XiMax(4) = 45##
        XiMin(5) = 27## : XiMax(5) = 45##

        REDIM StartingXiMin(1 TO Nd%), StartingXiMax(1 TO Nd%) : FOR i% = 1 TO Nd% : StartingXiMin(i%) = XiMin(i%) : StartingXiMax(i%) =
XiMax(i%) : NEXT i%

      CASE "TRIPOD"  '(2-D)

        Nd% = 2 : Np% = 8

        REDIM XiMin(1 TO Nd%), XiMax(1 TO Nd%) : FOR i% = 1 TO Nd% : XiMin(i%) = -100## : XiMax(i%) = 100## : NEXT i%

        REDIM StartingXiMin(1 TO Nd%), StartingXiMax(1 TO Nd%) : FOR i% = 1 TO Nd% : StartingXiMin(i%) = XiMin(i%) : StartingXiMax(i%) =
XiMax(i%) : NEXT i%

      CASE "ROSENBROCKF6"  '(10-D)

        Nd% = 10
        Np% = 40

        REDIM XiOffset(1 TO Nd%)

'        XiOffset(1)  =  81.0232##
'        XiOffset(2)  = -48.3950##
'        XiOffset(3)  =  19.2316##
'        XiOffset(4)  =  -2.5231##
'        XiOffset(5)  =  70.4338##
'        XiOffset(6)  =  47.1774##
'        XiOffset(7)  =  -7.8358##
'        XiOffset(8)  = -86.6693##
'        XiOffset(9)  =  57.8532##
'        XiOffset(10) = 0##

'        XiOffset(1)  =   80##
'        XiOffset(2)  =  -50##
'        XiOffset(3)  =   20##
'        XiOffset(4)  =   -3##
'        XiOffset(5)  =   70##
'        XiOffset(6)  =   47##
'        XiOffset(7)  =   -8##
'        XiOffset(8)  =  -87##
'        XiOffset(9)  =   58##
'        XiOffset(10) =    0##

        XiOffset(1)  =    5##
        XiOffset(2)  =  -25##
        XiOffset(3)  =    5##
        XiOffset(4)  =  -15##
        XiOffset(5)  =    5##
        XiOffset(6)  =  -25##
        XiOffset(7)  =   25##
        XiOffset(8)  =   -5##
        XiOffset(9)  =    5##
        XiOffset(10) =  -15##

'        XiOffset(1)  =    0##
'        XiOffset(2)  =   81.0232##
'        XiOffset(3)  =  -48.3950##
'        XiOffset(4)  =   19.2316##
'        XiOffset(5)  =   -2.5231##
'        XiOffset(6)  =   70.4338##
'        XiOffset(7)  =   47.1774##
'        XiOffset(8)  =   -7.8358##
'        XiOffset(9)  =  -86.6693##
'        XiOffset(10) =   57.8532##
```

```
        REDIM XiMin(1 TO Nd%), XiMax(1 TO Nd%) : FOR i% = 1 TO Nd% : XiMin(i%) = -100## : XiMax(i%) = 100## : NEXT i%

        REDIM StartingXiMin(1 TO Nd%), StartingXiMax(1 TO Nd%) : FOR i% = 1 TO Nd% : StartingXiMin(i%) = XiMin(i%) : StartingXiMax(i%) =
XiMax(i%) : NEXT i%

      CASE "COMPRESSIONSPRING" '(3-D)

        Nd% = 3 : Np% = 12

        REDIM XiMin(1 TO Nd%), XiMax(1 TO Nd%)

        XiMin(1) = 1##     : XiMax(1) = 70## 'integer values only!!
        XiMin(2) = 0.6##   : XiMax(2) =  3##
        XiMin(3) = 0.207## : XiMax(3) = 0.5##

        REDIM StartingXiMin(1 TO Nd%), StartingXiMax(1 TO Nd%) : FOR i% = 1 TO Nd% : StartingXiMin(i%) = XiMin(i%) : StartingXiMax(i%) =
XiMax(i%) : NEXT i%

      CASE "GEARTRAIN" '(4-D)

        Nd% = 4 : Np% = 16

        REDIM XiMin(1 TO Nd%), XiMax(1 TO Nd%) : FOR i% = 1 TO Nd% : XiMin(i%) = 12# : XiMax(i%) = 60## : NEXT i%

        REDIM StartingXiMin(1 TO Nd%), StartingXiMax(1 TO Nd%) : FOR i% = 1 TO Nd% : StartingXiMin(i%) = XiMin(i%) : StartingXiMax(i%) =
XiMax(i%) : NEXT i%

     CASE "F1" '(n-D)

        Nd% = 30 : Np% = 60

        REDIM XiMin(1 TO Nd%), XiMax(1 TO Nd%) : FOR i% = 1 TO Nd% : XiMin(i%) = -100## : XiMax(i%) = 100## : NEXT i%

        REDIM StartingXiMin(1 TO Nd%), StartingXiMax(1 TO Nd%) : FOR i% = 1 TO Nd% : StartingXiMin(i%) = XiMin(i%) : StartingXiMax(i%) =
XiMax(i%) : NEXT i%

     CASE "F2" '(n-D)

        Nd% = 30 : Np% = 60

        REDIM XiMin(1 TO Nd%), XiMax(1 TO Nd%) : FOR i% = 1 TO Nd% : XiMin(i%) = -10## : XiMax(i%) = 10# : NEXT i%

        REDIM StartingXiMin(1 TO Nd%), StartingXiMax(1 TO Nd%) : FOR i% = 1 TO Nd% : StartingXiMin(i%) = XiMin(i%) : StartingXiMax(i%) =
XiMax(i%) : NEXT i%

     CASE "F3" '(n-D)

        Nd% = 30 : Np% = 60

        REDIM XiMin(1 TO Nd%), XiMax(1 TO Nd%) : FOR i% = 1 TO Nd% : XiMin(i%) = -100## : XiMax(i%) = 100## : NEXT i%

        REDIM StartingXiMin(1 TO Nd%), StartingXiMax(1 TO Nd%) : FOR i% = 1 TO Nd% : StartingXiMin(i%) = XiMin(i%) : StartingXiMax(i%) =
XiMax(i%) : NEXT i%

     CASE "F4" '(n-D)

        Nd% = 30 : Np% = 60

        REDIM XiMin(1 TO Nd%), XiMax(1 TO Nd%) : FOR i% = 1 TO Nd% : XiMin(i%) = -100## : XiMax(i%) = 100## : NEXT i%

        REDIM StartingXiMin(1 TO Nd%), StartingXiMax(1 TO Nd%) : FOR i% = 1 TO Nd% : StartingXiMin(i%) = XiMin(i%) : StartingXiMax(i%) =
XiMax(i%) : NEXT i%

     CASE "F5" '(n-D)

        Nd% = 30 : Np% = 60

        REDIM XiMin(1 TO Nd%), XiMax(1 TO Nd%) : FOR i% = 1 TO Nd% : XiMin(i%) = -30## : XiMax(i%) = 30# : NEXT i%

        REDIM StartingXiMin(1 TO Nd%), StartingXiMax(1 TO Nd%) : FOR i% = 1 TO Nd% : StartingXiMin(i%) = XiMin(i%) : StartingXiMax(i%) =
XiMax(i%) : NEXT i%

     CASE "F6" '(n-D) STEP

        Nd% = 30 : Np% = 60

        REDIM XiMin(1 TO Nd%), XiMax(1 TO Nd%) : FOR i% = 1 TO Nd% : XiMin(i%) = -100## : XiMax(i%) = 100## : NEXT i%

        REDIM StartingXiMin(1 TO Nd%), StartingXiMax(1 TO Nd%) : FOR i% = 1 TO Nd% : StartingXiMin(i%) = XiMin(i%) : StartingXiMax(i%) =
XiMax(i%) : NEXT i%

     CASE "F7" '(n-D)

        Nd% = 30 : Np% = 60

        REDIM XiMin(1 TO Nd%), XiMax(1 TO Nd%) : FOR i% = 1 TO Nd% : XiMin(i%) = -1.28## : XiMax(i%) = 1.28## : NEXT i%

        REDIM StartingXiMin(1 TO Nd%), StartingXiMax(1 TO Nd%) : FOR i% = 1 TO Nd% : StartingXiMin(i%) = XiMin(i%) : StartingXiMax(i%) =
XiMax(i%) : NEXT i%

     CASE "F8" '(n-D)

        Nd% = 30 : Np% = 60

        REDIM XiMin(1 TO Nd%), XiMax(1 TO Nd%) : FOR i% = 1 TO Nd% : XiMin(i%) = -500## : XiMax(i%) = 500## : NEXT i%

        REDIM StartingXiMin(1 TO Nd%), StartingXiMax(1 TO Nd%) : FOR i% = 1 TO Nd% : StartingXiMin(i%) = XiMin(i%) : StartingXiMax(i%) =
XiMax(i%) : NEXT i%

     CASE "F9" '(n-D)

        Nd% = 30 : Np% = 60

        REDIM XiMin(1 TO Nd%), XiMax(1 TO Nd%) : FOR i% = 1 TO Nd% : XiMin(i%) = -5.12## : XiMax(i%) = 5.12## : NEXT i%

        REDIM StartingXiMin(1 TO Nd%), StartingXiMax(1 TO Nd%) : FOR i% = 1 TO Nd% : StartingXiMin(i%) = XiMin(i%) : StartingXiMax(i%) =
XiMax(i%) : NEXT i%

     CASE "F10" '(n-D) Ackley's Function

        Nd% = 30 : Np% = 60

        REDIM XiMin(1 TO Nd%), XiMax(1 TO Nd%) : FOR i% = 1 TO Nd% : XiMin(i%) = -32## : XiMax(i%) = 32## : NEXT i%

        REDIM StartingXiMin(1 TO Nd%), StartingXiMax(1 TO Nd%) : FOR i% = 1 TO Nd% : StartingXiMin(i%) = XiMin(i%) : StartingXiMax(i%) =
XiMax(i%) : NEXT i%
```

```
        CASE "F11" '(n-D)

            Nd% = 30 : Np% = 60

            REDIM XiMin(1 TO Nd%), XiMax(1 TO Nd%) : FOR i% = 1 TO Nd% : XiMin(i%) = -600## : XiMax(i%) = 600## : NEXT i%

            REDIM StartingXiMin(1 TO Nd%), StartingXiMax(1 TO Nd%) : FOR i% = 1 TO Nd% : StartingXiMin(i%) = XiMin(i%) : StartingXiMax(i%)
XiMax(i%) : NEXT i%

        CASE "F12" '(n-D) Penalized #1

            Nd% = 30 : Np% = 60

            REDIM XiMin(1 TO Nd%), XiMax(1 TO Nd%) : FOR i% = 1 TO Nd% : XiMin(i%) = -50## : XiMax(i%) = 50## : NEXT i%
'           REDIM XiMin(1 TO Nd%), XiMax(1 TO Nd%) : FOR i% = 1 TO Nd% : XiMin(i%) = -5## : XiMax(i%) = 5## : NEXT i% 'use this interval for second
run to improve performance
            REDIM StartingXiMin(1 TO Nd%), StartingXiMax(1 TO Nd%) : FOR i% = 1 TO Nd% : StartingXiMin(i%) = XiMin(i%) : StartingXiMax(i%)
XiMax(i%) : NEXT i%

        CASE "F13" '(n-D) Penalized #2

            Nd% = 30 : Np% = 60

            REDIM XiMin(1 TO Nd%), XiMax(1 TO Nd%) : FOR i% = 1 TO Nd% : XiMin(i%) = -50## : XiMax(i%) = 50## : NEXT i%

            REDIM StartingXiMin(1 TO Nd%), StartingXiMax(1 TO Nd%) : FOR i% = 1 TO Nd% : StartingXiMin(i%) = XiMin(i%) : StartingXiMax(i%)
XiMax(i%) : NEXT i%

        CASE "F14" '(2-D) Shekel's Foxholes

            Nd% = 2 : Np% = 8

            REDIM XiMin(1 TO Nd%), XiMax(1 TO Nd%) : FOR i% = 1 TO Nd% : XiMin(i%) = -65.536## : XiMax(i%) = 65.536## : NEXT i%

            REDIM StartingXiMin(1 TO Nd%), StartingXiMax(1 TO Nd%) : FOR i% = 1 TO Nd% : StartingXiMin(i%) = XiMin(i%) : StartingXiMax(i%)
XiMax(i%) : NEXT i%

        CASE "F15" '(4-D) Kowalik's Function

            Nd% = 4
            Np% = 16

            REDIM XiMin(1 TO Nd%), XiMax(1 TO Nd%) : FOR i% = 1 TO Nd% : XiMin(i%) = -5## : XiMax(i%) = 5## : NEXT i%

            REDIM StartingXiMin(1 TO Nd%), StartingXiMax(1 TO Nd%) : FOR i% = 1 TO Nd% : StartingXiMin(i%) = XiMin(i%) : StartingXiMax(i%)
XiMax(i%) : NEXT i%

        CASE "F16" '(2-D) Camel Back

            Nd% = 2 : Np% = 8

            REDIM XiMin(1 TO Nd%), XiMax(1 TO Nd%) : FOR i% = 1 TO Nd% : XiMin(i%) = -5## : XiMax(i%) = 5## : NEXT i%

            REDIM StartingXiMin(1 TO Nd%), StartingXiMax(1 TO Nd%) : FOR i% = 1 TO Nd% : StartingXiMin(i%) = XiMin(i%) : StartingXiMax(i%)
XiMax(i%) : NEXT i%

        CASE "F17" '(2-D) Branin

            Nd% = 2 : Np% = 8

            REDIM XiMin(1 TO Nd%), XiMax(1 TO Nd%) : XiMin(1) = -5## : XiMax(1) = 10## : XiMin(2) = 0## : XiMax(2) = 15##

            REDIM StartingXiMin(1 TO Nd%), StartingXiMax(1 TO Nd%) : FOR i% = 1 TO Nd% : StartingXiMin(i%) = XiMin(i%) : StartingXiMax(i%)
XiMax(i%) : NEXT i%

        CASE "F18" '(2-D) Goldstein-Price

            Nd% = 2 : Np% = 8

            REDIM XiMin(1 TO Nd%), XiMax(1 TO Nd%) : XiMin(1) = -2## : XiMax(1) = 2## : XiMin(2) = -2## : XiMax(2) = 2##

            REDIM StartingXiMin(1 TO Nd%), StartingXiMax(1 TO Nd%) : FOR i% = 1 TO Nd% : StartingXiMin(i%) = XiMin(i%) : StartingXiMax(i%)
XiMax(i%) : NEXT i%

        CASE "F19" '(3-D) Hartman's Family #1

            Nd% = 3 : Np% = 12

            REDIM XiMin(1 TO Nd%), XiMax(1 TO Nd%) : FOR i% = 1 TO Nd% : XiMin(i%) = 0## : XiMax(i%) = 1## : NEXT i%

            REDIM StartingXiMin(1 TO Nd%), StartingXiMax(1 TO Nd%) : FOR i% = 1 TO Nd% : StartingXiMin(i%) = XiMin(i%) : StartingXiMax(i%)
XiMax(i%) : NEXT i%

        CASE "F20" '(6-D) Hartman's Family #2

            Nd% = 6 : Np% = 24

            REDIM XiMin(1 TO Nd%), XiMax(1 TO Nd%) : FOR i% = 1 TO Nd% : XiMin(i%) = 0## : XiMax(i%) = 1## : NEXT i%

            REDIM StartingXiMin(1 TO Nd%), StartingXiMax(1 TO Nd%) : FOR i% = 1 TO Nd% : StartingXiMin(i%) = XiMin(i%) : StartingXiMax(i%)
XiMax(i%) : NEXT i%

        CASE "F21" '(4-D) Shekel's Family m=5

            Nd% = 4 : Np% = 16

            REDIM XiMin(1 TO Nd%), XiMax(1 TO Nd%) : FOR i% = 1 TO Nd% : XiMin(i%) = 0## : XiMax(i%) = 10## : NEXT i%

            REDIM StartingXiMin(1 TO Nd%), StartingXiMax(1 TO Nd%) : FOR i% = 1 TO Nd% : StartingXiMin(i%) = XiMin(i%) : StartingXiMax(i%)
XiMax(i%) : NEXT i%

        CASE "F22" '(4-D) Shekel's Family m=7

            Nd% = 4 : Np% = 16

            REDIM XiMin(1 TO Nd%), XiMax(1 TO Nd%) : FOR i% = 1 TO Nd% : XiMin(i%) = 0## : XiMax(i%) = 10## : NEXT i%

            REDIM StartingXiMin(1 TO Nd%), StartingXiMax(1 TO Nd%) : FOR i% = 1 TO Nd% : StartingXiMin(i%) = XiMin(i%) : StartingXiMax(i%)
XiMax(i%) : NEXT i%

        CASE "F23" '(4-D) Shekel's Family m=10

            Nd% = 4 : Np% = 16
```

```
            REDIM XiMin(1 TO Nd%), XiMax(1 TO Nd%) : FOR i% = 1 TO Nd% : XiMin(i%) = 0## : XiMax(i%) = 10## : NEXT i%

            REDIM StartingXiMin(1 TO Nd%), StartingXiMax(1 TO Nd%) : FOR i% = 1 TO Nd% : StartingXiMin(i%) = XiMin(i%) : StartingXiMax(i%) =
XiMax(i%) : NEXT i%

        CASE "PBM_1" '2-D

            Nd%                     = 2
            NumProbesPerDimension%  = 2 '4 '20
            Np%                     = NumProbesPerDimension%*Nd%

            Nt&     = 100
            G       = 2##
            Alpha   = 2##
            Beta    = 2##
            DeltaT  = 1##
            Frep    = 0.5##

            PlaceInitialProbes$   = "UNIFORM ON-AXIS"
            InitialAcceleration$  = "ZERO"
            RepositionFactor$     = "VARIABLE" '"FIXED"

            Np% = NumProbesPerDimension%*Nd%

            REDIM XiMin(1 TO Nd%), XiMax(1 TO Nd%)

            XiMin(1) = 0.5## : XiMax(1) = 3## 'dipole length, L, in Wavelengths
            XiMin(2) = 0## : XiMax(2) = Pi2 'polar angle, Theta, in Radians

            REDIM StartingXiMin(1 TO Nd%), StartingXiMax(1 TO Nd%) : FOR i% = 1 TO Nd% : StartingXiMin(i%) = XiMin(i%) : StartingXiMax(i%) =
XiMax(i%) : NEXT i%

            NN% = FREEFILE : OPEN "INFILE.DAT" FOR OUTPUT AS #NN% : PRINT #NN%,"PBM1.NEC" : PRINT #NN%,"PBM1.OUT" : CLOSE #NN% 'NEC Input/Output
Files

        CASE "PBM_2" '2-D

            AddNoiseToPBM2$ = "NO" '"YES" '"NO" '"YES"

            Nd%                     = 2
            NumProbesPerDimension% = 4 '20
            Np%                     = NumProbesPerDimension%*Nd%

            Nt&     = 100
            G       = 2##
            Alpha   = 2##
            Beta    = 2##
            DeltaT  = 1##
            Frep    = 0.5##

            PlaceInitialProbes$   = "UNIFORM ON-AXIS"
            InitialAcceleration$  = "ZERO"
            RepositionFactor$     = "VARIABLE" '"FIXED"

            Np% = NumProbesPerDimension%*Nd%

            REDIM XiMin(1 TO Nd%), XiMax(1 TO Nd%)

            XiMin(1) = 5## : XiMax(1) = 15## 'dipole separation, D, in Wavelengths
            XiMin(2) = 0## : XiMax(2) = Pi   'polar angle, Theta, in Radians

            REDIM StartingXiMin(1 TO Nd%), StartingXiMax(1 TO Nd%) : FOR i% = 1 TO Nd% : StartingXiMin(i%) = XiMin(i%) : StartingXiMax(i%) =
XiMax(i%) : NEXT i%

            NN% = FREEFILE : OPEN "INFILE.DAT" FOR OUTPUT AS #NN% : PRINT #NN%,"PBM2.NEC" : PRINT #NN%,"PBM2.OUT" : CLOSE #NN%

        CASE "PBM_3" '2-D

            Nd%                     = 2
            NumProbesPerDimension% = 4 '20
            Np%                     = NumProbesPerDimension%*Nd%

            Nt&     = 100
            G       = 2##
            Alpha   = 2##
            Beta    = 2##
            DeltaT  = 1##
            Frep    = 0.5##

            PlaceInitialProbes$   = "UNIFORM ON-AXIS"
            InitialAcceleration$  = "ZERO"
            RepositionFactor$     = "VARIABLE" '"FIXED"

            Np% = NumProbesPerDimension%*Nd%

            REDIM XiMin(1 TO Nd%), XiMax(1 TO Nd%)

            XiMin(1) = 0## : XiMax(1) = 4## 'Phase Parameter, Beta (0-4)
            XiMin(2) = 0## : XiMax(2) = Pi  'polar angle, Theta, in Radians

            REDIM StartingXiMin(1 TO Nd%), StartingXiMax(1 TO Nd%) : FOR i% = 1 TO Nd% : StartingXiMin(i%) = XiMin(i%) : StartingXiMax(i%) =
XiMax(i%) : NEXT i%

            NN% = FREEFILE : OPEN "INFILE.DAT" FOR OUTPUT AS #NN% : PRINT #NN%,"PBM3.NEC" : PRINT #NN%,"PBM3.OUT" : CLOSE #NN%

        CASE "PBM_4" '2-D

            Nd%                     = 2
            NumProbesPerDimension% = 4 '6 '2 '4 '20
            Np%                     = NumProbesPerDimension%*Nd%

            Nt&     = 100
            G       = 2##
            Alpha   = 2##
            Beta    = 2##
            DeltaT  = 1##
            Frep    = 0.5##

            PlaceInitialProbes$   = "UNIFORM ON-AXIS"
            InitialAcceleration$  = "ZERO"
            RepositionFactor$     = "VARIABLE" '"FIXED"

            Np% = NumProbesPerDimension%*Nd%

            REDIM XiMin(1 TO Nd%), XiMax(1 TO Nd%)

            XiMin(1) = 0.5##   : XiMax(1) = 1.5##  'ARM LENGTH (NOT Total Length), wavelengths (0.5-1.5)
```

```
            XiMin(2) = Pi/18## : XiMax(2) = Pi/2## 'Inner angle, Alpha, in Radians (Pi/18-Pi/2)

            REDIM StartingXiMin(1 TO Nd%), StartingXiMax(1 TO Nd%) : FOR i% = 1 TO Nd% : StartingXiMin(i%) = XiMin(i%) : StartingXiMax(i%) =
XiMax(i%) : NEXT i%

            NN% = FREEFILE : OPEN "INFILE.DAT" FOR OUTPUT AS #NN% : PRINT #NN%,"PBM4.NEC" : PRINT #NN%,"PBM4.OUT" : CLOSE #NN%

        CASE "PBM_5"

            NumCollinearElements% =  6 '30 'EVEN or ODD: 6,7,10,13,16,24 used by PBM

            Nd%                   = NumCollinearElements% - 1
            NumProbesPerDimension% = 4 '20
            Np%                   = NumProbesPerDimension%*Nd%

            Nt&     = 100
            G       = 2##
            Alpha   = 2##
            Beta    = 2##
            DeltaT  = 1##
            Frep    = 0.5##

            PlaceInitialProbes$   = "UNIFORM ON-AXIS"
            InitialAcceleration$  = "ZERO"
            RepositionFactor$     = "VARIABLE" '"FIXED"

            Nd% = NumCollinearElements% - 1

            Np% = NumProbesPerDimension%*Nd%

            REDIM XiMin(1 TO Nd%), XiMax(1 TO Nd%) : FOR i% = 1 TO Nd% : XiMin(i%) = 0.5## : XiMax(i%) = 1.5## : NEXT i%

            REDIM StartingXiMin(1 TO Nd%), StartingXiMax(1 TO Nd%) : FOR i% = 1 TO Nd% : StartingXiMin(i%) = XiMin(i%) : StartingXiMax(i%) =
XiMax(i%) : NEXT i%

            NN% = FREEFILE : OPEN "INFILE.DAT" FOR OUTPUT AS #NN% : PRINT #NN%,"PBM4.NEC" : PRINT #NN%,"PBM5.OUT" : CLOSE #NN%

        CASE "LD_MONO" 'Loaded Monopole

            IF DIR$("NEC2D_200_02-22-2011.EXE") = "" THEN
                MSGBOX("WARNING! NEC2D executable NEC2D_200_02-22-2011.EXE not found!"+CHR$(13)+CHR$(13)+_
                "                RUN TERMINATED."+CHR$(13)+CHR$(13)+_
                "           Output files are LOADED1.DAT, LOADED1.NEC")
                EXIT SUB
            END IF

            Nd% = 2 : Np% = 4

            REDIM XiMin(1 TO Nd%), XiMax(1 TO Nd%)

            XiMin(1) = 0##    : XiMax(1) = 1000## 'load resistor value

            XiMin(2) = 0.05## : XiMax(2) = 10.65## 'height of resistor along monopole

            REDIM StartingXiMin(1 TO Nd%), StartingXiMax(1 TO Nd%) : FOR i% = 1 TO Nd% : StartingXiMin(i%) = XiMin(i%) : StartingXiMax(i%) =
XiMax(i%) : NEXT i%

            NN% = FREEFILE : OPEN "INFILE.DAT" FOR OUTPUT AS #NN% : PRINT #NN%,"LD_MONO.NEC" : PRINT #NN%,"LD_MONO.OUT" : CLOSE #NN%
' ==========================================================================================================================
'  NOTE - DON'T FORGET TO ADD NEW TEST FUNCTIONS TO FUNCTION ObjectiveFunction() ABOVE !!
' ==========================================================================================================================

    END SELECT

    IF Nd% > 100 THEN Nt& = MIN(Nt&,200) 'to avoid array dimensioning problems

    DiagLength = 0## : FOR i% = 1 TO Nd% : DiagLength = DiagLength + (XiMax(i%)-XiMin(i%))^2 : NEXT i% : DiagLength = SQR(DiagLength) 'compute length
of decision space principal diagonal
END SUB 'GetFunctionRunParameters()
'--------------------------------

FUNCTION ParrottF4(R(),Nd%,p%,j&) 'Parrott F4 (1-D)

'MAXIMUM = 1 AT ~0.0796875... WITH ZERO OFFSET (SEEMS TO WORK BEST WITH JUST 3 PROBES, BUT NOT ALLOWED IN THIS VERSION...)

'References:

'Beasley, D., D. R. Bull, and R. R. Martin, "A Sequential Niche Technique for Multimodal
'Function Optimization," Evol. Comp. (MIT Press), vol. 1, no. 2, 1993, pp. 101-125
'(online at http://citeseer.ist.psu.edu/beasley93sequential.html).

'Parrott, D., and X. Li, "Locating and Tracking Multiple Dynamic Optima on a Particle Swarm
'Model Using Speciation," IEEE Trans. Evol. Computation, vol. 10, no. 4, Aug. 2006, pp. 440-458.

LOCAL Z, x, offset AS EXT

    offset = 0##

    x = R(p%,1,j&)

    Z = EXP(-2##*LOG(2##)*((x-0.08##-offset)/0.854##)^2)*(SIN(5##*Pi*((x-offset)^0.75##-0.05##)))^6 'WARNING! This is a NATURAL LOG, NOT Log10!!!

    ParrottF4 = Z
END FUNCTION 'ParrottF4()
'----------------------------

FUNCTION SGO(R(),Nd%,p%,j&) 'SGO Function (2-D)

'MAXIMUM = ~130.8323226... @ ~(-2.8362075...,-2.8362075...) WITH ZERO OFFSET.

'Reference:

'Hsiao, Y., Chuang, C., Jiang, J., and Chien, C., "A Novel Optimization Algorithm: Space
'Gravitational Optimization," Proc. of 2005 IEEE International Conference on Systems, Man,
'and Cybernetics, 3, 2323-2328. (2005)

    LOCAL x1, x2, Z, t1, t2, SGOx1offset, SGOx2offset AS EXT

    SGOx1offset = 0## : SGOx2offset = 0##

'    SGOx1offset = 40## : SGOx2offset = 10##
```

```
    x1 = R(p%,1,j&) - SGOx1offset : x2 = R(p%,2,j&) - SGOx2offset

    t1 = x1^4 - 16##*x1^2 + 0.5##*x1 : t2 = x2^4 - 16##*x2^2 + 0.5##*x2

    Z = t1 + t2

    SGO = -Z

END FUNCTION 'SGO()
'------------------

FUNCTION GoldsteinPrice(R(),Nd%,p%,j&) 'Goldstein-Price Function (2-D)

'MAXIMUM = -3 @ (0,-1) WITH ZERO OFFSET.

'Reference:

'Cui, Z., Zeng, J., and Sun, G. (2006) 'A Fast Particle Swarm Optimization,' Int'l. J.
'Innovative Computing, Information and Control, vol. 2, no. 6, December, pp. 1365-1380.

    LOCAL Z, x1, x2, offset1, offset2, t1, t2 AS EXT

    offset1 = 0## : offset2 = 0##

'   offset1 = 20## : offset2 = -10##

    x1 = R(p%,1,j&)-offset1 : x2 = R(p%,2,j&)-offset2

    t1 = 1##+(x1+x2+1##)^2*(19##-14##*x1+3##*x1^2-14##*x2+6##*x1*x2+3##*x2^2)

    t2 = 30##+(2##*x1-3##*x2)^2*(18##-32##*x1+12##*x1^2+48##*x2-36##*x1*x2+27##*x2^2)

    Z  = t1*t2

    GoldsteinPrice = -Z

END FUNCTION 'GoldsteinPrice()
'-----------

FUNCTION StepFunction(R(),Nd%,p%,j&) 'Step Function (n-D)

'MAXIMUM VALUE = 0 @ [Offset]^n.

'Reference:

'Yao, X., Liu, Y., and Lin, G., "Evolutionary Programming Made Faster,"
'IEEE Trans. Evolutionary Computation, Vol. 3, No. 2, 82-102, Jul. 1999.

    LOCAL Offset, Z AS EXT

    LOCAL i%

    Z = 0## : Offset = 0## '75.123## '0##

    FOR i% = 1 TO Nd%

        IF Nd% = 2 AND i% = 1 THEN Offset = 75 '75##

        IF Nd% = 2 AND i% = 2 THEN Offset = 35 '30 '35##

        Z = Z + INT((R(p%,i%,j&)-Offset) + 0.5##)^2

    NEXT i%

    StepFunction = -Z

END FUNCTION 'StepFunction()
'-----------

FUNCTION Schwefel226(R(),Nd%,p%,j&) 'Schwefel Problem 2.26 (n-D)

'MAXIMUM = 12,569.5 @ [420.8687]^30 (30-D CASE).

'Reference:

'Yao, X., Liu, Y., and Lin, G., "Evolutionary Programming Made Faster,"
'IEEE Trans. Evolutionary Computation, Vol. 3, No. 2, 82-102, Jul. 1999.

    LOCAL Z, Xi AS EXT

    LOCAL i%

    Z = 0##

    FOR i% = 1 TO Nd%

        Xi = R(p%,i%,j&)

        Z = Z + Xi*SIN(SQR(ABS(Xi)))

    NEXT i%

    Schwefel226 = Z

END FUNCTION 'SCHWEFEL226()
'-----------

FUNCTION Colville(R(),Nd%,p%,j&) 'Colville Function (4-D)

'MAXIMUM = 0 @ (1,1,1,1) WITH ZERO OFFSET.

'Reference:

'Doo-Hyun, and Se-Young, O., "A New Mutation Rule for Evolutionary Programming Motivated from
'Backpropagation Learning," IEEE Trans. Evolutionary Computation, Vol. 4, No. 2, pp. 188-190,
'July 2000.

    LOCAL Z, x1, x2, x3, x4, offset AS EXT

    offset = 0## '7.123##
```

```
    x1 = R(p%,1,j&)-offset : x2 = R(p%,2,j&)-offset : x3 = R(p%,3,j&)-offset : x4 = R(p%,4,j&)-offset

    Z =   100##*(x2-x1^2)^2 + (1##-x1)^2  + _
           90##*(x4-x3^2)^2 + (1##-x3)^2  + _
          10.1##*((x2-1##)^2 + (x4-1##)^2) + _
          19.8##*(x2-1##)*(x4-1##)

    Colville = -Z

END FUNCTION 'Colville()

'-----------

FUNCTION Griewank(R(),Nd%,p%,j&) 'Griewank (n-D)

'Max of zero at (0,...,0)

'Reference: Yao, X., Liu, Y., and Lin, G., "Evolutionary Programming Made Faster,"
'IEEE Trans. Evolutionary Computation, Vol. 3, No. 2, 82-102, Jul. 1999.

    LOCAL Offset, Sum, Prod, Z, Xi AS EXT

    LOCAL i%

    Sum = 0## : Prod = 1##

    Offset = 75.123##

    FOR i% = 1 TO Nd%

        Xi = R(p%,i%,j&) - offset

        Sum = Sum + Xi^2

        Prod = Prod*COS(Xi/SQR(i%))

    NEXT i%

    Z = Sum/4000## - Prod + 1##

    Griewank = -Z

END FUNCTION 'Griewank()

'-----------

FUNCTION Himmelblau(R(),Nd%,p%,j&) 'Himmelblau (2-D)

    LOCAL Z, x1, x2, offset AS EXT

    offset = 0##

    x1 = R(p%,1,j&)-offset : x2 = R(p%,2,j&)-offset

    Z = 200## - (x1^2 + x2 -11##)^2 - (x1+x2^2-7##)^2

    Himmelblau = Z

END FUNCTION 'Himmelblau()

'-----------

FUNCTION Rosenbrock(R(),Nd%,p%,j&) 'Rosenbrock (n-D)

'MAXIMUM = 0 @ [1,...,1]^n (n-D CASE).

'Reference: Yao, X., Liu, Y., and Lin, G., "Evolutionary Programming Made Faster,"
'IEEE Trans. Evolutionary Computation, Vol. 3, No. 2, 82-102, Jul. 1999.

    LOCAL Z, Xi, Xi1 AS EXT

    LOCAL i%

    Z = 0##

    FOR i% = 1 TO Nd%-1

        Xi  = R(p%,i%,j&) : Xi1 = R(p%,i%+1,j&)

        Z = Z + 100##*(Xi1-Xi^2)^2 + (Xi-1##)^2

    NEXT i%

    Rosenbrock = -Z

END FUNCTION 'ROSENBROCK()

'-----------

FUNCTION Sphere(R(),Nd%,p%,j&) 'Sphere (n-D)

'MAXIMUM = 0 @ [0,...,0]^n (n-D CASE).

'Reference: Yao, X., Liu, Y., and Lin, G., "Evolutionary Programming Made Faster,"
'IEEE Trans. Evolutionary Computation, Vol. 3, No. 2, 82-102, Jul. 1999.

    LOCAL Z, Xi, Xi1 AS EXT

    LOCAL i%

    Z = 0##

    FOR i% = 1 TO Nd%

        Xi  = R(p%,i%,j&)

        Z = Z + Xi^2

    NEXT i%

    Sphere = -Z

END FUNCTION 'SPHERE()
```

```
'-----------

FUNCTION HimmelblauNLO(R(),Nd%,p%,j&) 'Himmelblau non-linear optimization (5-D)

'MAXIMUM ~ 31025.5562644972 @ (78.0,33.0,27.0709971052,45.0,44.9692425501)

'Reference: "Constrained Optimization using CODEQ," Mahamed G.H. Omran & Ayed Salman,
'Chaos, Solitons and Tractals, 42(2009), 662-668

    LOCAL Z, x1, x2, x3, x4, x5, g1, g2, g3 AS EXT

    Z = 1E4200

    x1 = R(p%,1,j&) : x2 = R(p%,2,j&) : x3 = R(p%,3,j&) : x4 = R(p%,4,j&) : x5 = R(p%,5,j&)

    g1 = 85.334407## + 0.0056858##*x2*x5 + 0.00026##*x1*x4  - 0.0022053##*x3*x5

    g2 = 80.51249##  + 0.0071317##*x2*x5 + 0.0029955##*x1*x2 + 0.0021813##*x3*x3

    g3 = 9.300961##  + 0.0047026##*x3*x5 + 0.0012547##*x1*x3 + 0.0019085##*x3*x4

    IF g1## < 0 OR g1 > 92## OR g2 < 90## OR g2 > 110## OR g3 < 20## OR g3 > 25## THEN GOTO ExitHimmelblauNLO

    Z = 5.3578547##*x3*x3 + 0.8356891##*x1*x5 + 37.29329##*x1 - 40792.141##

ExitHimmelblauNLO:

    HimmelblauNLO = -Z

END FUNCTION 'HimmelblauNLO()

'-----------

FUNCTION Tripod(R(),Nd%,p%,j&) 'Tripod (2-D)

'MAXIMUM = 0 at (0,-50)

'Reference: "Appendix: A mini-benchmark," Maurice Clerc

    LOCAL Z, x1, x2, s1, s2, t1, t2, t3 AS EXT

    x1 = R(p%,1,j&) : x2 = R(p%,2,j&)

    s1 = Sign(x1) : s2 = Sign(x2)

    t1 = (1##-s2)*(ABS(x1)+ABS(x2+50##))

    t2 = 0.5##*(1##+s2)*(1##-s1)*(1##+ABS(x1+50##)+ABS(x2-50##))

    t3 = (1##+s1)*(2##+ABS(x1-50##)+ABS(x2-50##))

    Z = 0.5##*(t1 + t2 + t3)

    Tripod = -Z

END FUNCTION 'Tripod()

'--------------------

FUNCTION Sign(X)

LOCAL Z AS EXT

    Z = 1## : IF X =< 0## THEN Z = -1##

    Sign = Z

END FUNCTION

'-----------

FUNCTION RosenbrockF6(R(),Nd%,p%,j&) 'Rosenbrock F6 (10-D)

'WARNING !!  03-19-10  THIS FUNCTION CONTAINS ERRORS.  SEE CLERC's EMAIL!

'MAXIMUM = 394 at (0,-50)

'Reference: "Appendix: A mini-benchmark," Maurice Clerc (NOTE: Uses his notation...)

    LOCAL Z, Xi, Xi1, Zi, Zi1, Sum AS EXT

    LOCAL i%

    Sum = 0##

    FOR i% = 2 TO Nd%

        Xi = R(p%,i%,j&) : Xi1 = R(p%,i%-1,j&)

        Zi = Xi - Xioffset(i%) + 1## : Zi1 = Xi1 - Xioffset(i%-1) + 1##

        Sum = Sum + 100##*(Zi1^2-Zi)^2  + (Zi1-1##)^2

    NEXT i%

    Z = 390## + Sum

    RosenbrockF6 = -Z

END FUNCTION 'RosenbrockF6()

'--------------------

FUNCTION CompressionSpring(R(),Nd%,p%,j&) 'Compression Spring (3-D)

'MAXIMUM = 394 at (0,-50)

'Reference: "Appendix: A mini-benchmark," Maurice Clerc (NOTE: Uses his notation...)

LOCAL Z, x1, x2, x3, g1, g2, g3, g4, g5, Cf, Fmax, S, Lf, Lmax, SigP, SigPM, Fp, K, SigW AS EXT

    Z = 1E4200

    x1 = ROUND(R(p%,1,j&),0) : x2 = R(p%,2,j&)  : x3 = ROUND(R(p%,3,j&),3)

    Cf = 1## + 0.75##*x3/(x2-x3)+0.615##*x3/x2
```

```
    Fmax = 1000## : S = 189000## : Lmax = 14## : SigPM = 6## : Fp = 300## : SigW = 1.25##

    K = 11.5##*1E6*x3^4/(8##*x1*x2^3)

    Lf = Fmax/K + 1.05##*(x1+2##)*x3

    SigP = Fp/K

    g1 = 8##*Cf*Fmax*x2/(Pi*x3^3) - S

    g2 = Lf - Lmax

    g3 = SigP - SigPM

    g4 = SigP -Fp/K  'WARNING!  03-19-10.  THIS IS SATISFIED EXACTLY (SEE CLERC's EMAIL - TYPO IN HIS BENCHMARKS)

    g5 = SigW - (Fmax-Fp)/K

    IF g1 > 0## OR g2 > 0## OR g3 > 0## OR g4 > 0## OR g5 > 0## THEN GOTO ExitCompressionSpring

    Z = Pi^2*x2*x3^2*(x1+1##)/4##

ExitCompressionSpring:

    CompressionSpring = -Z

END FUNCTION 'CompressionSpring
'--------------------

FUNCTION GearTrain(R(),Nd%,p%,j&) 'GearTrain (4-D)
'MAXIMUM = 394 at (0,-50)

'Reference: "Appendix: A mini-benchmark," Maurice Clerc (NOTE: Uses his notation...)

LOCAL Z, x1, x2, x3, x4 AS EXT

    x1 = ROUND(R(p%,1,j&),0) : x2 = ROUND(R(p%,2,j&),0)
    x3 = ROUND(R(p%,3,j&),0) : x4 = ROUND(R(p%,4,j&),0)

    Z = (1##/6.931##-x1*x2/(x3*x4))^2

    GearTrain = -Z

END FUNCTION 'GearTrain
'--------------------

FUNCTION F1(R(),Nd%,p%,j&) 'F1 (n-D)
'MAXIMUM = ZERO (n-D CASE).
'Reference:

    LOCAL Z, Xi AS EXT

    LOCAL i%

    Z = 0##

    FOR i% = 1 TO Nd%

        Xi = R(p%,i%,j&)

        Z = Z + Xi^2

    NEXT i%

    F1 = -Z

END FUNCTION 'F1
'-----------

FUNCTION F2(R(),Nd%,p%,j&) 'F2 (n-D)
'MAXIMUM = ZERO (n-D CASE).
'Reference:

    LOCAL Sum, prod, Z, Xi AS EXT

    LOCAL i%

    Z = 0## : Sum = 0## : Prod = 1##

    FOR i% = 1 TO Nd%

        Xi = R(p%,i%,j&)

        Sum  = Sum+ ABS(Xi)

        Prod = Prod*ABS(Xi)

    NEXT i%

    Z = Sum + Prod

    F2 = -Z

END FUNCTION 'F2
'-----------

FUNCTION F3(R(),Nd%,p%,j&) 'F3 (n-D)
'MAXIMUM = ZERO (n-D CASE).
'Reference:

    LOCAL Z, Xk, Sum AS EXT

    LOCAL i%, k%

    Z = 0##
```

```
        FOR i% = 1 TO Nd%
            Sum = 0##
            FOR k% = 1 TO i%
                Xk = R(p%,k%,j&)
                Sum = Sum + Xk
            NEXT k%
            Z = Z + Sum^2
        NEXT i%
        F3 = -Z
END FUNCTION 'F3

'-----------

FUNCTION F4(R(),Nd%,p%,j&) 'F4 (n-D)
'MAXIMUM = ZERO (n-D CASE).
'Reference:
        LOCAL Z, Xi, MaxXi AS EXT
        LOCAL i%
        MaxXi = -1E4200
        FOR i% = 1 TO Nd%
            Xi = R(p%,i%,j&)
            IF ABS(Xi) >= MaxXi THEN MaxXi = ABS(Xi)
        NEXT i%
        F4 = -MaxXi
END FUNCTION 'F4
'-----------
FUNCTION F5(R(),Nd%,p%,j&) 'F5 (n-D)
'MAXIMUM = ZERO (n-D CASE).
'Reference:
        LOCAL Z, Xi, XiPlus1 AS EXT
        LOCAL i%
        Z = 0##
        FOR i% = 1 TO Nd%-1
            Xi      = R(p%,i%,j&)
            XiPlus1 = R(p%,i%+1,j&)
            Z = Z + (100##*(XiPlus1-Xi^2)^2+(Xi-1##))^2
        NEXT i%
        F5 = -Z
END FUNCTION 'F5
'-----------
FUNCTION F6(R(),Nd%,p%,j&) 'F6 (n-D STEP)
'MAXIMUM VALUE = 0 @ [offset]^n.
'Reference:
'Yao, X., Liu, Y., and Lin, G., "Evolutionary Programming Made Faster,"
'IEEE Trans. Evolutionary Computation, Vol. 3, No. 2, 82-102, Jul. 1999.

        LOCAL Z AS EXT
        LOCAL i%
        Z = 0##
        FOR i% = 1 TO Nd%
            Z = Z + INT(R(p%,i%,j&) + 0.5##)^2
        NEXT i%
        F6 = -Z
END FUNCTION 'F6
'-----------
FUNCTION F7(R(),Nd%,p%,j&) 'F7
'MAXIMUM VALUE = 0 @ [offset]^n.
'Reference:
'Yao, X., Liu, Y., and Lin, G., "Evolutionary Programming Made Faster,"
'IEEE Trans. Evolutionary Computation, Vol. 3, No. 2, 82-102, Jul. 1999.

        LOCAL Z, Xi AS EXT
```

```
    LOCAL i%

    Z = 0##

    FOR i% = 1 TO Nd%

        Xi = R(p%,i%,j&)

        Z = Z + i%*Xi^4

    NEXT i%

    F7 = -Z - RandomNum(0##,1##)

END FUNCTION 'F7

'-----------

FUNCTION F8(R(),Nd%,p%,j&) '(n-D) F8 [Schwefel Problem 2.26]

'MAXIMUM = 12,569.5 @ [420.8687]^30 (30-D CASE).

'Reference:

'Yao, X., Liu, Y., and Lin, G., "Evolutionary Programming Made Faster,"
'IEEE Trans. Evolutionary Computation, Vol. 3, No. 2, 82-102, Jul. 1999.

    LOCAL Z, Xi AS EXT

    LOCAL i%

    Z = 0##

    FOR i% = 1 TO Nd%

        Xi = R(p%,i%,j&)

        Z  = Z - Xi*SIN(SQR(ABS(Xi)))

    NEXT i%

    F8 = -Z

END FUNCTION 'F8

'-----------

FUNCTION F9(R(),Nd%,p%,j&) '(n-D) F9 [Rastrigin]

'MAXIMUM = ZERO (n-D CASE).

'Reference:

'Yao, X., Liu, Y., and Lin, G., "Evolutionary Programming Made Faster,"
'IEEE Trans. Evolutionary Computation, Vol. 3, No. 2, 82-102, Jul. 1999.

    LOCAL Z, Xi AS EXT

    LOCAL i%

    Z = 0##

    FOR i% = 1 TO Nd%

        Xi = R(p%,i%,j&)

        Z  = Z + (Xi^2 - 10##*COS(TwoPi*Xi) + 10##)^2

    NEXT i%

    F9 = -Z

END FUNCTION 'F9

'-----------

FUNCTION F10(R(),Nd%,p%,j&) '(n-D) F10 [Ackley's Function]

'MAXIMUM = ZERO (n-D CASE).

'Reference:

'Yao, X., Liu, Y., and Lin, G., "Evolutionary Programming Made Faster,"
'IEEE Trans. Evolutionary Computation, Vol. 3, No. 2, 82-102, Jul. 1999.

    LOCAL Z, Xi, Sum1, Sum2 AS EXT

    LOCAL i%

    Z = 0## : Sum1 = 0## : Sum2 = 0##

    FOR i% = 1 TO Nd%

        Xi   = R(p%,i%,j&)

        Sum1 = Sum1 + Xi^2

        Sum2 = Sum2 + COS(TwoPi*Xi)

    NEXT i%

    Z = -20##*EXP(-0.2##*SQR(Sum1/Nd%)) - EXP(Sum2/Nd%) + 20## + e

    F10 = -Z

END FUNCTION 'F10

'-----------

FUNCTION F11(R(),Nd%,p%,j&) '(n-D) F11

'MAXIMUM = ZERO (n-D CASE).

'Reference:

'Yao, X., Liu, Y., and Lin, G., "Evolutionary Programming Made Faster,"
'IEEE Trans. Evolutionary Computation, Vol. 3, No. 2, 82-102, Jul. 1999.
```

```
      LOCAL Z, Xi, Sum, Prod AS EXT

      LOCAL i%

      Z = 0## : Sum = 0## : Prod = 1##

      FOR i% = 1 TO Nd%

            Xi   = R(p%,i%,j&)

            Sum  = Sum + (Xi-100##)^2

            Prod = Prod*COS((Xi-100##)/SQR(i%))

      NEXT i%

      Z = Sum/4000## - Prod + 1##

      F11 = -Z

END FUNCTION 'F11

'-----

FUNCTION u(Xi,a,k,m)

LOCAL Z AS EXT

      Z = 0##

      SELECT CASE Xi

            CASE > a  : Z = k*(Xi-a)^m

            CASE < -a : Z = k*(-Xi-a)^m

      END SELECT

      u = Z

END FUNCTION

'-----------

FUNCTION F12(R(),Nd%,p%,j&)  '(n-D) F12, Penalized #1

'Ref: Yao(1999).  Max=0 @ (-1,-1,...,-1), -50=<xi=<50.

      LOCAL Offset, Sum1, Sum2, Z, X1, Y1, Xn, Yn, Xi, Yi, XiPlus1, YiPlus1 AS EXT

      LOCAL i%, m%, A$

      X1 = R(p%,1,j&)    : Y1 = 1## + (X1+1##)/4##

      Xn = R(p%,Nd%,j&) : Yn = 1## + (Xn+1##)/4##

      Sum1 = 0##

      FOR i% = 1 TO Nd%-1

            Xi      = R(p%,i%,j&)   : Yi      = 1## + (Xi+1##)/4##

            XiPlus1 = R(p%,i%+1,j&): YiPlus1 = 1## + (Xi+1##+1##)/4##

            Sum1 = Sum1 + (Yi-1##)^2*(1##+10##*(SIN(Pi*YiPlus1))^2)

      NEXT i%

      Sum1 = Sum1 + 10##*(SIN(Pi*Y1))^2 + (Yn-1##)^2

      Sum1 = Pi*Sum1/Nd%

      Sum2 = 0##

      FOR i% = 1 TO Nd%

            Xi = R(p%,i%,j&)

            Sum2 = Sum2 + u(Xi,10##,100##,4##)

      NEXT i%

      Z = Sum1 + Sum2

      F12 = -Z

END FUNCTION 'F12()

'------------------

FUNCTION F13(R(),Nd%,p%,j&)  '(n-D) F13, Penalized #2

'Ref: Yao(1999).  Max=0 @ (1,1,...,1), -50=<xi=<50.

      LOCAL Offset, Sum1, Sum2, Z, Xi, Xn, XiPlus1, X1 AS EXT

      LOCAL i%, m%, A$

      X1 = R(p%,1,j&) : Xn = R(p%,Nd%,j&)

      Sum1 = 0##

      FOR i% = 1 TO Nd%-1

            Xi      = R(p%,i%,j&) : XiPlus1 = R(p%,i%+1,j&)

            Sum1 = Sum1 + (Xi-1##)^2*(1##+(SIN(3##*Pi*XiPlus1))^2)

      NEXT i%

      Sum1 = Sum1 + (SIN(Pi*3##*X1))^2 +(Xn-1##)^2*(1##+(SIN(TwoPi*Xn))^2)

      Sum2 = 0##

      FOR i% = 1 TO Nd%
```

```
        Xi = R(p%,i%,j&)

        Sum2 = Sum2 + u(Xi,5##,100##,4##)

    NEXT i%

    Z = Sum1/10## + Sum2

    F13 = -Z

END FUNCTION 'F13()
'------------------
SUB FillArrayAij  'needed for function F14, Shekel's Foxholes

    Aij(1,1)=-32## : Aij(1,2)=-16## : Aij(1,3)=0## : Aij(1,4)=16## : Aij(1,5)=32##
    Aij(1,6)=-32## : Aij(1,7)=-16## : Aij(1,8)=0## : Aij(1,9)=16## : Aij(1,10)=32##
    Aij(1,11)=-32## : Aij(1,12)=-16## : Aij(1,13)=0## : Aij(1,14)=16## : Aij(1,15)=32##
    Aij(1,16)=-32## : Aij(1,17)=-16## : Aij(1,18)=0## : Aij(1,19)=16## : Aij(1,20)=32##
    Aij(1,21)=-32## : Aij(1,22)=-16## : Aij(1,23)=0## : Aij(1,24)=16## : Aij(1,25)=32##

    Aij(2,1)=-32## : Aij(2,2)=-32## : Aij(2,3)=-32## : Aij(2,4)=-32## : Aij(2,5)=-32##
    Aij(2,6)=-16## : Aij(2,7)=-16## : Aij(2,8)=-16## : Aij(2,9)=-16## : Aij(2,10)=-16##
    Aij(2,11)=0## : Aij(2,12)=0## : Aij(2,13)=0## : Aij(2,14)=0## : Aij(2,15)=0##
    Aij(2,16)=16## : Aij(2,17)=16## : Aij(2,18)=16## : Aij(2,19)=16## : Aij(2,20)=16##
    Aij(2,21)=32## : Aij(2,22)=32## : Aij(2,23)=32## : Aij(2,24)=32## : Aij(2,25)=32##

END SUB
'-----
FUNCTION F14(R(),Nd%,p%,j&) 'F14 (2-D) Shekel's Foxholes (INVERTED...)

    LOCAL Sum1, Sum2, Z, Xi AS EXT

    LOCAL i%, jj%

    Sum1 = 0##

    FOR jj% = 1 TO 25
        Sum2 = 0##
        FOR i% = 1 TO 2
            Xi = R(p%,i%,j&)
            Sum2 = Sum2 + (Xi-Aij(i%,jj%))^6
        NEXT i%
        Sum1 = Sum1 + 1##/(jj%+Sum2)
    NEXT j%
    Z = 1##/(0.002##+Sum1)

    F14 = -Z

END FUNCTION 'F14
'-----------
FUNCTION F16(R(),Nd%,p%,j&) 'F16 (2-D) 6-Hump Camel-Back

    LOCAL x1, x2, Z AS EXT

    x1 = R(p%,1,j&) : x2 = R(p%,2,j&)

    Z = 4##*x1^2 - 2.1##*x1^4 + x1^6/3## + x1*x2 - 4*x2^2 + 4*x2^4

    F16 = -Z

END FUNCTION 'F16
'-----------
FUNCTION F15(R(),Nd%,p%,j&) 'F15 (4-D) Kowalik's Function
'Global maximum = -0.0003075 @ (0.1928,0.1908,0.1231,0.1358)

    LOCAL x1, x2, x3, x4, Num, Denom, Z, Aj(), Bj() AS EXT

    LOCAL jj%

    REDIM Aj(1 TO 11), Bj(1 TO 11)

    Aj(1)  = 0.1957## : Bj(1)  = 1##/0.25##
    Aj(2)  = 0.1947## : Bj(2)  = 1##/0.50##
    Aj(3)  = 0.1735## : Bj(3)  = 1##/1.00##
    Aj(4)  = 0.1600## : Bj(4)  = 1##/2.00##
    Aj(5)  = 0.0844## : Bj(5)  = 1##/4.00##
    Aj(6)  = 0.0627## : Bj(6)  = 1##/6.00##
    Aj(7)  = 0.0456## : Bj(7)  = 1##/8.00##
    Aj(8)  = 0.0342## : Bj(8)  = 1##/10.0##
    Aj(9)  = 0.0323## : Bj(9)  = 1##/12.0##
    Aj(10) = 0.0235## : Bj(10) = 1##/14.0##
    Aj(11) = 0.0246## : Bj(11) = 1##/16.0##

    Z = 0##

    x1 = R(p%,1,j&) : x2 = R(p%,2,j&) : x3 = R(p%,3,j&) : x4 = R(p%,4,j&)

    FOR jj% = 1 TO 11
        Num   = x1*(Bj(jj%)^2+Bj(jj%)*x2)
        Denom = Bj(jj%)^2+Bj(jj%)*x3+x4
        Z = Z + (Aj(jj%)-Num/Denom)^2
    NEXT jj%
    F15 = -Z

END FUNCTION 'F15
```

```
'----------

FUNCTION F17(R(),Nd%,p%,j&) 'F17, (2-D) Branin

'Global maximum = -0.398 @ (-3.142.12.275), (3.142,2.275), (9.425,2.425)

    LOCAL x1, x2, Z AS EXT

    x1 = R(p%,1,j&) : x2 = R(p%,2,j&)

    Z = (x2-5.1##*x1^2/(4##*Pi^2)+5##*x1/Pi-6##)^2 + 10##*(1##-1##/(8##*Pi))*COS(x1) + 10##

    F17 = -Z

END FUNCTION 'F17
'----------

FUNCTION F18(R(),Nd%,p%,j&) 'Goldstein-Price 2-D Test Function

'Global maximum = -3 @ (0,-1)

    LOCAL Z, x1, x2, t1, t2 AS EXT

    x1 = R(p%,1,j&) : x2 = R(p%,2,j&)

    t1 = 1##+(x1+x2+1##)^2*(19##-14##*x1+3##*x1^2-14##*x2+6##*x1*x2+3##*x2^2)

    t2 = 30##+(2##*x1-3##*x2)^2*(18##-32##*x1+12##*x1^2+48##*x2-36##*x1*x2+27##*x2^2)

    Z = t1*t2

    F18 = -Z

END FUNCTION 'F18()
'----------

FUNCTION F19(R(),Nd%,p%,j&) 'F19 (3-D) Hartman's Family #1

'Global maximum = 3.86 @ (0.114,0.556,0.852)

    LOCAL Xi, Z, Sum, Aji(), Cj(), Pji() AS EXT

    LOCAL i%, jj%, m%

    REDIM Aji(1 TO 4, 1 TO 3), Cj(1 TO 4), Pji(1 TO 4, 1 TO 3)
    Aji(1,1) = 3.0## : Aji(1,2) = 10## : Aji(1,3) = 30## : Cj(1) = 1.0##
    Aji(2,1) = 0.1## : Aji(2,2) = 10## : Aji(2,3) = 35## : Cj(2) = 1.2##
    Aji(3,1) = 3.0## : Aji(3,2) = 10## : Aji(3,3) = 30## : Cj(3) = 3.0##
    Aji(4,1) = 0.1## : Aji(4,2) = 10## : Aji(4,3) = 35## : Cj(4) = 3.2##

    Pji(1,1) = 0.36890## : Pji(1,2) = 0.1170## : Pji(1,3) = 0.2673##
    Pji(2,1) = 0.46990## : Pji(2,2) = 0.4387## : Pji(2,3) = 0.7470##
    Pji(3,1) = 0.10910## : Pji(3,2) = 0.8732## : Pji(3,3) = 0.5547##
    Pji(4,1) = 0.03815## : Pji(4,2) = 0.5743## : Pji(4,3) = 0.8828##

    Z = 0##

    FOR jj% = 1 TO 4

        Sum = 0##

        FOR i% = 1 TO 3

            Xi = R(p%,i%,j&)

            Sum = Sum + Aji(jj%,i%)*(Xi-Pji(jj%,i%))^2

        NEXT i%

        Z = Z + Cj(jj%)*EXP(-Sum)

    NEXT jj%

    F19 = Z

END FUNCTION 'F19
'----------

FUNCTION F20(R(),Nd%,p%,j&) 'F20 (6-D) Hartman's Family #2

'Global maximum = 3.32 @ (0.201,0.150,0.477,0.275,0.311,0.657)

    LOCAL Xi, Z, Sum, Aji(), Cj(), Pji() AS EXT

    LOCAL i%, jj%, m%

    REDIM Aji(1 TO 4, 1 TO 6), Cj(1 TO 4), Pji(1 TO 4, 1 TO 6)

    Aji(1,1) = 10.0## : Aji(1,2) = 3.00## : Aji(1,3) = 17.0## : Cj(1) = 1.0##
    Aji(2,1) = 0.05## : Aji(2,2) = 10.0## : Aji(2,3) = 17.0## : Cj(2) = 1.2##
    Aji(3,1) = 3.00## : Aji(3,2) = 3.50## : Aji(3,3) = 1.70## : Cj(3) = 3.0##
    Aji(4,1) = 17.0## : Aji(4,2) = 8.00## : Aji(4,3) = 0.05## : Cj(4) = 3.2##

    Aji(1,4) = 3.5## : Aji(1,5) = 1.7## : Aji(1,6) =  8##
    Aji(2,4) = 0.1## : Aji(2,5) =    8## : Aji(2,6) = 14##
    Aji(3,4) =  10## : Aji(3,5) =  17## : Aji(3,6) =  8##
    Aji(4,4) =  10## : Aji(4,5) = 0.1## : Aji(4,6) = 14##

    Pji(1,1) = 0.13120## : Pji(1,2) = 0.1696## : Pji(1,3) = 0.5569##
    Pji(2,1) = 0.23290## : Pji(2,2) = 0.4135## : Pji(2,3) = 0.8307##
    Pji(3,1) = 0.23480## : Pji(3,2) = 0.1415## : Pji(3,3) = 0.3522##
    Pji(4,1) = 0.40470## : Pji(4,2) = 0.8828## : Pji(4,3) = 0.8732##

    Pji(1,4) = 0.01240## : Pji(1,5) = 0.8283## : Pji(1,6) = 0.5886##
    Pji(2,4) = 0.37360## : Pji(2,5) = 0.1004## : Pji(2,6) = 0.9991##
    Pji(3,4) = 0.28830## : Pji(3,5) = 0.3047## : Pji(3,6) = 0.6650##
    Pji(4,4) = 0.57430## : Pji(4,5) = 0.1091## : Pji(4,6) = 0.0381##

    Z = 0##

    FOR jj% = 1 TO 4
```

```
            Sum = 0##

            FOR i% = 1 TO 6

                Xi = R(p%,i%,j&)

                Sum = Sum + Aji(jj%,i%)*(Xi-Pji(jj%,i%))^2

            NEXT i%

            Z = Z + Cj(jj%)*EXP(-Sum)

        NEXT jj%

        F20 = Z

END FUNCTION 'F20

'-----------

FUNCTION F21(R(),Nd%,p%,j&) 'F21 (4-D) Shekel's Family m=5

'Global maximum = 10

        LOCAL Xi, Z, Sum, Aji(), Cj() AS EXT

        LOCAL i%, jj%, m%

        m% = 5 : REDIM Aji(1 TO m%, 1 TO 4), Cj(1 TO m%)

        Aji(1,1) =  4## : Aji(1,2) =   4## : Aji(1,3) =  4## : Aji(1,4) =   4## : Cj(1) = 0.1##
        Aji(2,1) =  1## : Aji(2,2) =   1## : Aji(2,3) =  1## : Aji(2,4) =   1## : Cj(2) = 0.2##
        Aji(3,1) =  8## : Aji(3,2) =   8## : Aji(3,3) =  8## : Aji(3,4) =   8## : Cj(3) = 0.2##
        Aji(4,1) =  6## : Aji(4,2) =   6## : Aji(4,3) =  6## : Aji(4,4) =   6## : Cj(4) = 0.4##
        Aji(5,1) =  3## : Aji(5,2) =   7## : Aji(5,3) =  3## : Aji(5,4) =   7## : Cj(5) = 0.4##

        Z = 0##

        FOR jj% = 1 TO m%  'NOTE:  Index jj% is used to avoid same variable name as j&

            Sum = 0##

            FOR i% = 1 TO 4 'Shekel's family is 4-D only

                Xi = R(p%,i%,j&)

                Sum = Sum + (Xi-Aji(jj%,i%))^2

            NEXT i%

            Z = Z + 1##/(Sum + Cj(jj%))

        NEXT jj%

        F21 = Z

END FUNCTION 'F21

'-----------

FUNCTION F22(R(),Nd%,p%,j&) 'F22 (4-D) Shekel's Family m=7

'Global maximum = 10

        LOCAL Xi, Z, Sum, Aji(), Cj() AS EXT

        LOCAL i%, jj%, m%

        m% = 7 : REDIM Aji(1 TO m%, 1 TO 4), Cj(1 TO m%)

        Aji(1,1) =  4## : Aji(1,2) =   4## : Aji(1,3) =  4## : Aji(1,4) =   4## : Cj(1) = 0.1##
        Aji(2,1) =  1## : Aji(2,2) =   1## : Aji(2,3) =  1## : Aji(2,4) =   1## : Cj(2) = 0.2##
        Aji(3,1) =  8## : Aji(3,2) =   8## : Aji(3,3) =  8## : Aji(3,4) =   8## : Cj(3) = 0.2##
        Aji(4,1) =  6## : Aji(4,2) =   6## : Aji(4,3) =  6## : Aji(4,4) =   6## : Cj(4) = 0.4##
        Aji(5,1) =  3## : Aji(5,2) =   7## : Aji(5,3) =  3## : Aji(5,4) =   7## : Cj(5) = 0.4##
        Aji(6,1) =  2## : Aji(6,2) =   9## : Aji(6,3) =  2## : Aji(6,4) =   9## : Cj(6) = 0.6##
        Aji(7,1) =  5## : Aji(7,2) =   5## : Aji(7,3) =  3## : Aji(7,4) =   3## : Cj(7) = 0.3##

        Z = 0##

        FOR jj% = 1 TO m%  'NOTE:  Index jj% is used to avoid same variable name as j&

            Sum = 0##

            FOR i% = 1 TO 4 'Shekel's family is 4-D only

                Xi = R(p%,i%,j&)

                Sum = Sum + (Xi-Aji(jj%,i%))^2

            NEXT i%

            Z = Z + 1##/(Sum + Cj(jj%))

        NEXT jj%

        F22 = Z

END FUNCTION 'F22

'-----------

FUNCTION F23(R(),Nd%,p%,j&) 'F23 (4-D) Shekel's Family m=10

'Global maximum = 10

        LOCAL Xi, Z, Sum, Aji(), Cj() AS EXT

        LOCAL i%, jj%, m%

        m% = 10 : REDIM Aji(1 TO m%, 1 TO 4), Cj(1 TO m%)

        Aji(1,1) =  4## : Aji(1,2) =   4## : Aji(1,3) =  4## : Aji(1,4) =   4## : Cj(1) = 0.1##
        Aji(2,1) =  1## : Aji(2,2) =   1## : Aji(2,3) =  1## : Aji(2,4) =   1## : Cj(2) = 0.2##
        Aji(3,1) =  8## : Aji(3,2) =   8## : Aji(3,3) =  8## : Aji(3,4) =   8## : Cj(3) = 0.2##
        Aji(4,1) =  6## : Aji(4,2) =   6## : Aji(4,3) =  6## : Aji(4,4) =   6## : Cj(4) = 0.4##
        Aji(5,1) =  3## : Aji(5,2) =   7## : Aji(5,3) =  3## : Aji(5,4) =   7## : Cj(5) = 0.4##
```

```
Aji(6,1)  = 2##  : Aji(6,2)  =   9##  : Aji(6,3)  = 2##  : Aji(6,4)  =   9##  : Cj(6)  = 0.6##
Aji(7,1)  = 5##  : Aji(7,2)  =   5##  : Aji(7,3)  = 3##  : Aji(7,4)  =   3##  : Cj(7)  = 0.3##
Aji(8,1)  = 8##  : Aji(8,2)  =   1##  : Aji(8,3)  = 8##  : Aji(8,4)  =   1##  : Cj(8)  = 0.7##
Aji(9,1)  = 6##  : Aji(9,2)  =   2##  : Aji(9,3)  = 6##  : Aji(9,4)  =   2##  : Cj(9)  = 0.5##
Aji(10,1) = 7##  : Aji(10,2) = 3.6## : Aji(10,3) = 7##  : Aji(10,4) = 3.6## : Cj(10) = 0.5##

        Z = 0##

        FOR jj% = 1 TO m%   'NOTE:  Index jj% is used to avoid same variable name as j&

            Sum = 0##

            FOR i% = 1 TO 4 'Shekel's family is 4-D only

                Xi = R(p%,i%,j&)

                Sum = Sum + (Xi-Aji(jj%,i%))^2

            NEXT i%

            Z = Z + 1##/(Sum + Cj(jj%))

        NEXT jj%

        F23 = Z

END FUNCTION 'F23

'=================================== END FUNCTION DEFINITIONS ===================================
SUB Plot2bbestProbeTrajectories(NumTrajectories%,M(),R(),Np%,Nd%,LastStep&,FunctionName$)

LOCAL TrajectoryNumber%, ProbeNumber%, StepNumber&, N%, M%, ProcID???

LOCAL MaximumFitness, MinimumFitness AS EXT

LOCAL BestProbeThisStep%()

LOCAL BestFitnessThisStep(), TempFitness() AS EXT

LOCAL Annotation$, xCoord$, yCoord$, GnuPlotEXE$, PlotwithLines$

        Annotation$   = ""

        PlotwithLines$ = "YES"  '"NO"

        NumTrajectories% = MIN(Np%,NumTrajectories%)

        GnuPlotEXE$ = "wgnuplot.exe"
'   --------------- Get Min/Max Fitnesses ----------------
        MaximumFitness = M(1,0) : MinimumFitness = M(1,0)  'Note:  M(p%,j&)

        FOR StepNumber& = 0 TO LastStep&

            FOR ProbeNumber% = 1 TO Np%

                IF M(ProbeNumber%,StepNumber&) >= MaximumFitness THEN MaximumFitness = M(ProbeNumber%,StepNumber&)

                IF M(ProbeNumber%,StepNumber&) =< MinimumFitness THEN MinimumFitness = M(ProbeNumber%,StepNumber&)

            NEXT ProbeNumber%

        NEXT StepNumber%
'   ------------ Copy Fitness Array M() into TempFitness to Preserve M() ---------------
        REDIM TempFitness(1 TO Np%, 0 TO LastStep&)

        FOR StepNumber& = 0 TO LastStep&

            FOR ProbeNumber% = 1 TO Np%

                TempFitness(ProbeNumber%,StepNumber&) = M(ProbeNumber%,StepNumber&)

            NEXT ProbeNumber%

        NEXT StepNumber%
'   ----------- LOOP ON TRAJECTORIES -----------
        FOR TrajectoryNumber% = 1 TO NumTrajectories%
'       --------------- Get Trajectory Coordinate Data -----------------
            REDIM BestFitnessThisStep(0 TO LastStep&), BestProbeThisStep(0 TO LastStep&)

            FOR StepNumber& = 0 TO LastStep&

                BestFitnessThisStep(StepNumber&) = TempFitness(1,StepNumber&)

                FOR ProbeNumber% = 1 TO Np%

                    IF TempFitness(ProbeNumber%,StepNumber&) >= BestFitnessThisStep(StepNumber&) THEN

                        BestFitnessThisStep(StepNumber&) = TempFitness(ProbeNumber%,StepNumber&)

                        BestProbeThisStep%(StepNumber&)  = ProbeNumber%

                    END IF

                NEXT ProbeNumber%

            NEXT StepNumber&
'   ----- Create Plot Data File -----
        N% = FREEFILE

        SELECT CASE TrajectoryNumber%

            CASE 1 : OPEN "t1" FOR OUTPUT AS #N%
            CASE 2 : OPEN "t2" FOR OUTPUT AS #N%
            CASE 3 : OPEN "t3" FOR OUTPUT AS #N%
            CASE 4 : OPEN "t4" FOR OUTPUT AS #N%
            CASE 5 : OPEN "t5" FOR OUTPUT AS #N%
```

```
        CASE 6   : OPEN "t6"  FOR OUTPUT AS #N%
        CASE 7   : OPEN "t7"  FOR OUTPUT AS #N%
        CASE 8   : OPEN "t8"  FOR OUTPUT AS #N%
        CASE 9   : OPEN "t9"  FOR OUTPUT AS #N%
        CASE 10  : OPEN "t10" FOR OUTPUT AS #N%

    END SELECT

'    ----------- Write Plot File Data ------------

    FOR StepNumber& = 0 TO LastStep&

      PRINT #N%, USING$("######.########
######.########",R(BestProbeThisStep%(StepNumber&),1,StepNumber&),R(BestProbeThisStep%(StepNumber&),2,StepNumber&))

      TempFitness(BestProbeThisStep%(StepNumber&),StepNumber&) = MinimumFitness 'so that same max will not be found for next trajectory

    NEXT StepNumber%

    CLOSE #N%

  NEXT TrajectoryNumber%

'    ----------------------- Plot Trajectories --------------------------

  CALL CreateGNUplotINIfile(0.13##*ScreenWidth&,0.18##*ScreenHeight&,0.7##*ScreenHeight&,0.7##*ScreenHeight&)

  Annotation$ = ""

  N% = FREEFILE

  OPEN "cmd2d.gp" FOR OUTPUT AS #N%

    PRINT #N%, "set xrange ["+REMOVE$(STR$(XiMin(1)),ANY" ")+":"+REMOVE$(STR$(XiMax(1)),ANY" ")+"]"
    PRINT #N%, "set yrange ["+REMOVE$(STR$(XiMin(2)),ANY" ")+":"+REMOVE$(STR$(XiMax(2)),ANY" ")+"]"

    'PRINT #N%, "set label "  + Quote$ + Annotation$ + Quote$ + " at graph " + xCoord$ + "," + yCoord$
    PRINT #N%, "set grid xtics " + "10"
    PRINT #N%, "set grid ytics " + "10"
    PRINT #N%, "set grid mxtics"
    PRINT #N%, "set grid mytics"
    PRINT #N%, "show grid"
    PRINT #N%, "set title "  + Quote$ + "2D "+ FunctionName$+" TRAJECTORIES OF PROBES WITH BEST\nFITNESSES (ORDERED BY FITNESS)" + "\n" + RunID$ + Quote$
    PRINT #N%, "set xlabel "  + Quote$ + "x1\n\n"                     + Quote$
    PRINT #N%, "set ylabel "  + Quote$ + "\nx2"                       + Quote$

    IF PlotWithLines$ = "YES" THEN

      SELECT CASE NumTrajectories%

        CASE 1  : PRINT #N%, "plot "+Quote$+"t1"+Quote$+" w l lw 3"
        CASE 2  : PRINT #N%, "plot "+Quote$+"t1"+Quote$+" w l lw 3,"+Quote$+"t2"+Quote$+" w l"
        CASE 3  : PRINT #N%, "plot "+Quote$+"t1"+Quote$+" w l lw 3,"+Quote$+"t2"+Quote$+" w l,"+Quote$+"t3"+Quote$+" w l"
        CASE 4  : PRINT #N%, "plot "+Quote$+"t1"+Quote$+" w l lw 3,"+Quote$+"t2"+Quote$+" w l,"+Quote$+"t3"+Quote$+" w
l,"+Quote$+"t4"+Quote$+" w l"
        CASE 5  : PRINT #N%, "plot "+Quote$+"t1"+Quote$+" w l lw 3,"+Quote$+"t2"+Quote$+" w l,"+Quote$+"t3"+Quote$+" w
l,"+Quote$+"t4"+Quote$+" w l,"+Quote$+"t5"+Quote$+" w l"
        CASE 6  : PRINT #N%, "plot "+Quote$+"t1"+Quote$+" w l lw 3,"+Quote$+"t2"+Quote$+" w l,"+Quote$+"t3"+Quote$+" w
l,"+Quote$+"t4"+Quote$+" w l,"+Quote$+"t5"+Quote$+" w l,"+Quote$+"t6"+Quote$+" w l"
        CASE 7  : PRINT #N%, "plot "+Quote$+"t1"+Quote$+" w l lw 3,"+Quote$+"t2"+Quote$+" w l,"+Quote$+"t3"+Quote$+" w
l,"+Quote$+"t4"+Quote$+" w l,"+Quote$+"t5"+Quote$+" w l,"+Quote$+"t6"+Quote$+" w l,"+_
                                             Quote$+"t7"+Quote$+" w l"
        CASE 8  : PRINT #N%, "plot "+Quote$+"t1"+Quote$+" w l lw 3,"+Quote$+"t2"+Quote$+" w l,"+Quote$+"t3"+Quote$+" w
l,"+Quote$+"t4"+Quote$+" w l,"+Quote$+"t5"+Quote$+" w l,"+Quote$+"t6"+Quote$+" w l,"+_
                                             Quote$+"t7"+Quote$+" w l,"+      Quote$+"t8"+Quote$+" w l"
        CASE 9  : PRINT #N%, "plot "+Quote$+"t1"+Quote$+" w l lw 3,"+Quote$+"t2"+Quote$+" w l,"+Quote$+"t3"+Quote$+" w
l,"+Quote$+"t4"+Quote$+" w l,"+Quote$+"t5"+Quote$+" w l,"+Quote$+"t6"+Quote$+" w l,"+_
                                             Quote$+"t7"+Quote$+" w l,"+      Quote$+"t8"+Quote$+" w l,"+Quote$+"t9"+Quote$+" w l"
        CASE 10 : PRINT #N%, "plot "+Quote$+"t1"+Quote$+" w l lw 3,"+Quote$+"t2"+Quote$+" w l,"+Quote$+"t3"+Quote$+" w
l,"+Quote$+"t4"+Quote$+" w l,"+Quote$+"t5"+Quote$+" w l,"+Quote$+"t6"+Quote$+" w l,"+_
                                             Quote$+"t7"+Quote$+" w l,"+      Quote$+"t8"+Quote$+" w l,"+Quote$+"t9"+Quote$+" w
l,"+Quote$+"t10"+Quote$+" w l"
      END SELECT

    ELSE

      SELECT CASE NumTrajectories%

        CASE 1  : PRINT #N%, "plot "+Quote$+"t1"+Quote$+" lw 2"
        CASE 2  : PRINT #N%, "plot "+Quote$+"t1"+Quote$+" lw 2,"+Quote$+"t2"+Quote$
        CASE 3  : PRINT #N%, "plot "+Quote$+"t1"+Quote$+" lw 2,"+Quote$+"t2"+Quote$+" ,"+Quote$+"t3"+Quote$
        CASE 4  : PRINT #N%, "plot "+Quote$+"t1"+Quote$+" lw 2,"+Quote$+"t2"+Quote$+" ,"+Quote$+"t3"+Quote$+" ,"+Quote$+"t4"+Quote$
        CASE 5  : PRINT #N%, "plot "+Quote$+"t1"+Quote$+" lw 2,"+Quote$+"t2"+Quote$+" ,"+Quote$+"t3"+Quote$+" ,"+Quote$+"t4"+Quote$+"
,"+Quote$+"t5"+Quote$
        CASE 6  : PRINT #N%, "plot "+Quote$+"t1"+Quote$+" lw 2,"+Quote$+"t2"+Quote$+" ,"+Quote$+"t3"+Quote$+" ,"+Quote$+"t4"+Quote$+"
,"+Quote$+"t5"+Quote$+" ,"+Quote$+"t6"+Quote$
        CASE 7  : PRINT #N%, "plot "+Quote$+"t1"+Quote$+" lw 2,"+Quote$+"t2"+Quote$+" ,"+Quote$+"t3"+Quote$+" ,"+Quote$+"t4"+Quote$+"
,"+Quote$+"t5"+Quote$+" ,"+Quote$+"t6"+Quote$+" ,"+_
                                             Quote$+"t7"+Quote$
        CASE 8  : PRINT #N%, "plot "+Quote$+"t1"+Quote$+" lw 2,"+Quote$+"t2"+Quote$+" ,"+Quote$+"t3"+Quote$+" ,"+Quote$+"t4"+Quote$+"
,"+Quote$+"t5"+Quote$+" ,"+Quote$+"t6"+Quote$+" ,"+_
                                             Quote$+"t7"+Quote$+" ,"+      Quote$+"t8"+Quote$
        CASE 9  : PRINT #N%, "plot "+Quote$+"t1"+Quote$+" lw 2,"+Quote$+"t2"+Quote$+" ,"+Quote$+"t3"+Quote$+" ,"+Quote$+"t4"+Quote$+"
,"+Quote$+"t5"+Quote$+" ,"+Quote$+"t6"+Quote$+" ,"+_
                                             Quote$+"t7"+Quote$+" ,"+      Quote$+"t8"+Quote$+" ,"+Quote$+"t9"+Quote$
        CASE 10 : PRINT #N%, "plot "+Quote$+"t1"+Quote$+" lw 2,"+Quote$+"t2"+Quote$+" ,"+Quote$+"t3"+Quote$+" ,"+Quote$+"t4"+Quote$+"
,"+Quote$+"t5"+Quote$+" ,"+Quote$+"t6"+Quote$+" ,"+_
                                             Quote$+"t7"+Quote$+" ,"+      Quote$+"t8"+Quote$+" ,"+Quote$+"t9"+Quote$+" ,"+Quote$+"t10"+Quote$
      END SELECT

    END IF

    CLOSE #N%

    ProcID??? = SHELL(GnuPlotEXE$+" cmd2d.gp -") : CALL Delay(1##)

END SUB 'Plot2DbestProbeTrajectories()

'----

SUB Plot2DindividualProbeTrajectories(NumTrajectories%,M(),R(),Np%,Nd%,LastStep&,FunctionName$)

LOCAL ProbeNumber%, StepNumber&, N%, ProcID???

LOCAL Annotation$, xCoord$, yCoord$, GnuPlotEXE$, PlotwithLines$
```

```basic
    NumTrajectories% = MIN(Np%,NumTrajectories%)

    Annotation$    = ""

    PlotwithLines$ = "YES" '"NO"

    GnuPlotEXE$ = "wgnuplot.exe"
'  ------------- LOOP ON PROBES ---------------

    FOR ProbeNumber% = 1 TO MIN(NumTrajectories%,Np%)
'  ----- Create Plot Data File -----

    N% = FREEFILE

    SELECT CASE ProbeNumber%

        CASE 1  : OPEN "p1"  FOR OUTPUT AS #N%
        CASE 2  : OPEN "p2"  FOR OUTPUT AS #N%
        CASE 3  : OPEN "p3"  FOR OUTPUT AS #N%
        CASE 4  : OPEN "p4"  FOR OUTPUT AS #N%
        CASE 5  : OPEN "p5"  FOR OUTPUT AS #N%
        CASE 6  : OPEN "p6"  FOR OUTPUT AS #N%
        CASE 7  : OPEN "p7"  FOR OUTPUT AS #N%
        CASE 8  : OPEN "p8"  FOR OUTPUT AS #N%
        CASE 9  : OPEN "p9"  FOR OUTPUT AS #N%
        CASE 10 : OPEN "p10" FOR OUTPUT AS #N%
        CASE 11 : OPEN "p11" FOR OUTPUT AS #N%
        CASE 12 : OPEN "p12" FOR OUTPUT AS #N%
        CASE 13 : OPEN "p13" FOR OUTPUT AS #N%
        CASE 14 : OPEN "p14" FOR OUTPUT AS #N%
        CASE 15 : OPEN "p15" FOR OUTPUT AS #N%
        CASE 16 : OPEN "p16" FOR OUTPUT AS #N%

    END SELECT

'  ----------- Write Plot File Data ------------

    FOR StepNumber& = 0 TO LastStep&

        PRINT #N%, USING$("#####.####### #####.#######",R(ProbeNumber%,1,StepNumber&),R(ProbeNumber%,2,StepNumber&))

    NEXT StepNumber%

    CLOSE #N%

    NEXT ProbeNumber%
'  ---------------------------------------------- Plot Trajectories ----------------------------------------------

'usage:  CALL CreateGNUplotINIfile(PlotWindowULC_X%,PlotWindowULC_Y%,PlotWindowWidth%,PlotWindowHeight%)

    CALL CreateGNUplotINIfile(0.17##*ScreenWidth&,0.22##*ScreenHeight&,0.7##*ScreenWidth&,0.7##*ScreenHeight&)

    Annotation$ = ""

    N% = FREEFILE

    OPEN "cmd2d.gp" FOR OUTPUT AS #N%

        PRINT #N%, "set xrange ["+REMOVE$(STR$(XiMin(1)),ANY"" )+":"+REMOVE$(STR$(XiMax(1)),ANY" ")+"]"
        PRINT #N%, "set yrange ["+REMOVE$(STR$(XiMin(2)),ANY" ")+":"+REMOVE$(STR$(XiMax(2)),ANY" ")+"]"

        PRINT #N%, "set grid xtics " + "10"
        PRINT #N%, "set grid ytics " + "10"
        PRINT #N%, "set grid mxtics"
        PRINT #N%, "set grid mytics"
        PRINT #N%, "show grid"
        PRINT #N%, "set title " + Quote$ + "2D "+ FunctionName$+" INDIVIDUAL PROBE TRAJECTORIES\n(ORDERED BY PROBE #)" + "\n" + RunID$ + Quote$
        PRINT #N%, "set xlabel " + Quote$ + "x1\n\n"                       + Quote$
        PRINT #N%, "set ylabel " + Quote$ + "\nx2"                         + Quote$

        IF PlotwithLines$ = "YES" THEN

            SELECT CASE NumTrajectories%

                CASE 1  : PRINT #N%, "plot "+Quote$+"p1"  +Quote$+" w l lw l"
                CASE 2  : PRINT #N%, "plot "+Quote$+"p1"  +Quote$+" w l lw l,"+Quote$+"p2"+Quote$+" w l"
                CASE 3  : PRINT #N%, "plot "+Quote$+"p1"  +Quote$+" w l lw l,"+Quote$+"p2"+Quote$+" w l,"+Quote$+"p3"+Quote$+" w l"
                CASE 4  : PRINT #N%, "plot "+Quote$+"p1"  +Quote$+" w l lw l,"+Quote$+"p2"+Quote$+" w l,"+Quote$+"p3"+Quote$+" w
l,"+Quote$+"p4"+Quote$+" w l"
                CASE 5  : PRINT #N%, "plot "+Quote$+"p1"  +Quote$+" w l lw l,"+Quote$+"p2"+Quote$+" w l,"+Quote$+"p3"+Quote$+" w
l,"+Quote$+"p4"+Quote$+" w l,"+Quote$+"p5"+Quote$+" w l"
                CASE 6  : PRINT #N%, "plot "+Quote$+"p1"  +Quote$+" w l lw l,"+Quote$+"p2"+Quote$+" w l,"+Quote$+"p3"+Quote$+" w
l,"+Quote$+"p4"+Quote$+" w l,"+Quote$+"p5"+Quote$+" w l,"+Quote$+"p6"+Quote$+" w l"
                CASE 7  : PRINT #N%, "plot "+Quote$+"p1"  +Quote$+" w l lw l,"+Quote$+"p2"+Quote$+" w l,"+Quote$+"p3"+Quote$+" w
l,"+Quote$+"p4"+Quote$+" w l,"+Quote$+"p5"+Quote$+" w l,"+Quote$+"p6"+Quote$+" w l,+_
                                      Quote$+"p7"     +Quote$+" w l"
                CASE 8  : PRINT #N%, "plot "+Quote$+"p1"  +Quote$+" w l lw l,"+Quote$+"p2"+Quote$+" w l,"+Quote$+"p3"+Quote$+" w
l,"+Quote$+"p4"+Quote$+" w l,"+Quote$+"p5"+Quote$+" w l,"+Quote$+"p6"+Quote$+" w l,+_
                                      Quote$+"p7"     +Quote$+" w l,"     +Quote$+"p8"+Quote$+" w l"
                CASE 9  : PRINT #N%, "plot "+Quote$+"p1"  +Quote$+" w l lw l,"+Quote$+"p2"+Quote$+" w l,"+Quote$+"p3"+Quote$+" w
l,"+Quote$+"p4"+Quote$+" w l,"+Quote$+"p5"+Quote$+" w l,"+Quote$+"p6"+Quote$+" w l,+_
                                      Quote$+"p7"     +Quote$+" w l,"     +Quote$+"p8"+Quote$+" w l,"+Quote$+"p9"+Quote$+" w l"

                CASE 10 : PRINT #N%, "plot "+Quote$+"p1"  +Quote$+" w l lw l,"+Quote$+"p2"+Quote$+" w l,"+Quote$+"p3"+Quote$+" w l,"+Quote$+"p4"
+Quote$+" w l,"+Quote$+"p5"+Quote$+" w l,"+Quote$+"p6"+Quote$+" w l,+_
                                      Quote$+"p7"     +Quote$+" w l,"     +Quote$+"p8"+Quote$+" w l,"+Quote$+"p9"+Quote$+" w l,"+Quote$+"p10"+Quote$+" w l"

                CASE 11 : PRINT #N%, "plot "+Quote$+"p1"  +Quote$+" w l lw l,"+Quote$+"p2"+Quote$+" w l,"+Quote$+"p3"+Quote$+" w l,"+Quote$+"p4"
+Quote$+" w l,"+Quote$+"p5"+Quote$+" w l,"+Quote$+"p6"+Quote$+" w l,+_
                                      Quote$+"p7"     +Quote$+" w l,+_     Quote$+"p8"+Quote$+" w l,"+Quote$+"p9"+Quote$+" w
l,"+Quote$+"p10"+Quote$+" w l,"+Quote$+"p11"+Quote$+" w l"

                CASE 12 : PRINT #N%, "plot "+Quote$+"p1"  +Quote$+" w l lw l,"+Quote$+"p2"+Quote$+" w l,"+Quote$+"p3"+Quote$+" w l,"+Quote$+"p4"
+Quote$+" w l,"+Quote$+"p5"+Quote$+" w l,"+Quote$+"p6"+Quote$+" w l,+_
                                      Quote$+"p7"     +Quote$+" w l,"     +Quote$+"p8"+Quote$+" w l,"+Quote$+"p9"+Quote$+" w
l,"+Quote$+"p10"+Quote$+" w l,"+Quote$+"p11"+Quote$+" w l,"+Quote$+"p12"+Quote$+" w l"

                CASE 13 : PRINT #N%, "plot "+Quote$+"p1"  +Quote$+" w l lw l,"+Quote$+"p2"+Quote$+" w l,"+Quote$+"p3"+Quote$+" w l,"+Quote$+"p4"
+Quote$+" w l,"+Quote$+"p5"+Quote$+" w l,"+Quote$+"p6"+Quote$+" w l,+_
                                      Quote$+"p7"     +Quote$+" w l,"     +Quote$+"p8"+Quote$+" w l,"+Quote$+"p9"+Quote$+" w
l,"+Quote$+"p10"+Quote$+" w l,"+Quote$+"p11"+Quote$+" w l,"+Quote$+"p12"+Quote$+" w l,+_
                                      Quote$+"p13"    +Quote$+" w l"
```

```
                CASE 14 : PRINT #N%, "plot "+Quote$+"p1"   +Quote$+" w l lw 1,"+Quote$+"p2"   +Quote$+" w l,"+Quote$+"p3" +Quote$+" w l,"+Quote$+"p4"
            +Quote$+" w l,"+Quote$+"p5"   +Quote$+"   w l,"+Quote$+"p6"   +Quote$+" w l,"+_
                                          +Quote$+"p7"   +Quote$+" w l,"+_
            l,"+Quote$+"p10"+Quote$+" w l,"+Quote$+"p11"+Quote$+" w l,"+Quote$+"p12"+Quote$+" w l,"+_
                                          +Quote$+"p13"   +Quote$+" w l,"      +Quote$+"p14"+Quote$+" w l"

                CASE 15 : PRINT #N%, "plot "+Quote$+"p1"   +Quote$+" w l lw 1,"+Quote$+"p2"   +Quote$+" w l,"+Quote$+"p3" +Quote$+" w l,"+Quote$+"p4"
            +Quote$+" w l,"+Quote$+"p5"   +Quote$+"   w l,"+Quote$+"p6"   +Quote$+" w l,"+_
                                          +Quote$+"p7"   +Quote$+" w l,"+_
            l,"+Quote$+"p10"+Quote$+" w l,"+Quote$+"p11"+Quote$+" w l,"+Quote$+"p12"+Quote$+" w l,"+_
                                          +Quote$+"p13"   +Quote$+" w l,"      +Quote$+"p14"+Quote$+" w l,"+Quote$+"p15"+Quote$+" w l"

                CASE 16 : PRINT #N%, "plot "+Quote$+"p1"   +Quote$+" w l lw 1,"+Quote$+"p2"   +Quote$+" w l,"+Quote$+"p3" +Quote$+" w l,"+Quote$+"p4"
            +Quote$+" w l,"+Quote$+"p5"   +Quote$+"   w l,"+Quote$+"p6"   +Quote$+" w l,"+_
                                          +Quote$+"p7"   +Quote$+" w l,"+_
            l,"+Quote$+"p10"+Quote$+" w l,"+Quote$+"p11"+Quote$+" w l,"+Quote$+"p12"+Quote$+" w l,"+_
                                          +Quote$+"p13"   +Quote$+" w l,"      +Quote$+"p14"+Quote$+" w l,"+Quote$+"p15"+Quote$+" w
            l,"+Quote$+"p16"+Quote$+" w l"
                END SELECT

            ELSE

                SELECT CASE NumTrajectories%
                    CASE 1  : PRINT #N%, "plot "+Quote$+"p1"+Quote$+" lw 1"
                    CASE 2  : PRINT #N%, "plot "+Quote$+"p1"+Quote$+" lw 1,"+Quote$+"p2"+Quote$
                    CASE 3  : PRINT #N%, "plot "+Quote$+"p1"+Quote$+" lw 1,"+Quote$+"p2"+Quote$+"  ,"+Quote$+"p3"+Quote$
                    CASE 4  : PRINT #N%, "plot "+Quote$+"p1"+Quote$+" lw 1,"+Quote$+"p2"+Quote$+"  ,"+Quote$+"p3"+Quote$+"  ,"+Quote$+"p4"+Quote$
                    CASE 5  : PRINT #N%, "plot "+Quote$+"p1"+Quote$+" lw 1,"+Quote$+"p2"+Quote$+"  ,"+Quote$+"p3"+Quote$+"  ,"+Quote$+"p4"+Quote$+"
            ,"+Quote$+"p5"+Quote$
                    CASE 6  : PRINT #N%, "plot "+Quote$+"p1"+Quote$+" lw 1,"+Quote$+"p2"+Quote$+"  ,"+Quote$+"p3"+Quote$+"  ,"+Quote$+"p4"+Quote$+"
            ,"+Quote$+"p5"+Quote$+"  ,"+Quote$+"p6"+Quote$
                    CASE 7  : PRINT #N%, "plot "+Quote$+"p1"+Quote$+" lw 1,"+Quote$+"p2"+Quote$+"  ,"+Quote$+"p3"+Quote$+"  ,"+Quote$+"p4"+Quote$+"
            ,"+Quote$+"p5"+Quote$+"  ,"+Quote$+"p6"+Quote$+"  ,"+_
                                          +Quote$+"p7"+Quote$
                    CASE 8  : PRINT #N%, "plot "+Quote$+"p1"+Quote$+" lw 1,"+Quote$+"p2"+Quote$+"  ,"+Quote$+"p3"+Quote$+"  ,"+Quote$+"p4"+Quote$+"
            ,"+Quote$+"p5"+Quote$+"  ,"+Quote$+"p6"+Quote$+"  ,"+_
                                          +Quote$+"p7"+Quote$+"  ,"+Quote$+"p8"+Quote$
                    CASE 9  : PRINT #N%, "plot "+Quote$+"p1"+Quote$+" lw 1,"+Quote$+"p2"+Quote$+"  ,"+Quote$+"p3"+Quote$+"  ,"+Quote$+"p4"+Quote$+"
            ,"+Quote$+"p5"+Quote$+"  ,"+Quote$+"p6"+Quote$+"  ,"+_
                                          +Quote$+"p7"+Quote$+"  ,"      +Quote$+"p8"+Quote$+"  ,"+Quote$+"p9"+Quote$
                    CASE 10 : PRINT #N%, "plot "+Quote$+"p1"+Quote$+" lw 1,"+Quote$+"p2"+Quote$+"  ,"+Quote$+"p3"+Quote$+"  ,"+Quote$+"p4"  +Quote$+"
            ,"+Quote$+"p5"+Quote$+"  ,"+Quote$+"p6"+Quote$+"  ,"+_
                                          +Quote$+"p7"+Quote$+"  ,"      +Quote$+"p8"   +Quote$+"  ,"+Quote$+"p9"+Quote$+"  ,"+Quote$+"p10"+Quote$
                    CASE 11 : PRINT #N%, "plot "+Quote$+"p1"+Quote$+" lw 1,"+Quote$+"p2"+Quote$+"  ,"+Quote$+"p3"+Quote$+"  ,"+Quote$+"p4"  +Quote$+"
            ,"+Quote$+"p5"  +Quote$+"  ,"+Quote$+"p6"+Quote$+"  ,"+_
                                          +Quote$+"p7"+Quote$+"  ,"      +Quote$+"p8"   +Quote$+"  ,"+Quote$+"p9"+Quote$+"  ,"+Quote$+"p10"+Quote$+"
            ,"+Quote$+"p11"+Quote$
                    CASE 12 : PRINT #N%, "plot "+Quote$+"p1"+Quote$+" lw 1,"+Quote$+"p2"+Quote$+"  ,"+Quote$+"p3"+Quote$+"  ,"+Quote$+"p4"  +Quote$+"
            ,"+Quote$+"p5"  +Quote$+"  ,"+Quote$+"p6"  +Quote$+"  ,"+_
                                          +Quote$+"p7"+Quote$+"  ,"      +Quote$+"p8"   +Quote$+"  ,"+Quote$+"p9"+Quote$+"  ,"+Quote$+"p10"+Quote$+"
            ,"+Quote$+"p11"+Quote$+"  ,"+Quote$+"p12"+Quote$
                    CASE 13 : PRINT #N%, "plot "+Quote$+"p1"   +Quote$+" lw 1,"+Quote$+"p2"+Quote$+"  ,"+Quote$+"p3"+Quote$+"  ,"+Quote$+"p4"  +Quote$+"
            ,"+Quote$+"p5"  +Quote$+"  ,"      +Quote$+"p6"   +Quote$+"  ,"+_
                                          +Quote$+"p7"   +Quote$+"  ,"      +Quote$+"p8"   +Quote$+"  ,"+Quote$+"p9"+Quote$+"  ,"+Quote$+"p10"+Quote$+"
            ,"+Quote$+"p11"+Quote$+"  ,"  +Quote$+"p12"+Quote$+"  ,"+_
                                          +Quote$+"p13"+Quote$
                    CASE 14 : PRINT #N%, "plot "+Quote$+"p1"   +Quote$+" lw 1,"+Quote$+"p2"+Quote$+"  ,"+Quote$+"p3"+Quote$+"  ,"+Quote$+"p4"  +Quote$+"
            ,"+Quote$+"p5"  +Quote$+"  ,"      +Quote$+"p6"   +Quote$+"  ,"+_
                                          +Quote$+"p7"   +Quote$+"  ,"      +Quote$+"p8"   +Quote$+"  ,"+Quote$+"p9"+Quote$+"  ,"+Quote$+"p10"+Quote$+"
            ,"+Quote$+"p11"+Quote$+"  ,"  +Quote$+"p12"+Quote$+"  ,"+_
                                          +Quote$+"p13"+Quote$+"  ,"      +Quote$+"p14"+Quote$
                    CASE 15 : PRINT #N%, "plot "+Quote$+"p1"   +Quote$+" lw 1,"+Quote$+"p2"+Quote$+"  ,"+Quote$+"p3" +Quote$+"  ,"+Quote$+"p4"  +Quote$+"
            ,"+Quote$+"p5"  +Quote$+"  ,"      +Quote$+"p6"   +Quote$+"  ,"+_
                                          +Quote$+"p7"   +Quote$+"  ,"      +Quote$+"p8"   +Quote$+"  ,"+Quote$+"p9" +Quote$+"  ,"+Quote$+"p10"+Quote$+"
            ,"+Quote$+"p11"+Quote$+"  ,"  +Quote$+"p12"+Quote$+"  ,"+_
                                          +Quote$+"p13"+Quote$+"  ,"      +Quote$+"p14"+Quote$+"  ,"+Quote$+"p15"+Quote$
                    CASE 16 : PRINT #N%, "plot "+Quote$+"p1"   +Quote$+" lw 1,"+Quote$+"p2"+Quote$+"  ,"+Quote$+"p3"+Quote$+"  ,"+Quote$+"p4"  +Quote$+"
            ,"+Quote$+"p5"  +Quote$+"  ,"      +Quote$+"p6"   +Quote$+"  ,"+_
                                          +Quote$+"p7"   +Quote$+"  ,"      +Quote$+"p8"   +Quote$+"  ,"+Quote$+"p9" +Quote$+"  ,"+Quote$+"p10"+Quote$+"
            ,"+Quote$+"p11"+Quote$+"  ,"  +Quote$+"p12"+Quote$+"  ,"+_
                                          +Quote$+"p13"+Quote$+"  ,"      +Quote$+"p14"+Quote$+"  ,"+Quote$+"p15"+Quote$+"  ,"+Quote$+"p16"+Quote$
                END SELECT

            END IF

        CLOSE #N%

        ProcID??? = SHELL(GnuPlotEXE$+" cmd2d.gp -") : CALL Delay(1##)

END SUB 'Plot2DindividualProbeTrajectories()

'----

SUB Plot3DbestProbeTrajectories(NumTrajectories%,M(),R(),Np%,Nd%,LastStep&,FunctionName$) 'XYZZY

LOCAL TrajectoryNumber%, ProbeNumber%, StepNumber&, N%, M%, ProcID???

LOCAL MaximumFitness, MinimumFitness AS EXT

LOCAL BestProbeThisStep%()

LOCAL BestFitnessThisStep(), TempFitness() AS EXT

LOCAL Annotation$, xCoord$, yCoord$, zCoord$, GnuPlotEXE$, PlotwithLines$

    Annotation$    = ""

    PlotwithLines$ = "NO" '"YES" '"NO"

    NumTrajectories% = MIN(Np%,NumTrajectories%)

    GnuPlotEXE$ = "wgnuplot.exe"
'    --------------- Get Min/Max Fitnesses -----------------
    MaximumFitness = M(1,0) : MinimumFitness = M(1,0)  'Note:  M(p%,j&)

    FOR StepNumber& = 0 TO LastStep&
```

```
        FOR ProbeNumber% = 1 TO Np%

            IF M(ProbeNumber%,StepNumber&) >= MaximumFitness THEN MaximumFitness = M(ProbeNumber%,StepNumber&)

            IF M(ProbeNumber%,StepNumber&) =< MinimumFitness THEN MinimumFitness = M(ProbeNumber%,StepNumber&)

        NEXT ProbeNumber%

    NEXT StepNumber%

'    ------------- Copy Fitness Array M() into TempFitness to Preserve M() ----------------

    REDIM TempFitness(1 TO Np%, 0 TO LastStep&)

    FOR StepNumber& = 0 TO LastStep&

        FOR ProbeNumber% = 1 TO Np%

            TempFitness(ProbeNumber%,StepNumber&) = M(ProbeNumber%,StepNumber&)

        NEXT ProbeNumber%

    NEXT StepNumber%

'    ----------- LOOP ON TRAJECTORIES -----------

    FOR TrajectoryNumber% = 1 TO NumTrajectories%

'        -------------- Get Trajectory Coordinate Data -----------------

        REDIM BestFitnessThisStep(0 TO LastStep&), BestProbeThisStep%(0 TO LastStep&)

        FOR StepNumber& = 0 TO LastStep&

            BestFitnessThisStep(StepNumber&) = TempFitness(1,StepNumber&)

            FOR ProbeNumber% = 1 TO Np%

                IF TempFitness(ProbeNumber%,StepNumber&) >= BestFitnessThisStep(StepNumber&) THEN

                    BestFitnessThisStep(StepNumber&) = TempFitness(ProbeNumber%,StepNumber&)

                    BestProbeThisStep%(StepNumber&)  = ProbeNumber%

                END IF

            NEXT ProbeNumber%

        NEXT StepNumber&

'    ----- Create Plot Data File -----

    N% = FREEFILE

    SELECT CASE TrajectoryNumber%

        CASE 1  : OPEN "t1"  FOR OUTPUT AS #N%
        CASE 2  : OPEN "t2"  FOR OUTPUT AS #N%
        CASE 3  : OPEN "t3"  FOR OUTPUT AS #N%
        CASE 4  : OPEN "t4"  FOR OUTPUT AS #N%
        CASE 5  : OPEN "t5"  FOR OUTPUT AS #N%
        CASE 6  : OPEN "t6"  FOR OUTPUT AS #N%
        CASE 7  : OPEN "t7"  FOR OUTPUT AS #N%
        CASE 8  : OPEN "t8"  FOR OUTPUT AS #N%
        CASE 9  : OPEN "t9"  FOR OUTPUT AS #N%
        CASE 10 : OPEN "t10" FOR OUTPUT AS #N%

    END SELECT

'    ----------- Write Plot File Data ------------

    FOR StepNumber& = 0 TO LastStep&

        PRINT #N%, USING$("######.######## ######.########
######.########",R(BestProbeThisStep%(StepNumber&),1,StepNumber&),R(BestProbeThisStep%(StepNumber&),2,StepNumber&),R(BestProbeThisStep%(StepNumber&),
3,StepNumber&))+CHR$(13)

        TempFitness(BestProbeThisStep%(StepNumber&),StepNumber&) = MinimumFitness 'so that same max will not be found for next trajectory

    NEXT StepNumber%

    CLOSE #N%

    NEXT TrajectoryNumber%

'    ------------------------- Plot Trajectories -------------------------

'CALL CreateGNUplotINIfile(0.1##*Screenwidth&,0.25##*ScreenHeight&,0.6##*ScreenHeight&,0.6##*ScreenHeight&)

    Annotation$ = ""

    N% = FREEFILE

    OPEN "cmd3d.gp" FOR OUTPUT AS #N%

    PRINT #N%, "set pm3d"
    PRINT #N%, "show pm3d"
    PRINT #N%, "set hidden3d"
    PRINT #N%, "set view 45, 45, 1, 1"

    PRINT #N%, "unset colorbox"

    PRINT #N%, "set xrange [" + REMOVE$(STR$(XiMin(1)),ANY"" ) + ":" + REMOVE$(STR$(XiMax(1)),ANY"" ) + "]"
    PRINT #N%, "set yrange [" + REMOVE$(STR$(XiMin(2)),ANY"" ) + ":" + REMOVE$(STR$(XiMax(2)),ANY"" ) + "]"
    PRINT #N%, "set zrange [" + REMOVE$(STR$(XiMin(3)),ANY"" ) + ":" + REMOVE$(STR$(XiMax(3)),ANY"" ) + "]"

    PRINT #N%, "set grid xtics ytics ztics"
    PRINT #N%, "show grid"
    PRINT #N%, "set title " + Quote$ + "3D " + FunctionName$ + " PROBE TRAJECTORIES" + "\n" + RunID$ + Quote$
    PRINT #N%, "set xlabel " + Quote$ + "x1"                 + Quote$
    PRINT #N%, "set ylabel " + Quote$ + "x2"                 + Quote$
    PRINT #N%, "set zlabel " + Quote$ + "x3"                 + Quote$

    IF PlotWithLines$ = "YES" THEN

        SELECT CASE NumTrajectories%
```

```
        CASE 1  : PRINT #N%, "splot "+Quote$+"t1"+Quote$+" w l lw 3"
        CASE 2  : PRINT #N%, "splot "+Quote$+"t1"+Quote$+" w l lw 3,"+Quote$+"t2"+Quote$+" w l"
        CASE 3  : PRINT #N%, "splot "+Quote$+"t1"+Quote$+" w l lw 3,"+Quote$+"t2"+Quote$+" w l,"+Quote$+"t3"+Quote$+" w l"
        CASE 4  : PRINT #N%, "splot "+Quote$+"t1"+Quote$+" w l lw 3,"+Quote$+"t2"+Quote$+" w l,"+Quote$+"t3"+Quote$+" w l,"+Quote$+"t4"+Quote$+"
w l"
        CASE 5  : PRINT #N%, "splot "+Quote$+"t1"+Quote$+" w l lw 3,"+Quote$+"t2"+Quote$+" w l,"+Quote$+"t3"+Quote$+" w l,"+Quote$+"t4"+Quote$+"
w l,"+Quote$+"t5"+Quote$+" w l"
        CASE 6  : PRINT #N%, "splot "+Quote$+"t1"+Quote$+" w l lw 3,"+Quote$+"t2"+Quote$+" w l,"+Quote$+"t3"+Quote$+" w l,"+Quote$+"t4"+Quote$+"
w l,"+Quote$+"t5"+Quote$+" w l,"+Quote$+"t6"+Quote$+" w l"
        CASE 7  : PRINT #N%, "splot "+Quote$+"t1"+Quote$+" w l lw 3,"+Quote$+"t2"+Quote$+" w l,"+Quote$+"t3"+Quote$+" w l,"+Quote$+"t4"+Quote$+"
w l,"+Quote$+"t5"+Quote$+" w l,"+Quote$+"t6"+Quote$+" w l,"+
                             Quote$+"t7"+Quote$+" w l"
        CASE 8  : PRINT #N%, "splot "+Quote$+"t1"+Quote$+" w l lw 3,"+Quote$+"t2"+Quote$+" w l,"+Quote$+"t3"+Quote$+" w l,"+Quote$+"t4"+Quote$+"
w l,"+Quote$+"t5"+Quote$+" w l,"+Quote$+"t6"+Quote$+" w l,"+
                             Quote$+"t7"+Quote$+" w l,"        +Quote$+"t8"+Quote$+" w l"
        CASE 9  : PRINT #N%, "splot "+Quote$+"t1"+Quote$+" w l lw 3,"+Quote$+"t2"+Quote$+" w l,"+Quote$+"t3"+Quote$+" w l,"+Quote$+"t4"+Quote$+"
w l,"+Quote$+"t5"+Quote$+" w l,"+Quote$+"t6"+Quote$+" w l,"+
                             Quote$+"t7"+Quote$+" w l,"        +Quote$+"t8"+Quote$+" w l,"+Quote$+"t9"+Quote$+" w l"
        CASE 10 : PRINT #N%, "splot "+Quote$+"t1"+Quote$+" w l lw 3,"+Quote$+"t2"+Quote$+" w l,"+Quote$+"t3"+Quote$+" w l,"+Quote$+"t4"+Quote$+"
w l,"+Quote$+"t5"+Quote$+" w l,"+Quote$+"t6"+Quote$+" w l,"+
                             Quote$+"t7"+Quote$+" w l,"        +Quote$+"t8"+Quote$+" w l,"+Quote$+"t9"+Quote$+" w l,"+Quote$+"t10"+Quote$+"
w l"
    END SELECT

  ELSE

    SELECT CASE NumTrajectories%

        CASE 1  : PRINT #N%, "splot "+Quote$+"t1"+Quote$+" lw 2"
        CASE 2  : PRINT #N%, "splot "+Quote$+"t1"+Quote$+" lw 2,"+Quote$+"t2"+Quote$
        CASE 3  : PRINT #N%, "splot "+Quote$+"t1"+Quote$+" lw 2,"+Quote$+"t2"+Quote$+" ,"+Quote$+"t3"+Quote$
        CASE 4  : PRINT #N%, "splot "+Quote$+"t1"+Quote$+" lw 2,"+Quote$+"t2"+Quote$+" ,"+Quote$+"t3"+Quote$+" ,"+Quote$+"t4"+Quote$
        CASE 5  : PRINT #N%, "splot "+Quote$+"t1"+Quote$+" lw 2,"+Quote$+"t2"+Quote$+" ,"+Quote$+"t3"+Quote$+" ,"+Quote$+"t4"+Quote$+"
,"+Quote$+"t5"+Quote$
        CASE 6  : PRINT #N%, "splot "+Quote$+"t1"+Quote$+" lw 2,"+Quote$+"t2"+Quote$+" ,"+Quote$+"t3"+Quote$+" ,"+Quote$+"t4"+Quote$+"
,"+Quote$+"t5"+Quote$+" ,"+Quote$+"t6"+Quote$
        CASE 7  : PRINT #N%, "splot "+Quote$+"t1"+Quote$+" lw 2,"+Quote$+"t2"+Quote$+" ,"+Quote$+"t3"+Quote$+" ,"+Quote$+"t4"+Quote$+"
,"+Quote$+"t5"+Quote$+" ,"+Quote$+"t6"+Quote$+"
                             Quote$+"t7"+Quote$
        CASE 8  : PRINT #N%, "splot "+Quote$+"t1"+Quote$+" lw 2,"+Quote$+"t2"+Quote$+" ,"+Quote$+"t3"+Quote$+" ,"+Quote$+"t4"+Quote$+"
,"+Quote$+"t5"+Quote$+" ,"+Quote$+"t6"+Quote$+"
                             Quote$+"t7"+Quote$+" ,"     +Quote$+"t8"+Quote$
        CASE 9  : PRINT #N%, "splot "+Quote$+"t1"+Quote$+" lw 2,"+Quote$+"t2"+Quote$+" ,"+Quote$+"t3"+Quote$+" ,"+Quote$+"t4"+Quote$+"
,"+Quote$+"t5"+Quote$+" ,"+Quote$+"t6"+Quote$+"
                             Quote$+"t7"+Quote$+" ,"     +Quote$+"t8"+Quote$+" ,"+Quote$+"t9"+Quote$+" "
        CASE 10 : PRINT #N%, "splot "+Quote$+"t1"+Quote$+" lw 2,"+Quote$+"t2"+Quote$+" ,"+Quote$+"t3"+Quote$+" ,"+Quote$+"t4"+Quote$+"
,"+Quote$+"t5"+Quote$+" ,"+Quote$+"t6"+Quote$+"
                             Quote$+"t7"+Quote$+" ,"     +Quote$+"t8"+Quote$+" ,"+Quote$+"t9"+Quote$+" ,"+Quote$+"t10"+Quote$
    END SELECT

  END IF

  CLOSE #N%

  ProcID??? = SHELL(GnuPlotEXE$+" cmd3d.gp -") : CALL Delay(1##)

END SUB 'Plot3DbestProbeTrajectories()

'-----------

FUNCTION HasDAVGsaturated$(Nsteps&,j&,Np%,Nd%,M(),R(),DiagLength)

LOCAL A$

LOCAL k&

LOCAL SumOfDavg, DavgStepJ AS EXT

LOCAL DavgSatTOL AS EXT

    A$ = "NO"

    DavgSatTOL = 0.0005## 'tolerance for DAVG saturation

    IF j& < Nsteps& + 10 THEN GOTO ExitHasDAVGsaturated 'execute at least 10 steps after averaging interval before performing this check

    DavgStepJ = DavgThisStep(j&,Np%,Nd%,M(),R(),DiagLength)

    SumOfDavg = 0##

    FOR k& = j&-Nsteps&+1 TO j& 'check this step and previous (Nsteps&-1) steps

        SumOfDavg = SumOfDavg + DavgThisStep(k&,Np%,Nd%,M(),R(),DiagLength)

    NEXT k&

    IF ABS(SumOfDavg/Nsteps&-DavgStepJ) =< DavgSatTOL THEN A$ = "YES" 'saturation if (avg value - last value) are within TOL

ExitHasDAVGsaturated:

    HasDAVGsaturated$ = A$

END FUNCTION 'HasDAVGsaturated$()

'-----------

FUNCTION OscillationInDavg$(j&,Np%,Nd%,M(),R(),DiagLength)

LOCAL A$

LOCAL k&, NumSlopeChanges%

    A$ = "NO"

    NumSlopeChanges% = 0

    IF j& < 15 THEN GOTO ExitDavgOscillation 'wait at least 15 steps

    FOR k& = j&-10 TO j&-1 'check previous ten steps

        IF (DavgThisStep(k&,Np%,Nd%,M(),R(),DiagLength)-DavgThisStep(k&-1,Np%,Nd%,M(),R(),DiagLength))* _
           (DavgThisStep(k&+1,Np%,Nd%,M(),R(),DiagLength)-DavgThisStep(k&,Np%,Nd%,M(),R(),DiagLength)) < 0## THEN INCR NumSlopeChanges%

    NEXT j&

    IF NumSlopeChanges% >= 3 THEN A$ = "YES"
```

```
ExitDavgOscillation:
    OscillationInDavg$ = A$
END FUNCTION 'OscillationInDavg()
'------
FUNCTION DavgThisStep(j&,Np%,Nd%,M(),R(),DiagLength)
LOCAL BestFitness, TotalDistanceAllProbes, SumSQ AS EXT
LOCAL p%, k&, N%, i%, BestProbeNumber%, BestTimeStep&
'  ---------- Best Probe #, etc. -----------
    FOR k& = 0 TO j&
        BestFitness = M(1,k&)
        FOR p% = 1 TO Np%
            IF M(p%,k&) >= BestFitness THEN
                BestFitness = M(p%,k&) : BestProbeNumber% = p% : BestTimeStep& = k&
            END IF
        NEXT p% 'probe #
    NEXT k& 'time step
'  --------- Average Distance to Best Probe -----------
    TotalDistanceAllProbes = 0##
    FOR p% = 1 TO Np%
        SumSQ = 0##
        FOR i% = 1 TO Nd%
            SumSQ = SumSQ + (R(BestProbeNumber%,i%,BestTimeStep&)-R(p%,i%,j&))^2 'do not exclude p%=BestProbeNumber%(j&) from sum because it adds
zero
        NEXT i%
        TotalDistanceAllProbes = TotalDistanceAllProbes + SQR(SumSQ)
    NEXT p%
    DavgThisStep = TotalDistanceAllProbes/(DiagLength*(Np%-1)) 'but exclude best prove from average
END FUNCTION 'DavgThisStep()
'-----------
SUB
PlotBestFitnessEvolution(Nd%,Np%,LastStep&,G,DeltaT,Alpha,Beta,Frep,Mbest(),PlaceInitialProbes$,InitialAcceleration$,RepositionFactor$,FunctionName$,
Gamma)
LOCAL BestFitness(), GlobalBestFitness AS EXT
LOCAL PlotAnnotation$, PlotTitle$
LOCAL p%, j&, N%
    REDIM BestFitness(0 TO LastStep&)
    CALL
GetPlotAnnotation(PlotAnnotation$,Nd%,Np%,LastStep&,G,DeltaT,Alpha,Beta,Frep,Mbest(),PlaceInitialProbes$,InitialAcceleration$,RepositionFactor$,Funct
ionName$,Gamma)
    GlobalBestFitness = Mbest(1,0)
    FOR j& = 0 TO LastStep&
'        BestFitness(j&) = Mbest(1,j&) 'orig code 03-23-2010
        BestFitness(j&) = -1E4200 'added 03-23-2010
        FOR p% = 1 TO Np%
            IF Mbest(p%,j&) >= BestFitness(j&)   THEN BestFitness(j&)   = Mbest(p%,j&)
            IF Mbest(p%,j&) >= GlobalBestFitness THEN GlobalBestFitness = Mbest(p%,j&)
        NEXT p% 'probe #
    NEXT j& 'time step
    N% = FREEFILE
    OPEN "Fitness" FOR OUTPUT AS #N%
        FOR j& = 0 TO LastStep&
'            PRINT #N%, USING$("###### ##.######^^^^^",j&,BestFitness(j&))
            PRINT #N%, USING$("###### #######.########",j&,BestFitness(j&))
        NEXT j&

    CLOSE #N%

    PlotAnnotation$ = PlotAnnotation$ + "Best Fitness = " + REMOVE$(STR$(ROUND(GlobalBestFitness,8)),ANY" ")
    PlotTitle$ = "Best Fitness vs Time Step\n" + "[" + REMOVE$(STR$(Np%),ANY" ") + " probes, "+REMOVE$(STR$(LastStep&),ANY" ")+" time steps]"
    CALL CreateGNuplotINIfile(0.1##*Screenwidth&,0.1##*ScreenHeight&,0.6##*Screenwidth&,0.6##*ScreenHeight&)
    CALL TwoDplot("Fitness","Best Fitness","0.7","0.7","Time Step\n\n.","\n\nBest Fitness(X)", _
                  "","","","","","","","wgnuplot.exe"," with lines linewidth 2",PlotAnnotation$)
END SUB 'PlotBestFitnessEvolution()
'------
```

```
SUB
PlotAverageDistance(Nd%,Np%,LastStep&,G,DeltaT,Alpha,Beta,Frep,Mbest(),PlaceInitialProbes$,InitialAcceleration$,RepositionFactor$,FunctionName$,R(),D
iagLength,Gamma)

LOCAL Davg(), BestFitness(), TotalDistanceAllProbes, SumSQ AS EXT

LOCAL PlotAnnotation$, PlotTitle$

LOCAL p%, j&, N%, i%, BestProbeNumber%(), BestTimeStep&()

    REDIM Davg(0 TO LastStep&), BestFitness(0 TO LastStep&), BestProbeNumber%(0 TO LastStep&), BestTimeStep&(0 TO LastStep&)
    CALL
GetPlotAnnotation(PlotAnnotation$,Nd%,Np%,LastStep&,G,DeltaT,Alpha,Beta,Frep,Mbest(),PlaceInitialProbes$,InitialAcceleration$,RepositionFactor$,Funct
ionName$,Gamma)
' ---------- Best Probe #, etc. ----------

    FOR j& = 0 TO LastStep&

        BestFitness(j&) = Mbest(1,j&)

        FOR p% = 1 TO Np%

            IF Mbest(p%,j&) >= BestFitness(j&) THEN

                BestFitness(j&) = Mbest(p%,j&) : BestProbeNumber%(j&) = p% : BestTimeStep&(j&) = j& 'only probe number is used at this time, but
other data are computed for possible future use.

            END IF

            NEXT p% 'probe #

        NEXT j& 'time step

    N% = FREEFILE
' --------- Average Distance to Best Probe ----------

    FOR j& = 0 TO LastStep&

        TotalDistanceAllProbes = 0##

        FOR p% = 1 TO Np%

            SumSQ = 0##

            FOR i% = 1 TO Nd%

                SumSQ = SumSQ + (R(BestProbeNumber%(j&),i%)-R(p%,i%,j&))^2 'do not exclude p%=BestProbeNumber%(j&) from sum because it adds zero

            NEXT i%

            TotalDistanceAllProbes = TotalDistanceAllProbes + SQR(SumSQ)

        NEXT p%

        Davg(j&) = TotalDistanceAllProbes/(DiagLength*(Np%-1)) 'but exclude best prove from average

    NEXT j&
' ---------- Create Plot Data File ----------

    OPEN "Davg" FOR OUTPUT AS #N%

        FOR j& = 0 TO LastStep&

            PRINT #N%, USING$("###### #######.######",j&,Davg(j&))

        NEXT j&

    CLOSE #N%

    PlotTitle$ = "Average Distance of " + REMOVE$(STR$(Np%-1),ANY" ") + " Probes to Best Probe\nNormalized to Size of Decision Space\n" + _
                 "[" + REMOVE$(STR$(Np%),ANY" ") + " probes, " + REMOVE$(STR$(LastStep&),ANY" ") + " time steps]"

    CALL CreateGNuplotINIfile(0.2##*Screenwidth&,0.2##*ScreenHeight&,0.6##*Screenwidth&,0.6##*ScreenHeight&)

    CALL TwoDplot("Davg",PlotTitle$,"0.7","0.9","Time Step\n\n.",".",\n\n<b>/Ldiag", _
                  "","","","","","","wgnuplot.exe", with lines linewidth 2",PlotAnnotation$)

END SUB 'PlotAverageDistance()
'------

SUB
GetPlotAnnotation(PlotAnnotation$,Nd%,Np%,LastStep&,G,DeltaT,Alpha,Beta,Frep,Mbest(),PlaceInitialProbes$,InitialAcceleration$,RepositionFactor$,Funct
ionName$,Gamma)

LOCAL A$

    A$ = "" : IF PlaceInitialProbes$ = "UNIFORM ON-AXIS" AND Nd% > 1 THEN A$ = " ("+REMOVE$(STR$(Np%/Nd%),ANY" ") + "/axis)"

    PlotAnnotation$ = RunID$ + "\n" + _
                      FunctionName$ + " Function" + " ("+ FormatInteger$(Nd%) + "-D) \n"    +_
                      FormatInteger$(Np%) + " probes"         + A$ + "\n" +_
                      "G = " + FormatFP$(G,2)                 + "\n" +_
                      "Alpha = "      + FormatFP$(Alpha,1)    + "\n" +_
                      "Beta = "       + FormatFP$(Beta,1)     + "\n" +_
                      "DelT = "       + FormatFP$(DeltaT,1)   + "\n" +_
                      "Gamma = "      + FormatFP$(Gamma,3)    + "\n" +_
                      "Init Probes "  + PlaceInitialProbes$   + "\n" +_
                      "Init Accel "   + InitialAcceleration$  + "\n" +_
                      "Frep "         + RepositionFactor$ + "\n"

END SUB
'------

SUB
PlotBestProbeVsTimeStep(Nd%,Np%,LastStep&,G,DeltaT,Alpha,Beta,Frep,Mbest(),PlaceInitialProbes$,InitialAcceleration$,RepositionFactor$,FunctionName$,G
amma)

LOCAL BestFitness AS EXT

LOCAL PlotAnnotation$, PlotTitle$
```

```
LOCAL p%, j&, N%, BestProbeNumber%()

    REDIM BestProbeNumber%(0 TO LastStep&)

    CALL
GetPlotAnnotation(PlotAnnotation$,Nd%,Np%,LastStep&,G,DeltaT,Alpha,Beta,Frep,Mbest(),PlaceInitialProbes$,InitialAcceleration$,RepositionFactor$,Funct
ionName$,Gamma)

    FOR j& = 0 TO LastStep&

        Bestfitness = Mbest(1,j&)

        FOR p% = 1 TO Np%

            IF Mbest(p%,j&) >= BestFitness THEN

                BestFitness = Mbest(p%,j&) : BestProbeNumber%(j&) = p%

            END IF

        NEXT p% 'probe #

    NEXT j& 'time step

    N% = FREEFILE

    OPEN "Best Probe" FOR OUTPUT AS #N%

        FOR j& = 0 TO LastStep&

            PRINT #N%, USING$("###### #####",j&,BestProbeNumber%(j&))

        NEXT j&

    CLOSE #N%

    PlotTitle$ = "Best Probe Number vs Time Step\n" + "[" +REMOVE$(STR$(Np%),ANY" ") + " probes, " + REMOVE$(STR$(LastStep&),ANY" ") + " time steps]"

    CALL CreateGNuPlotINIfile(0.15#*ScreenWidth&,0.15#*ScreenHeight&,0.6#*ScreenWidth&,0.6#*ScreenHeight&)

'USAGE: CALL
TwoDplot(PlotFileName$,PlotTitle$,xCoord$,yCoord$,XaxisLabel$,YaxisLabel$,LogXaxis$,LogYaxis$,xMin$,xMax$,yMin$,yMax$,xTics$,yTics$,GnuPlotEXE$,LineT
ype$,Annotation$)

        CALL TwoDplot("Best Probe",PlotTitle$,"0.7","0.7","Time Step\n.",".\n\nBest Probe #","","","","","0",NoSpaces$(Np%+1,0),"","","wgnuplot.exe","
pt 8 ps .5 lw 1",PlotAnnotation$) 'pt, pointtype; ps, pointsize; lw, linewidth

END SUB 'PlotBestProbevsTimeStep()

'------

FUNCTION FormatInteger$(M%) : FormatInteger$ = REMOVE$(STR$(M%),ANY" ") : END FUNCTION

'------

FUNCTION FormatFP$(X,Ndigits%)

LOCAL A$

    IF X = 0## THEN

        A$ = "0." : GOTO ExitFormatFP

    END IF

    A$ = REMOVE$(STR$(ROUND(ABS(X),Ndigits%)),ANY" ")

    IF ABS(X) < 1## THEN

        IF X > 0## THEN

            A$ = "0" + A$

        ELSE

            A$ = "-0" + A$

        END IF

    ELSE

        IF X < 0## THEN A$ = "-" + A$

    END IF

ExitFormatFP:

    FormatFP$ = A$

END FUNCTION

'-----------

SUB InitialProbeDistribution(Np%,Nd%,Nt&,R(),PlaceInitialProbes$,Gamma)

LOCAL DeltaXi, Delx1, Delx2, Di AS EXT

LOCAL NumProbesPerDimension%, p%, i%, k%, NumX1points%, NumX2points%, x1pointNum%, x2pointNum%, A$

    SELECT CASE PlaceInitialProbes$

        CASE "UNIFORM ON-AXIS"

            IF Nd% > 1 THEN

                NumProbesPerDimension% = Np%\Nd% 'even #

            ELSE

                NumProbesPerDimension% = Np%

            END IF

            FOR i% = 1 TO Nd%
```

```
                    FOR p% = 1 TO Np%

                        R(p%,i%,0) = XiMin(i%) + Gamma$*(XiMax(i%)-XiMin(i%))

                    NEXT Np%

                NEXT i%

                FOR i% = 1 TO Nd%  'place probes probe line-by-probe line (i% is dimension number)

                    DeltaXi = (XiMax(i%)-XiMin(i%))/(NumProbesPerDimension%-1)

                    FOR k% = 1 TO NumProbesPerDimension

                        p% = k% + NumProbesPerDimension%*(i%-1) 'probe #

                        R(p%,i%,0) = XiMin(i%) + (k%-1)*DeltaXi

                    NEXT k%

                NEXT i%

            CASE "UNIFORM ON-DIAGONAL"

                FOR p% = 1 TO Np%

                    FOR i% = 1 TO Nd%

                        DeltaXi = (XiMax(i%)-XiMin(i%))/(Np%-1)

                        R(p%,i%,0) = XiMin(i%) + (p%-1)*DeltaXi

                    NEXT i%

                NEXT p%

            CASE "2D GRID"

                NumProbesPerDimension% = SQR(Np%) : NumX1points% = NumProbesPerDimension% : NumX2points% = NumX1points% 'broken down for possible future
use

                Delx1 = (XiMax(1)-XiMin(1))/(NumX1points%-1)

                Delx2 = (XiMax(2)-XiMin(2))/(NumX2points%-1)

                FOR x1pointNum% = 1 TO NumX1points%

                    FOR x2pointNum% = 1 TO NumX2points%

                        p% = NumX1points%*(x1pointNum%-1)+x2pointNum% 'probe #

                        R(p%,1,0) = XiMin(1) + Delx1*(x1pointNum%-1) 'x1 coord
                        R(p%,2,0) = XiMin(2) + Delx2*(x2pointNum%-1) 'x2 coord

                    NEXT x2pointNum%

                NEXT x1pointNum%

            CASE "RANDOM"

                FOR p% = 1 TO Np%

                    FOR i% = 1 TO Nd%

                        R(p%,i%,0) = XiMin(i%) + RandomNum(0##,1##)*(XiMax(i%)-XiMin(i%))

                    NEXT i%

                NEXT p%

        END SELECT

    END SUB 'InitialProbeDistribution()

    '------

    SUB
    ChangeRunParameters(NumProbesPerDimension%,Np%,Nd%,Nt&,G,Alpha,Beta,DeltaT,Frep,PlaceInitialProbes$,InitialAcceleration$,RepositionFactor$,FunctionNa
    me$) 'THIS PROCEDURE NOT USED

    LOCAL A$, DefaultValue$

        A$ = INPUTBOX$("# dimensions?","Change # Dimensions ("+FunctionName$+")",NoSpaces$(Nd%+0,0)) : Nd%    = VAL(A$) : IF Nd% < 1 OR Nd% > 500 THEN
    Nd% = 2

        IF Nd% > 1 THEN NumProbesPerDimension% = 2*((NumProbesPerDimension%+1)\2) 'require an even # probes on each probe line to avoid overlapping at
    origin (in symmetrical spaces at least...)

        IF Nd% = 1 THEN NumProbesPerDimension% = MAX(NumProbesPerDimension%,3)   'at least 3 probes on x-axis for 1-D functions

        Np% = NumProbesPerDimension%*Nd%

        A$ = INPUTBOX$("# time steps?","Change # Steps ("+FunctionName$+")",NoSpaces$(Nt&+0,0)) : Nt&    = VAL(A$) : IF Nt& < 3
    THEN Nt& = 50

        A$ = INPUTBOX$("Grav Const G?","Change G ("+FunctionName$+")",NoSpaces$(G,2))           : G    = VAL(A$) : IF G < -100## OR G > 100##
    THEN G = 2##

        A$ = INPUTBOX$("Alpha?","Change Alpha ("+FunctionName$+")",NoSpaces$(Alpha,2))          : Alpha = VAL(A$) : IF Alpha < -50## OR Alpha > 50##
    THEN Alpha = 2##

        A$ = INPUTBOX$("Beta?","Change Beta ("+FunctionName$+")",NoSpaces$(Beta,2))             : Beta = VAL(A$) : IF Beta  < -50## OR Beta  > 50##
    THEN Beta = 2##

        A$ = INPUTBOX$("Delta T","Change Delta-T ("+FunctionName$+")",NoSpaces$(DeltaT,2))       : DeltaT = VAL(A$) : IF DeltaT =< 0##
    THEN DeltaT = 1##

        A$ = INPUTBOX$("Frep [0-1]?","Change Frep ("+FunctionName$+")",NoSpaces$(Frep,3))        : Frep = VAL(A$) : IF Frep < 0##    OR Frep > 1##
    THEN Frep = 0.5##

    '  ----------- Initial Probe Distribution -----------

        SELECT CASE PlaceInitialProbes$
            CASE "UNIFORM ON-AXIS"     : DefaultValue$ = "1"
            CASE "UNIFORM ON-DIAGONAL" : DefaultValue$ = "2"
            CASE "2D GRID"             : DefaultValue$ = "3"
            CASE "RANDOM"              : DefaultValue$ = "4"
```

```
        END SELECT

        A$ = INPUTBOX$("Initial Probes?"+CHR$(13)+"1 - UNIFORM ON-AXIS"+CHR$(13)+"2 - UNIFORM ON-DIAGONAL"+CHR$(13)+"3 - 2D GRID"+CHR$(13)+"4 -
RANDOM","Initial Probe Distribution ("+FunctionName$+")",DefaultValue$)

        IF VAL(A$) < 1 OR VAL(A$) > 4 THEN A$ = "1"

        SELECT CASE VAL(A$)
            CASE 1 : PlaceInitialProbes$ = "UNIFORM ON-AXIS"
            CASE 2 : PlaceInitialProbes$ = "UNIFORM ON-DIAGONAL"
            CASE 3 : PlaceInitialProbes$ = "2D GRID"
            CASE 4 : PlaceInitialProbes$ = "RANDOM"
        END SELECT

        IF Nd% = 1  AND PlaceInitialProbes$ = "UNIFORM ON-DIAGONAL" THEN PlaceInitialProbes$ = "UNIFORM ON-AXIS" 'cannot do diagonal in 1-D space

        IF Nd% <> 2 AND PlaceInitialProbes$ = "2D GRID" THEN PlaceInitialProbes$ = "UNIFORM ON-AXIS" '2D grid is available only in 2 dimensions!

    '    ----------- Initial Acceleration -----------------

        SELECT CASE InitialAcceleration$
            CASE "ZERO"   : DefaultValue$ = "1"
            CASE "FIXED"  : DefaultValue$ = "2"
            CASE "RANDOM" : DefaultValue$ = "3"
        END SELECT

        A$ = INPUTBOX$("Initial Acceleration?"+CHR$(13)+"1 - ZERO"+CHR$(13)+"2 - FIXED"+CHR$(13)+"3 - RANDOM","Initial Acceleration
("+FunctionName$+")",DefaultValue$)

        IF VAL(A$) < 1 OR VAL(A$) > 3 THEN A$ = "1"

        SELECT CASE VAL(A$)
            CASE 1 : InitialAcceleration$ = "ZERO"
            CASE 2 : InitialAcceleration$ = "FIXED"
            CASE 3 : InitialAcceleration$ = "RANDOM"
        END SELECT

    '    ----------- Reposition Factor ---------------

        SELECT CASE RepositionFactor$
            CASE "FIXED"    : DefaultValue$ = "1"
            CASE "VARIABLE" : DefaultValue$ = "2"
            CASE "RANDOM"   : DefaultValue$ = "3"
        END SELECT

        A$ = INPUTBOX$("Reposition Factor?"+CHR$(13)+"1 - FIXED"+CHR$(13)+"2 - VARIABLE"+CHR$(13)+"3 - RANDOM","Retrieve Probes
("+FunctionName$+")",DefaultValue$)

        IF VAL(A$) < 1 OR VAL(A$) > 3 THEN A$ = "1"

        SELECT CASE VAL(A$)
            CASE 1 : RepositionFactor$ = "FIXED"
            CASE 2 : RepositionFactor$ = "VARIABLE"
            CASE 3 : RepositionFactor$ = "RANDOM"
        END SELECT

END SUB 'ChangeRunParameters()

'------

FUNCTION NoSpaces$(X,NumDigits%) :  NoSpaces$ = REMOVE$(STR$(X,NumDigits%),ANY" ") : END FUNCTION

'-----------

FUNCTION TerminateNowForSaturation$(j&,Nd%,Np%,Nt&,G,DeltaT,Alpha,Beta,R(),A(),M())

LOCAL A$, i&, p%, NumStepsForAveraging&

LOCAL BestFitness, AvgFitness, FitnessTOL AS EXT 'terminate if avg fitness does not change over NumStepsForAveraging& time steps

    FitnessTOL = 0.00001## : NumStepsForAveraging& = 10

    A$ = "NO"

    IF j& >= NumStepsForAveraging+10 THEN 'wait until step 10 to start checking for fitness saturation

        AvgFitness = 0##

        FOR i& = j&-NumStepsForAveraging&+1 TO j& 'avg fitness over current step & previous NumStepsForAveraging&-1 steps

            BestFitness = M(1,i&)

            FOR p% = 1 TO Np%

                IF M(p%,i&) >= BestFitness THEN BestFitness = M(p%,i&)

            NEXT p%

            AvgFitness = AvgFitness + BestFitness

        NEXT i&

        AvgFitness = AvgFitness/NumStepsForAveraging&

        IF ABS(AvgFitness-BestFitness) < FitnessTOL THEN A$ = "YES" 'compare avg fitness to best fitness at this step
    END IF

    TerminateNowForSaturation$ = A$

END FUNCTION 'TerminateNowForSaturation$()

'-----------

FUNCTION MagVector(V(),N%) 'returns magnitude of Nx1 column vector V

LOCAL SumSq AS EXT

LOCAL i%

    SumSQ = 0## : FOR i% = 1 TO N% : SumSQ = SumSQ + V(i%)^2 : NEXT i% : MagVector = SQR(SumSq)

END FUNCTION 'MagVector()

'---

FUNCTION UnitStep(X)
```

```
LOCAL Z AS EXT

    IF X < 0## THEN

        Z = 0##

    ELSE

        Z = 1##

    END IF

    UnitStep = Z

END FUNCTION 'UnitStep()
'---

SUB Plot1Dfunction(FunctionName$,R()) 'plots 1D function on-screen

LOCAL NumPoints%, i%, N%

LOCAL DeltaX, X AS EXT

    NumPoints% = 32001

    DeltaX = (XiMax(1)-XiMin(1))/(NumPoints%-1)

    N% = FREEFILE

    SELECT CASE FunctionName$

        CASE "ParrottF4" 'PARROTT F4 FUNCTION

            OPEN "ParrottF4" FOR OUTPUT AS #N%

                FOR i% = 1 TO NumPoints%

                    R(1,1,0) = XiMin(1) + (i%-1)*DeltaX

                    PRINT #N%, USING$("#.##### #.#####",R(1,1,0),ParrottF4(R(),1,1,0))

                NEXT i%

            CLOSE #N%

            CALL CreateGNUplotINIfile(0.2##*ScreenWidth&,0.2##*ScreenHeight&,0.6##*ScreenWidth&,0.6##*ScreenHeight&)

            CALL TwoDplot("ParrottF4","Parrott F4 Function","0.7","0.7","X\n\n.",".\n\nParrott F4(X)","","","0","1","0","1","","","wgnuplot.exe"," with lines linewidth 2","")

    END SELECT

END SUB
'------

SUB CLEANUP 'probe coordinate plot files

    IF DIR$("P1") <> "" THEN KILL "P1"
    IF DIR$("P2") <> "" THEN KILL "P2"
    IF DIR$("P3") <> "" THEN KILL "P3"
    IF DIR$("P4") <> "" THEN KILL "P4"
    IF DIR$("P5") <> "" THEN KILL "P5"
    IF DIR$("P6") <> "" THEN KILL "P6"
    IF DIR$("P7") <> "" THEN KILL "P7"
    IF DIR$("P8") <> "" THEN KILL "P8"
    IF DIR$("P9") <> "" THEN KILL "P9"
    IF DIR$("P10") <> "" THEN KILL "P10"
    IF DIR$("P11") <> "" THEN KILL "P11"
    IF DIR$("P12") <> "" THEN KILL "P12"
    IF DIR$("P13") <> "" THEN KILL "P13"
    IF DIR$("P14") <> "" THEN KILL "P14"
    IF DIR$("P15") <> "" THEN KILL "P15"

END SUB
'------

SUB Plot2Dfunction(FunctionName$,R())

LOCAL A$

LOCAL NumPoints%, i%, k%, N%

LOCAL Del1, Del2, Z AS EXT

    SELECT CASE FunctionName$

        CASE "LD_MONO" : Numpoints% = 5

        CASE "PBM_1","PBM_2","PBM_3","PBM_4","PBM_5" : NumPoints% = 25

        CASE ELSE : NumPoints% = 100

    END SELECT

    N% = FREEFILE : OPEN "TwoDplot.DAT" FOR OUTPUT AS #N%

    Del1 = (XiMax(1)-XiMin(1))/(NumPoints%-1) : Del2 = (XiMax(2)-XiMin(2))/(NumPoints%-1)

    FOR i% = 1 TO NumPoints%

        R(1,1,0) = XiMin(1) + (i%-1)*Del1 'x1 value

        FOR k% = 1 TO NumPoints%

            R(1,2,0) = XiMin(2) + (k%-1)*Del2 'x2 value

            Z = ObjectiveFunction(R(),2,1,0,FunctionName$)

            PRINT #N%, USING$("#####.##### #####.##### ######.#####^^^^",R(1,1,0),R(1,2,0),Z)

        NEXT k%

        PRINT #N%, ""
```

```
        NEXT i%

    CLOSE #N%

    CALL CreateGNUplotINIfile(0.1##*ScreenWidth&,0.1##*ScreenHeight&,0.6##*ScreenWidth&,0.6##*ScreenHeight&)

    A$ = "" : IF INSTR(FunctionName$,"PBM_") > 0 THEN A$ = "Coarse "

    CALL ThreeDplot2("TwoDplot.DAT",A$+"Plot of "+FunctionName$+" Function","","0.6","0.6","1.2", _
                    "x1","x2","z=F(x1,x2)","","","","wgnuplot.exe","","","","","")

END SUB
'------

    SUB TwoDplot3curves(NumCurves%,PlotFileName1$,PlotFileName2$,PlotFileName3$,PlotTitle$,Annotation$,xCoord$,yCoord$,XaxisLabel$,YaxisLabel$, _
                        LogXaxis$,LogYaxis$,xMin$,xMax$,yMin$,yMax$,xTics$,yTics$,GnuPlotEXE$)

        LOCAL N%

        LOCAL LineSize$

        LineSize$ = "2"

        N% = FREEFILE

        OPEN "cmd2d.gp" FOR OUTPUT AS #N%

            IF LogXaxis$ = "YES" AND LogYaxis$ = "NO"  THEN PRINT #N%, "set logscale x"
            IF LogXaxis$ = "NO"  AND LogYaxis$ = "YES" THEN PRINT #N%, "set logscale y"
            IF LogXaxis$ = "YES" AND LogYaxis$ = "YES" THEN PRINT #N%, "set logscale xy"

            IF xMin$ <> "" AND xMax$ <> "" THEN PRINT #N%, "set xrange ["+xMin$+":"+xMax$+"]"

            IF yMin$ <> "" AND yMax$ <> "" THEN  PRINT #N%, "set yrange ["+yMin$+":"+yMax$+"]"

            PRINT #N%, "set label "+Quote$+Annotation$+Quote$+" at graph "+xCoord$+","+yCoord$
            PRINT #N%, "set grid xtics"
            PRINT #N%, "set grid ytics"
            PRINT #N%, "set xtics "+xTics$
            PRINT #N%, "set ytics "+yTics$
            PRINT #N%, "set grid mxtics"
            PRINT #N%, "set grid mytics"
            PRINT #N%, "set title  "+Quote$+PlotTitle$+Quote$
            PRINT #N%, "set xlabel "+Quote$+XaxisLabel$+Quote$
            PRINT #N%, "set ylabel "+Quote$+YaxisLabel$+Quote$

            SELECT CASE NumCurves%

            CASE 1
            PRINT #N%, "plot " + Quote$ + PlotFileName1$ + Quote$ + " with lines linewidth " + LineSize$

            CASE 2
            PRINT #N%, "plot " + Quote$ + PlotFileName1$ + Quote$ + " with lines linewidth " + LineSize$+", " + _
                            Quote$ + PlotFileName2$ + Quote$ + " with lines linewidth " + LineSize$
            CASE 3
            PRINT #N%, "plot " + Quote$ + PlotFileName1$ + Quote$ + " with lines linewidth " + LineSize$+", " + _
                            Quote$ + PlotFileName2$ + Quote$ + " with lines linewidth " + LineSize$+", " + _
                            Quote$ + PlotFileName3$ + Quote$ + " with lines linewidth " + LineSize$
            END SELECT

        CLOSE #N%

        SHELL(GnuPlotEXE$+" cmd2d.gp -")

        CALL Delay(1##)

    END SUB 'TwoDplot3Curves()
'---

FUNCTION Fibonacci&&(N%) 'RETURNS Nth FIBONACCI NUMBER

LOCAL i%, Fn&&, Fn1&&, Fn2&&

LOCAL A$

    IF N% > 91 OR N% < 0 THEN

        MSGBOX("ERROR!  Fibonacci argument"+STR$(N%)+" > 91.  Out of range or < 0...") : EXIT FUNCTION

    END IF

    SELECT CASE N%

        CASE 0: Fn&& = 1

        CASE ELSE

            Fn&& = 0 : Fn2&& = 1 : i% = 0

            FOR i% = 1 TO N%

                Fn&& = Fn1&& + Fn2&&

                Fn1&& = Fn2&&

                Fn2&& = Fn&&

            NEXT i% 'LOOP

    END SELECT

    Fibonacci&& = Fn&&

END FUNCTION 'Fibonacci&&()
'----------

FUNCTION RandomNum(a,b) 'Returns random number X, a=< X < b.

    RandomNum = a + (b-a)*RND

END FUNCTION 'RandomNum()
'----------
```

```
FUNCTION GaussianDeviate(Mu,Sigma) 'returns NORMAL (Gaussian) random deviate with mean Mu and standard deviation Sigma (variance = Sigma^2)

'Refs: (1) Press, W.H., Flannery, B.P., Teukolsky, S.A., and Vetterling, W.T., "Numerical Recipes: The Art of Scientific Computing,"
'            §7.2, Cambridge University Press, Cambridge, UK, 1986.
'      (2) Shinzato, T., "Box Muller Method," 2007, http://www.sp.dis.titech.ac.jp/~shinzato/boxmuller.pdf

LOCAL s, t, Z AS EXT

    s = RND : t = RND

    Z = Mu + Sigma*SQR(-2##*LOG(s))*COS(TwoPi*t)

    GaussianDeviate = Z

END FUNCTION 'GaussianDeviate()

'-----------

    SUB ContourPlot(PlotFileName$,PlotTitle$,Annotation$,xCoord$,yCoord$,zCoord$, _
                    XaxisLabel$,YaxisLabel$,ZaxisLabel$,zMin$,zMax$,GnuPlotEXE$,A$)

        LOCAL N%

        N% = FREEFILE

        OPEN "cmd3d.gp" FOR OUTPUT AS #N%

            PRINT #N%, "show surface"
            PRINT #N%, "set hidden3d"
            IF zMin$ <> "" AND zMax$ <> "" THEN  PRINT #N%, "set zrange ["+zMin$+":"+zMax$+"]"
            PRINT #N%, "set label "+Quote$+AnnoTation$+Quote$+" at graph "+xCoord$+","+yCoord$+","+zCoord$
            PRINT #N%, "show label"
            PRINT #N%, "set grid xtics ytics ztics"
            PRINT #N%, "show grid"
            PRINT #N%, "set title "+Quote$+PlotTitle$+Quote$
            PRINT #N%, "set xlabel "+Quote$+XaxisLabel$+Quote$
            PRINT #N%, "set ylabel "+Quote$+YaxisLabel$+Quote$
            PRINT #N%, "set zlabel "+Quote$+ZaxisLabel$+Quote$
            PRINT #N%, "splot "+Quote$+PlotFileName$+Quote$+A$  '" notitle with linespoints' 'A$'" notitle with lines"
        CLOSE #N%

        SHELL(GnuPlotEXE$+" cmd3d.gp -")

    END SUB 'ContourPlot()

'---

    SUB ThreeDplot(PlotFileName$,PlotTitle$,Annotation$,xCoord$,yCoord$,zCoord$, _
                   XaxisLabel$,YaxisLabel$,ZaxisLabel$,zMin$,zMax$,GnuPlotEXE$,A$)

        LOCAL N%, ProcessID???

        N% = FREEFILE

        OPEN "cmd3d.gp" FOR OUTPUT AS #N%

            PRINT #N%, "set pm3d"
            PRINT #N%, "show pm3d"
            IF zMin$ <> "" AND zMax$ <> "" THEN  PRINT #N%, "set zrange ["+zMin$+":"+zMax$+"]"
            PRINT #N%, "set label "+Quote$+AnnoTation$+Quote$+" at graph "+xCoord$+","+yCoord$+","+zCoord$
            PRINT #N%, "show label"
            PRINT #N%, "set grid xtics ytics ztics"
            PRINT #N%, "show grid"
            PRINT #N%, "set title "+Quote$+PlotTitle$+Quote$
            PRINT #N%, "set xlabel "+Quote$+XaxisLabel$+Quote$
            PRINT #N%, "set ylabel "+Quote$+YaxisLabel$+Quote$
            PRINT #N%, "set zlabel "+Quote$+ZaxisLabel$+Quote$
            PRINT #N%, "splot "+Quote$+PlotFileName$+Quote$+A$+" notitle"' with lines"
        CLOSE #N%

        SHELL(GnuPlotEXE$+" cmd3d.gp -") : CALL Delay(1##)

    END SUB 'ThreeDplot()

'---

    SUB ThreeDplot2(PlotFileName$,PlotTitle$,Annotation$,xCoord$,yCoord$,zCoord$, _
                    XaxisLabel$,YaxisLabel$,ZaxisLabel$,zMin$,zMax$,A$,xStart$,xStop$,yStart$,yStop$)

        LOCAL N%

        N% = FREEFILE

        OPEN "cmd3d.gp" FOR OUTPUT AS #N%

            PRINT #N%, "set pm3d"
            PRINT #N%, "show pm3d"
            PRINT #N%, "set hidden3d"
            PRINT #N%, "set view 45, 45, 1, 1"

            IF zMin$ <> "" AND zMax$ <> "" THEN  PRINT #N%, "set zrange ["+zMin$+":"+zMax$+"]"

            PRINT #N%, "set xrange [" + xStart$ + ":" + xStop$ + "]"
            PRINT #N%, "set yrange [" + yStart$ + ":" + yStop$ + "]"

            PRINT #N%, "set label "    + Quote$  + AnnoTation$ + Quote$+" at graph "+xCoord$+","+yCoord$+","+zCoord$
            PRINT #N%, "show label"
            PRINT #N%, "set grid xtics ytics ztics"
            PRINT #N%, "show grid"
            PRINT #N%, "set title "  + Quote$+PlotTitle$    + Quote$
            PRINT #N%, "set xlabel " + Quote$+XaxisLabel$   + Quote$
            PRINT #N%, "set ylabel " + Quote$+YaxisLabel$   + Quote$
            PRINT #N%, "set zlabel " + Quote$+ZaxisLabel$   + Quote$
            PRINT #N%, "splot "      + Quote$+PlotFileName$ + Quote$ + A$ + " notitle with lines"
        CLOSE #N%

        SHELL(GnuPlotEXE$+" cmd3d.gp -")

    END SUB 'ThreeDplot2()

'---

    SUB TwoDplot2Curves(PlotFileName1$,PlotFileName2$,PlotTitle$,Annotation$,xCoord$,yCoord$,XaxisLabel$,YaxisLabel$, _
                        LogXaxis$,LogYaxis$,xMin$,xMax$,yMin$,yMax$,xTics$,yTics$,GnuPlotEXE$,LineSize)

        LOCAL N%, ProcessID???
```

```
        N% = FREEFILE

        OPEN "cmd2d.gp" FOR OUTPUT AS #N%
            'print #N%, "set output "+Quote$+"test.plt"+Quote$ 'tried this 3/11/06, didn't work...

            IF LogXaxis$ = "YES" AND LogYaxis$ = "NO"  THEN PRINT #N%, "set logscale x"
            IF LogXaxis$ = "NO"  AND LogYaxis$ = "YES" THEN PRINT #N%, "set logscale y"
            IF LogXaxis$ = "YES" AND LogYaxis$ = "YES" THEN PRINT #N%, "set logscale xy"

            IF xMin$ <> "" AND xMax$ <> "" THEN  PRINT #N%, "set xrange ["+xMin$+":"+xMax$+"]"

            IF yMin$ <> "" AND yMax$ <> "" THEN  PRINT #N%, "set yrange ["+yMin$+":"+yMax$+"]"

            PRINT #N%, "set label "+Quote$+Annotation$+Quote$+" at graph "+xCoord$+","+yCoord$
            PRINT #N%, "set grid xtics"
            PRINT #N%, "set grid ytics"
            PRINT #N%, "set xtics "+xTics$
            PRINT #N%, "set ytics "+yTics$
            PRINT #N%, "set grid mxtics"
            PRINT #N%, "set grid mytics"
            PRINT #N%, "set title "+Quote$+PlotTitle$+Quote$
            PRINT #N%, "set xlabel "+Quote$+XaxisLabel$+Quote$
            PRINT #N%, "set ylabel "+Quote$+YaxisLabel$+Quote$

            PRINT #N%, "plot "+Quote$+PlotFileName1$+Quote$+" with lines linewidth "+REMOVE$(STR$(LineSize),ANY) +" "+","+_
                       Quote$+PlotFileName2$+Quote$+" with points pointsize 0.05"+REMOVE$(STR$(LineSize),ANY) ")"

        CLOSE #N%

        ProcessID??? = SHELL(GnuPlotEXE$+" cmd2d.gp -") : CALL Delay(1##)

    END SUB 'TwoDplot2Curves()

'---

    SUB Probe2Dplots(ProbePlotsFileList$,PlotTitle$,Annotation$,xCoord$,yCoord$,XaxisLabel$,YaxisLabel$, _
                     LogXaxis$,LogYaxis$,xMin$,xMax$,yMin$,yMax$,xTics$,yTics$,GnuPlotEXE$)

        LOCAL N%, ProcessID???

        N% = FREEFILE

        OPEN "cmd2d.gp" FOR OUTPUT AS #N%

            IF LogXaxis$ = "YES" AND LogYaxis$ = "NO"  THEN PRINT #N%, "set logscale x"
            IF LogXaxis$ = "NO"  AND LogYaxis$ = "YES" THEN PRINT #N%, "set logscale y"
            IF LogXaxis$ = "YES" AND LogYaxis$ = "YES" THEN PRINT #N%, "set logscale xy"

            IF xMin$ <> "" AND xMax$ <> "" THEN  PRINT #N%, "set xrange ["+xMin$+":"+xMax$+"]"

            IF yMin$ <> "" AND yMax$ <> "" THEN  PRINT #N%, "set yrange ["+yMin$+":"+yMax$+"]"

            PRINT #N%, "set label "+Quote$+Annotation$+Quote$+" at graph "+xCoord$+","+yCoord$
            PRINT #N%, "set grid xtics"
            PRINT #N%, "set grid ytics"
            PRINT #N%, "set xtics "+xTics$
            PRINT #N%, "set ytics "+yTics$
            PRINT #N%, "set grid mxtics"
            PRINT #N%, "set grid mytics"
            PRINT #N%, "set title "+Quote$+PlotTitle$+Quote$
            PRINT #N%, "set xlabel "+Quote$+XaxisLabel$+Quote$
            PRINT #N%, "set ylabel "+Quote$+YaxisLabel$+Quote$

            PRINT #N%, ProbePlotsFileList$

        CLOSE #N%

        ProcessID??? = SHELL(GnuPlotEXE$+" cmd2d.gp -") : CALL Delay(1##)

    END SUB 'Probe2Dplots()

'---

SUB Show2Dprobes(R(),Np%,Nt&,j&,Frep,BestFitness,BestProbeNumber%,BestTimeStep&,FunctionName$,RepositionFactor$,Gamma)

    LOCAL N%, p%

    LOCAL A$, PlotFileName$, PlotTitle$, Symbols$

    LOCAL xMin$, xMax$, yMin$, yMax$

    LOCAL s1, s2, s3, s4 AS EXT

    PlotFileName$ = "Probes("+REMOVE$(STR$(j&)," ")+")"

    IF j& > 0 THEN 'PLOT PROBES AT THIS TIME STEP

        PlotTitle$ = "\nLOCATIONS OF "+REMOVE$(STR$(Np%),ANY) + " PROBES AT TIME STEP" + STR$(j&) + " / " + REMOVE$(STR$(Nt&),ANY) +") + "\n" +_
             "Fitness = "+REMOVE$(STR$(ROUND(BestFitness,3)),ANY) +") + ", Probe #" + REMOVE$(STR$(BestProbeNumber%),ANY) +") + " at Step #"+_
REMOVE$(STR$(BestTimeStep&),ANY) +") +_
                    " [Frep = "+REMOVE$(STR$(Frep,4),ANY) + " " + RepositionFactor$ + "]\n"

    ELSE 'PLOT INITIAL PROBE DISTRIBUTION

        PlotTitle$ = "\nLOCATIONS OF "+REMOVE$(STR$(Np%),ANY) + " INITIAL PROBES FOR " + FunctionName$ + " FUNCTION\n[gamma =
"+STR$(ROUND(Gamma,3))+"]\n"

    END IF

    N% = FREEFILE : OPEN PlotFileName$ FOR OUTPUT AS #N%

        FOR p% = 1 TO Np% : PRINT #N%, USING$("######.####    ######.####",R(p%,1,j&),R(p%,2,j&)) : NEXT p%

    CLOSE #N%

    s1 = 1.1## : s2 = 1.1## : s3 = 1.1## : s4 = 1.1## 'expand plots axes by 10%

    IF XiMin(1) > 0## THEN s1 = 0.9##
    IF XiMax(1) < 0## THEN s2 = 0.9##
    IF XiMin(2) > 0## THEN s3 = 0.9##
    IF XiMax(2) < 0## THEN s4 = 0.9##

    xMin$ = REMOVE$(STR$(s1*XiMin(1),2),ANY" ")
    xMax$ = REMOVE$(STR$(s2*XiMax(1),2),ANY" ")
    yMin$ = REMOVE$(STR$(s3*XiMin(2),2),ANY" ")
    yMax$ = REMOVE$(STR$(s4*XiMax(2),2),ANY" ")
```

```
    CALL TwoDplot(PlotFileName$,PlotTitle$,"0.6","0.7","x1\n\n","\nx2","NO","NO",xMin$,xMax$,yMin$,yMax$,"5","5","wgnuplot.exe"," pointsize 1
linewidth 2","")

    KILL PlotFileName$ 'erase plot data file after probes have been displayed

END SUB 'Show2Dprobes()

'---

SUB Show3Dprobes(R(),Np%,Nd%,Nt&,j&,Frep,BestFitness,BestProbeNumber%,BestTimeStep&,FunctionName$,RepositionFactor$,Gamma)

    LOCAL N%, p%, PlotWindowULC_X%, PlotWindowULC_Y%, PlotWindowwidth%, PlotWindowHeight%, PlotWindowOffset%

    LOCAL A$, PlotFileName$, PlotTitle$, Symbols$, Annotation$

    LOCAL xMin$, xMax$, yMin$, yMax$, zMin$, zMax$

    LOCAL s1, s2, s3, s4, s5, s6 AS EXT

    PlotFileName$ = "Probes("+REMOVE$(STR$(j&),ANY" ")+")"

    IF j& > 0 THEN 'PLOT PROBES AT THIS TIME STEP

        PlotTitle$ = "\nLOCATIONS OF "+REMOVE$(STR$(Np%),ANY" ") + " PROBES AT TIME STEP" + STR$(j&) + " / " + REMOVE$(STR$(Nt&),ANY" ") + "\n" +_
                     "Fitness = "+REMOVE$(STR$(ROUND(BestFitness,3)),ANY" ") + ", Probe #" + REMOVE$(STR$(BestProbeNumber%),ANY" ") + " at Step #"+
REMOVE$(STR$(BestTimeStep&),ANY" ") +_
                     "  [Frep = "+REMOVE$(STR$(Frep,4),ANY" ") + " " + RepositionFactor$ + "]\n"

    ELSE 'PLOT INITIAL PROBE DISTRIBUTION

        A$ = "" : IF Gamma > 0## AND Gamma < 1## THEN A$ = "0"

        PlotTitle$ = "\n"+REMOVE$(STR$(Np%),ANY" ") + "-PROBE IPD FOR FUNCTION " + FunctionName$ + ", GAMMA = "+A$+REMOVE$(STR$(ROUND(Gamma,3)),ANY"
")

    END IF

'   --------------- Probe Coordinates -----------------

    N% = FREEFILE : OPEN PlotFileName$ FOR OUTPUT AS #N%

        PRINT #N%, USING$("#####.#####     #####.#####     #####.#####",R(1,1,j&),R(1,2,j&),R(1,3,j&)) 'This line repeats Probe #1's coordinates.
It's necessary
                                                                                       'to deal with a plotting artifact in Gnuplot!
        FOR p% = 1 TO Np% : PRINT #N%, USING$("#####.#####     #####.#####     #####.#####",R(p%,1,j&),R(p%,2,j&),R(p%,3,j&)) : NEXT p%

    CLOSE #N%

'   ------------- Principal Diagonal -------------

    N% = FREEFILE : OPEN "diag" FOR OUTPUT AS #N%

        PRINT #N%, USING$("#####.#####     #####.#####     #####.#####",XiMin(1),XiMin(2),XiMin(3))
        PRINT #N%, ""
        PRINT #N%, USING$("#####.#####     #####.#####     #####.#####",XiMax(1),XiMax(2),XiMax(3))

    CLOSE #N%

'   ------------------- Probe Line #1 ------------------

    N% = FREEFILE : OPEN "probeline1" FOR OUTPUT AS #N%

        PRINT #N%, USING$("#####.#####     #####.#####     #####.#####",R(1,1,j&),R(1,2,j&),R(1,3,j&))
        PRINT #N%, ""
        PRINT #N%, USING$("#####.#####     #####.#####     #####.#####",R(Np%/Nd%,1,j&),R(Np%/Nd%,2,j&),R(Np%/Nd%,3,j&))

    CLOSE #N%

'   ------------------- Probe Line #2 ------------------

    N% = FREEFILE : OPEN "probeline2" FOR OUTPUT AS #N%

        PRINT #N%, USING$("#####.#####     #####.#####     #####.#####",R(1+Np%/Nd%,1,j&),R(1+Np%/Nd%,2,j&),R(1+Np%/Nd%,3,j&))
        PRINT #N%, ""
        PRINT #N%, USING$("#####.#####     #####.#####     #####.#####",R(2*Np%/Nd%,1,j&),R(2*Np%/Nd%,2,j&),R(2*Np%/Nd%,3,j&))

    CLOSE #N%

'   ------------------- Probe Line #3 ------------------

    N% = FREEFILE : OPEN "probeline3" FOR OUTPUT AS #N%

        PRINT #N%, USING$("#####.#####     #####.#####     #####.#####",R(1+2*Np%/Nd%,1,j&),R(1+2*Np%/Nd%,2,j&),R(1+2*Np%/Nd%,3,j&))
        PRINT #N%, ""
        PRINT #N%, USING$("#####.#####     #####.#####     #####.#####",R(3*Np%/Nd%,1,j&),R(3*Np%/Nd%,2,j&),R(3*Np%/Nd%,3,j&))

    CLOSE #N%

'   -------- RE-PLOT PROBE #1 BECAUSE OF SOME ARTIFACT THAT DROPS IT FROM PROBE LINE #1 ????? ------------------

    N% = FREEFILE : OPEN "probe1" FOR OUTPUT AS #N%

        PRINT #N%, USING$("#####.#####     #####.#####     #####.#####",R(1,1,j&),R(1,2,j&),R(1,3,j&))
        PRINT #N%, USING$("#####.#####     #####.#####     #####.#####",R(1,1,j&),R(1,2,j&),R(1,3,j&))

    CLOSE #N%

    s1 = 1.1## : s2 = s1 : s3 = s1 : s4 = s1 : s5 = s1 : s6 = s1 'expand plots axes by 10%

    IF XiMin(1) > 0## THEN s1 = 0.9##
    IF XiMax(1) < 0## THEN s2 = 0.9##
    IF XiMin(2) > 0## THEN s3 = 0.9##
    IF XiMax(2) < 0## THEN s4 = 0.9##
    IF XiMin(3) > 0## THEN s5 = 0.9##
    IF XiMax(3) < 0## THEN s6 = 0.9##

    xMin$ = REMOVE$(STR$(s1*XiMin(1),2),ANY" ")
    xMax$ = REMOVE$(STR$(s2*XiMax(1),2),ANY" ")
    yMin$ = REMOVE$(STR$(s3*XiMin(2),2),ANY" ")
    yMax$ = REMOVE$(STR$(s4*XiMax(2),2),ANY" ")
    zMin$ = REMOVE$(STR$(s5*XiMin(2),2),ANY" ")
    zMax$ = REMOVE$(STR$(s6*XiMax(2),2),ANY" ")

'USAGE: CALL ThreeDplot3(PlotFileName$,PlotTitle$,Annotation$,xCoord$,yCoord$,zCoord$,_
'                        xaxisLabel$,yaxisLabel$,zaxisLabel$,zMin$,zMax$,GnuPlotEXE$,xStart$,xStop$,yStart$,yStop$)

    PlotWindowULC_X% = 50 : PlotWindowULC_Y% = 50 : PlotWindowwidth% = 1000 : PlotWindowHeight% = 800
```

```
        PlotwindowOffset% = 100*Gamma

        CALL CreateGNuplotINIfile(PlotwindowULC_X%+PlotwindowOffset%,PlotwindowULC_Y%+PlotwindowOffset%,PlotwindowwidtH%,PlotwindowHeight%)

        CALL ThreeDplot3(PlotFileName$,PlotTitle$,Annotation$,"0.6","0.7","0.8", _
                         "x1","x2","x3",zMin$,zMax$,"wgnuplot.exe",xMin$,xMax$,yMin$,yMax$)

        'KILL PlotFileName$ 'erase plot data file after probes have been displayed

END SUB 'Show3Dprobes()

'---

        SUB ThreeDplot3(PlotFileName$,PlotTitle$,Annotation$,xCoord$,yCoord$,zCoord$, _
                        XaxisLabel$,YaxisLabel$,ZaxisLabel$,zMin$,zMax$,GnuPlotEXE$,xStart$,xStop$,yStart$,yStop$)

            LOCAL N%, ProcID???

            N% = FREEFILE

            OPEN "cmd3d.gp" FOR OUTPUT AS #N%

                PRINT #N%, "set pm3d"
                PRINT #N%, "show pm3d"
                PRINT #N%, "set hidden3d"
'               PRINT #N%, "set view 45, 45, 1, 1"

                PRINT #N%, "set view 45, 60, 1, 1"

                IF zMin$ <> "" AND zMax$ <> "" THEN  PRINT #N%, "set zrange ["+zMin$+":"+zMax$+"]"

                PRINT #N%, "set xrange [" + xStart$ + ":" + xStop$ + "]"
                PRINT #N%, "set yrange [" + yStart$ + ":" + yStop$ + "]"

                PRINT #N%, "set label "  + Quote$ + AnnoTation$ + Quote$+" at graph "+xCoord$+","+yCoord$+","+zCoord$
                PRINT #N%, "show label"
                PRINT #N%, "set grid xtics ytics ztics"
                PRINT #N%, "show grid"
                PRINT #N%, "set title "  + Quote$+PlotTitle$   + Quote$
                PRINT #N%, "set xlabel " + Quote$+XaxisLabel$  + Quote$
                PRINT #N%, "set ylabel " + Quote$+YaxisLabel$  + Quote$
                PRINT #N%, "set zlabel " + Quote$+ZaxisLabel$  + Quote$
                PRINT #N%, "unset colorbox"
'               print #N%, "set style fill"

                PRINT #N%, "splot "      + Quote$+PlotFileName$ + Quote$ + " notitle lw 1 pt 8," _
                                         + Quote$ + "diag"      + Quote$ + " notitle w l," _
                                         + Quote$ + "probeline1" + Quote$ + " notitle w l," _
                                         + Quote$ + "probeline2" + Quote$ + " notitle w l," _
                                         + Quote$ + "probeline3" + Quote$ + " notitle w l"

            CLOSE #N%

            ProcID??? = SHELL(GnuPlotEXE$+" cmd3d.gp -")

            CALL Delay(1##)

        END SUB 'ThreeDplot3()

'----

        SUB TwoDplot(PlotFileName$,PlotTitle$,xCoord$,yCoord$,XaxisLabel$,YaxisLabel$, _
                     LogXaxis$,LogYaxis$,xMin$,xMax$,yMin$,yMax$,xTics$,yTics$,GnuPlotEXE$,LineType$,Annotation$)

            LOCAL N%, ProcessID???

            N% = FREEFILE

            OPEN "cmd2d.gp" FOR OUTPUT AS #N%

                IF LogXaxis$ = "YES" AND LogYaxis$ = "NO"  THEN PRINT #N%, "set logscale x"
                IF LogXaxis$ = "NO"  AND LogYaxis$ = "YES" THEN PRINT #N%, "set logscale y"
                IF LogXaxis$ = "YES" AND LogYaxis$ = "YES" THEN PRINT #N%, "set logscale xy"

                IF xMin$ <> "" AND xMax$ <> "" THEN  PRINT #N%, "set xrange ["+xMin$+":"+xMax$+"]"
                IF yMin$ <> "" AND yMax$ <> "" THEN  PRINT #N%, "set yrange ["+yMin$+":"+yMax$+"]"

                PRINT #N%, "set label "     + Quote$ + Annotation$ + Quote$ + " at graph " + xCoord$ + "," + yCoord$
                PRINT #N%, "set grid xtics " + XTics$
                PRINT #N%, "set grid ytics " + yTics$
                PRINT #N%, "set grid mxtics"
                PRINT #N%, "set grid mytics"
                PRINT #N%, "show grid"
                PRINT #N%, "set title "  + Quote$+PlotTitle$+Quote$
                PRINT #N%, "set xlabel " + Quote$+XaxisLabel$+Quote$
                PRINT #N%, "set ylabel " + Quote$+YaxisLabel$+Quote$

                PRINT #N%, "plot "+Quote$+PlotFileName$+Quote$+" notitle"+LineType$

            CLOSE #N%

            ProcessID??? = SHELL(GnuPlotEXE$+" cmd2d.gp -") : CALL Delay(1##)

        END SUB 'TwoDplot()

'-----

        SUB CreateGNuplotINIfile(PlotwindowULC_X%,PlotwindowULC_Y%,PlotwindowwidtH%,PlotwindowHeight%)

        LOCAL N%, WinPath$, A$, B$, WindowsDirectory$

        WinPath$ = UCASE$(ENVIRON$("Path"))'DIR$("C:\WINDOWS",23)

        DO
            B$ = A$

            A$ = EXTRACT$(WinPath$,";")

            WinPath$ = REMOVE$(WinPath$,A$+";")

            IF RIGHT$(A$,7) = "WINDOWS" OR A$ = B$ THEN EXIT LOOP

            IF RIGHT$(A$,5) = "WINNT"   OR A$ = B$ THEN EXIT LOOP

        LOOP

        WindowsDirectory$ = A$
```

```
        N% = FREEFILE

'   ----------- WGNUPLOT.INPUT FILE -----------
        OPEN WindowsDirectory$+"\wgnuplot.ini" FOR OUTPUT AS #N%

            PRINT #N%,"[WGNUPLOT]"
            PRINT #N%,"TextOrigin=0 0"
            PRINT #N%,"TextSize=640 150"
            PRINT #N%,"TextFont=Terminal,9"
            PRINT #N%,"GraphOrigin="+REMOVE$(STR$(PlotWindowULC_X%),ANY" ")+" "+REMOVE$(STR$(PlotWindowULC_Y%),ANY" ")
            PRINT #N%,"GraphSize="+REMOVE$(STR$(PlotWindowWidth%),ANY" ")+" "+REMOVE$(STR$(PlotWindowHeight%),ANY" ")
            PRINT #N%,"GraphFont=Arial,10"
            PRINT #N%,"GraphColor=1"
            PRINT #N%,"GraphToTop=1"
            PRINT #N%,"GraphBackground=255 255 255"
            PRINT #N%,"Border=0 0 0 0"
            PRINT #N%,"Axis=192 192 192 2 2"
            PRINT #N%,"Line1=0 0 255 0 0"
            PRINT #N%,"Line2=0 255 0 0 1"
            PRINT #N%,"Line3=255 0 0 0 2"
            PRINT #N%,"Line4=255 0 255 0 3"
            PRINT #N%,"Line5=0 0 128 0 4"

        CLOSE #N%

        END SUB 'CreateGNUplotINIfile()

'------

        SUB Delay(NumSecs)

            LOCAL StartTime, StopTime AS EXT

            StartTime = TIMER

            DO UNTIL (StopTime-StartTime) >= NumSecs

                StopTime = TIMER

            LOOP

        END SUB 'Delay()

'-----

SUB MathematicalConstants
    EulerConst  = 0.5772156649015328606606512##
    Pi          = 3.141592653589793238462643##
    Pi2         = Pi/2##
    Pi4         = Pi/4##
    TwoPi       = 2##*Pi
    FourPi      = 4##*Pi
    e           = 2.718281828459045235360287##
    Root2       = 1.414213562373095048##
END SUB

'-----

SUB AlphabetAndDigits
    Alphabet$   = "ABCDEFGHIJKLMNOPQRSTUVWXYZabcdefghijklmnopqrstuvwxyz"
    Digits$     = "0123456789"
    RunID$      = DATE$ + ". " + TIME$
END SUB

'------

SUB SpecialSymbols
    Quote$              = CHR$(34) 'Quotation mark "
    SpecialCharacters$ = "'(),#:;/_"
END SUB

'-----

SUB EMConstants
    Mu0  = 4E-7##*Pi      'hy/meter
    Eps0 = 8.854##*1E-12  'fd/meter
    c    = 2.998E8##      'velocity of light, 1##/SQR(Mu0*Eps0) 'meters/sec
    eta0 = SQR(Mu0/Eps0) 'impedance of free space, ohms
END SUB

'------

SUB ConversionFactors
    Rad2Deg      = 180##/Pi
    Deg2Rad      = 1##/Rad2Deg
    Feet2Meters  = 0.3048##
    Meters2Feet  = 1##/Feet2Meters
    Inches2Meters = 0.0254##
    Meters2Inches = 1##/Inches2Meters
    Miles2Meters = 1609.344##
    Meters2Miles = 1##/Miles2Meters
    NautMi2Meters = 1852##
    Meters2NautMi = 1##/NautMi2Meters
END SUB

'------

SUB ShowConstants 'puts up msgbox showing all constants

LOCAL A$

A$ = _
"Mathematical Constants:"+CHR$(13)+_
"Euler const="+STR$(EulerConst)+CHR$(13)+_
"Pi="+STR$(Pi)+CHR$(13)+_
"Pi/2="+STR$(Pi2)+CHR$(13)+_
"Pi/4="+STR$(Pi4)+CHR$(13)+_
"2Pi="+STR$(TwoPi)+CHR$(13)+_
"4Pi="+STR$(FourPi)+CHR$(13)+_
"e="+STR$(e)+CHR$(13)+_
"Sqr2="+STR$(Root2)+CHR$(13)+CHR$(13)+_
"Alphabet, Digits & Special Characters:"+CHR$(13)+_
"Alphabet="+Alphabet$+CHR$(13)+_
"Digits="+Digits$+CHR$(13)+_
"quote="+Quote$+CHR$(13)+_
"Spec chars="+SpecialCharacters$+CHR$(13)+CHR$(13)+_
```

```
"E&M Constants:"+CHR$(13)+_
"Mu0="+STR$(Mu0)+CHR$(13)+_
"Eps0="+STR$(Eps0)+CHR$(13)+_
"C="+STR$(c)+CHR$(13)+_
"Eta0="+STR$(eta0)+CHR$(13)+CHR$(13)+_
"Conversion Factors:"+CHR$(13)+_
"Rad2Deg="+STR$(Rad2Deg)+CHR$(13)+_
"Deg2Rad="+STR$(Deg2Rad)+CHR$(13)+_
"Ft2meters="+STR$(Feet2Meters)+CHR$(13)+_
"Meters2Ft="+STR$(Meters2Feet)+CHR$(13)+_
"Inches2Meters="+STR$(Inches2Meters)+CHR$(13)+_
"Meters2Inches="+STR$(Meters2Inches)+CHR$(13)+_
"Miles2Meters="+STR$(Miles2Meters)+CHR$(13)+_
"Meters2Miles="+STR$(Meters2Miles)+CHR$(13)+_
"NautMi2Meters="+STR$(NautMi2Meters)+CHR$(13)+_
"Meters2NautMi="+STR$(Meters2NautMi)+CHR$(13)+CHR$(13)

MSGBOX(A$)

END SUB

'------

SUB DisplayRmatrix(Np%,Nd%,Nt&,R())

LOCAL p%, i%, j&, A$

    A$ = "Position Vector Matrix R()"+CHR$(13)

    FOR p% = 1 TO Np%

        FOR i% = 1 TO Nd%

            FOR j& = 0 TO Nt&

                A$ = A$ + "R("+STR$(p%)+", "+STR$(i%)+", "+STR$(j&)+ " )="+STR$(R(p%,i%,j&)) + CHR$(13)

            NEXT j&

        NEXT i%

    NEXT p%

    MSGBOX(A$)

END SUB

'------

SUB DisplayRmatrixThisTimeStep(Np%,Nd%,j&,R(),Gamma)

LOCAL p%, i%, A$, B$

    A$ = "Position Vector Matrix R() at step "+STR$(j&)+", Gamma ="+STR$(Gamma)+":"+CHR$(13)+CHR$(13)

    FOR p% = 1 TO Np%

        A$ = A$ + "Probe#"+REMOVE$(STR$(p%),ANY" ")+": "

        B$ = ""

        FOR i% = 1 TO Nd%

            B$ = B$ + "  " + USING$("####.##",R(p%,i%,j&))

        NEXT i%

        A$ = A$ + B$ + CHR$(13)

    NEXT p%

    MSGBOX(A$)

END SUB

'------

SUB DisplayAmatrix(Np%,Nd%,Nt&,A())

LOCAL p%, i%, j&, A$

    A$ = "Acceleration Vector Matrix A()"+CHR$(13)

    FOR p% = 1 TO Np%

        FOR i% = 1 TO Nd%

            FOR j& = 0 TO Nt&

                A$ = A$ + "A("+STR$(p%)+", "+STR$(i%)+", "+STR$(j&)+ " )="+STR$(A(p%,i%,j&)) + CHR$(13)

            NEXT j&

        NEXT i%

    NEXT p%

    MSGBOX(A$)

END SUB

'------

SUB DisplayAmatrixThisTimeStep(Np%,Nd%,j&,A())

LOCAL p%, i%, A$

    A$ = "Acceleration matrix A() at step "+STR$(j&)+":"+CHR$(13)

    FOR p% = 1 TO Np%

        FOR i% = 1 TO Nd%

            A$ = A$ + "A("+STR$(p%)+", "+STR$(i%)+", "+STR$(j&)+ " )="+STR$(A(p%,i%,j&)) + CHR$(13)

        NEXT i%
```

```
        NEXT p%

        MSGBOX(A$)

END SUB
'------

SUB DisplayMmatrix(Np%,Nt&,M())

LOCAL p%, j&, A$

        A$ = "Fitness Matrix M()"+CHR$(13)

        FOR p% = 1 TO Np%

            FOR j& = 0 TO Nt&

                A$ = A$ + "M("+STR$(p%)+", "+STR$(j&)+ " ) ="+STR$(M(p%,j&)) + CHR$(13)

            NEXT j&

        NEXT p%

        MSGBOX(A$)

END SUB
'------

SUB DisplayMbestMatrix(Np%,Nt&,Mbest())

LOCAL p%, j&, A$

        A$ = "Np= "+STR$(Np%)+"   Nt="+STR$(Nt&)+CHR$(13)+"Fitness Matrix Mbest()"+CHR$(13)

        FOR p% = 1 TO Np%

            FOR j& = 0 TO Nt&

                A$ = A$ + "Mbest("+STR$(p%)+", "+STR$(j&)+ " ) ="+STR$(Mbest(p%,j&)) + CHR$(13)

            NEXT j&

        NEXT p%

        MSGBOX(A$)

END SUB
'------

SUB DisplayMmatrixThisTimeStep(Np%,j&,M())

LOCAL p%, A$

        A$ = "Fitness matrix M() at step "+STR$(j&)+":"+CHR$(13)

        FOR p% = 1 TO Np%

                A$ = A$ + "M("+STR$(p%)+", "+STR$(j&)+ " ) ="+STR$(M(p%,j&)) + CHR$(13)

        NEXT p%

        MSGBOX(A$)

END SUB
'------

SUB DisplayXiMinMax(Nd%,XiMin(),XiMax())

LOCAL i%, A$

        A$ = ""

        FOR i% = 1 TO Nd%

                A$ = A$ + "XiMin("+STR$(i%)+" ) = "+STR$(XiMin(i%))+"   XiMax("+STR$(i%)+" ) = "+STR$(XiMax(i%)) + CHR$(13)

        NEXT i%

        MSGBOX(A$)

END SUB
'------

SUB DisplayRunParameters2(FunctionName$,Nd%,Np%,Nt&,G,DeltaT,Alpha,Beta,Frep,PlaceInitialProbes$,InitialAcceleration$,RepositionFactor$)

LOCAL A$

        A$ = "Function = "+ FunctionName$+CHR$(13)+_
             "Nd = "+STR$(Nd%)+CHR$(13)+_
             "Np = "+STR$(Np%)+CHR$(13)+_
             "Nt = "+STR$(Nt&)+CHR$(13)+_
             "G  = "+STR$(G)+CHR$(13)+_
             "DeltaT = "+STR$(DeltaT)+CHR$(13)+_
             "Alpha = "+STR$(Alpha)+CHR$(13)+_
             "Beta  = "+STR$(Beta)+CHR$(13)+_
             "Frep  = "+STR$(Frep)+CHR$(13)+_
             "Init Probes: "+PlaceInitialProbes$+CHR$(13)+_
             "Init Accel:  "+InitialAcceleration$+CHR$(13)+_
             "Retrive Method: "+RepositionFactor$+CHR$(13)

        MSGBOX(A$)

END SUB
'------

SUB
Tabulate1DprobeCoordinates(Max1DprobesPlotted%,Nd%,Np%,LastStep&,G,DeltaT,Alpha,Beta,Frep,R(),M(),PlaceInitialProbes$,InitialAcceleration$,Reposition
Factor$,FunctionName$,Gamma)

LOCAL N%, ProbeNum%, FileHeader$, A$, B$, C$, D$, E$, F$, H$, StepNum&, FieldNumber% 'kludgy, yes, but it accomplishes its purpose...
```

```
            CALL
GetPlotAnnotation(FileHeader$,Nd%,Np%,LastStep&,G,DeltaT,Alpha,Beta,Frep,M(),PlaceInitialProbes$,InitialAcceleration$,RepositionFactor$,FunctionName$
,Gamma)

        REPLACE "\n" WITH ", " IN FileHeader$

        FileHeader$ = LEFT$(FileHeader$,LEN(FileHeader$)-2)

        FileHeader$ = "PROBE COORDINATES" + CHR$(13) +_
                      "-----------------" + CHR$(13) + FileHeader$

        N% = FREEFILE : OPEN "ProbeCoordinates.DAT" FOR OUTPUT AS #N%

            A$ = "   Step #     " : B$ = "  ------    " : C$ = ""

            FOR ProbeNum% = 1 TO Np% 'create out data file header

                SELECT CASE ProbeNum%
                    CASE   1 TO   9 : E$ = ""    : F$ = "           " : H$ = "           "
                    CASE  10 TO  99 : E$ = "-"   : F$ = "          " : H$ = "          "
                    CASE 100 TO 999 : E$ = "--"  : F$ = "         " : H$ = "         "
                END SELECT

                A$ = A$ + "P" + NoSpaces$(ProbeNum%+0,0) + F$ 'note: adding zero to ProbeNum% necessary to convert to floating point...

                B$ = B$ + E$ + "--" + H$

                C$ = C$ + "######.###    "
'
                 C$ = C$ + "##.#######"

            NEXT ProbeNum%

            PRINT #N%, FileHeader$ + CHR$(13) : PRINT #N%, A$ : PRINT #N%, B$

            FOR StepNum& = 0 TO LastStep&

                D$ = USING$("######   ",StepNum&)

                    FOR ProbeNum% = 1 TO Np% : D$ = D$ + USING$(C$,R(ProbeNum%,1,StepNum&)) : NEXT ProbeNum%

                PRINT #N%, D$

            NEXT StepNum&

        CLOSE #N%

END SUB 'Tabulate1dprobeCoordinates()

'------

SUB
Plot1DprobePositions(Max1DprobesPlotted%,Nd%,Np%,LastStep&,G,DeltaT,Alpha,Beta,Frep,R(),M(),PlaceInitialProbes$,InitialAcceleration$,RepositionFactor
$,FunctionName$,Gamma)
        'plots on-screen 1D function probe positions vs time step if Np =< 10

LOCAL ProcessID???, N%, n1%, n2%, n3%, n4%, n5%, n6%, n7%, n8%, n9%, n10%, n11%, n12%, n13%, n14%, n15%, ProbeNum%, StepNum&, A$

LOCAL PlotAnnotation$

    IF Np% > Max1DprobesPlotted% THEN EXIT SUB

    CALL CLEANUP 'delete old "Px" plot files, if any

    ProbeNum% = 0

    DO 'create output data files, probe-by-probe
        INCR ProbeNum% : n1%  = FREEFILE : OPEN "P"+REMOVE$(STR$(ProbeNum%),ANY" ") FOR OUTPUT AS #n1%  : IF ProbeNum% = Np% THEN EXIT SUB
        INCR ProbeNum% : n2%  = FREEFILE : OPEN "P"+REMOVE$(STR$(ProbeNum%),ANY" ") FOR OUTPUT AS #n2%  : IF ProbeNum% = Np% THEN EXIT LOOP
        INCR ProbeNum% : n3%  = FREEFILE : OPEN "P"+REMOVE$(STR$(ProbeNum%),ANY" ") FOR OUTPUT AS #n3%  : IF ProbeNum% = Np% THEN EXIT LOOP
        INCR ProbeNum% : n4%  = FREEFILE : OPEN "P"+REMOVE$(STR$(ProbeNum%),ANY" ") FOR OUTPUT AS #n4%  : IF ProbeNum% = Np% THEN EXIT LOOP
        INCR ProbeNum% : n5%  = FREEFILE : OPEN "P"+REMOVE$(STR$(ProbeNum%),ANY" ") FOR OUTPUT AS #n5%  : IF ProbeNum% = Np% THEN EXIT LOOP
        INCR ProbeNum% : n6%  = FREEFILE : OPEN "P"+REMOVE$(STR$(ProbeNum%),ANY" ") FOR OUTPUT AS #n6%  : IF ProbeNum% = Np% THEN EXIT LOOP
        INCR ProbeNum% : n7%  = FREEFILE : OPEN "P"+REMOVE$(STR$(ProbeNum%),ANY" ") FOR OUTPUT AS #n7%  : IF ProbeNum% = Np% THEN EXIT LOOP
        INCR ProbeNum% : n8%  = FREEFILE : OPEN "P"+REMOVE$(STR$(ProbeNum%),ANY" ") FOR OUTPUT AS #n8%  : IF ProbeNum% = Np% THEN EXIT LOOP
        INCR ProbeNum% : n9%  = FREEFILE : OPEN "P"+REMOVE$(STR$(ProbeNum%),ANY" ") FOR OUTPUT AS #n9%  : IF ProbeNum% = Np% THEN EXIT LOOP
        INCR ProbeNum% : n10% = FREEFILE : OPEN "P"+REMOVE$(STR$(ProbeNum%),ANY" ") FOR OUTPUT AS #n10% : IF ProbeNum% = Np% THEN EXIT LOOP
        INCR ProbeNum% : n11% = FREEFILE : OPEN "P"+REMOVE$(STR$(ProbeNum%),ANY" ") FOR OUTPUT AS #n11% : IF ProbeNum% = Np% THEN EXIT LOOP
        INCR ProbeNum% : n12% = FREEFILE : OPEN "P"+REMOVE$(STR$(ProbeNum%),ANY" ") FOR OUTPUT AS #n12% : IF ProbeNum% = Np% THEN EXIT LOOP
        INCR ProbeNum% : n13% = FREEFILE : OPEN "P"+REMOVE$(STR$(ProbeNum%),ANY" ") FOR OUTPUT AS #n13% : IF ProbeNum% = Np% THEN EXIT LOOP
        INCR ProbeNum% : n14% = FREEFILE : OPEN "P"+REMOVE$(STR$(ProbeNum%),ANY" ") FOR OUTPUT AS #n14  : IF ProbeNum% = Np% THEN EXIT LOOP
        INCR ProbeNum% : n15% = FREEFILE : OPEN "P"+REMOVE$(STR$(ProbeNum%),ANY" ") FOR OUTPUT AS #n15% : IF ProbeNum% = Np% THEN EXIT LOOP
    LOOP

    ProbeNum% = 0

    DO 'output probe positions as a function of time step
        INCR ProbeNum% : FOR StepNum& = 0 TO LastStep& : PRINT #n1%,  USING$("######  ######.########",StepNum&,R(ProbeNum%,1,StepNum&)) : NEXT
StepNum& : IF ProbeNum% = Np% THEN EXIT LOOP
        INCR ProbeNum% : FOR StepNum& = 0 TO LastStep& : PRINT #n2%,  USING$("######  ######.########",StepNum&,R(ProbeNum%,1,StepNum&)) : NEXT
StepNum& : IF ProbeNum% = Np% THEN EXIT LOOP
        INCR ProbeNum% : FOR StepNum& = 0 TO LastStep& : PRINT #n3%,  USING$("######  ######.########",StepNum&,R(ProbeNum%,1,StepNum&)) : NEXT
StepNum& : IF ProbeNum% = Np% THEN EXIT LOOP
        INCR ProbeNum% : FOR StepNum& = 0 TO LastStep& : PRINT #n4%,  USING$("######  ######.########",StepNum&,R(ProbeNum%,1,StepNum&)) : NEXT
StepNum& : IF ProbeNum% = Np% THEN EXIT LOOP
        INCR ProbeNum% : FOR StepNum& = 0 TO LastStep& : PRINT #n5%,  USING$("######  ######.########",StepNum&,R(ProbeNum%,1,StepNum&)) : NEXT
StepNum& : IF ProbeNum% = Np% THEN EXIT LOOP
        INCR ProbeNum% : FOR StepNum& = 0 TO LastStep& : PRINT #n6%,  USING$("######  ######.########",StepNum&,R(ProbeNum%,1,StepNum&)) : NEXT
StepNum& : IF ProbeNum% = Np% THEN EXIT LOOP
        INCR ProbeNum% : FOR StepNum& = 0 TO LastStep& : PRINT #n7%,  USING$("######  ######.########",StepNum&,R(ProbeNum%,1,StepNum&)) : NEXT
StepNum& : IF ProbeNum% = Np% THEN EXIT LOOP
        INCR ProbeNum% : FOR StepNum& = 0 TO LastStep& : PRINT #n8%,  USING$("######  ######.########",StepNum&,R(ProbeNum%,1,StepNum&)) : NEXT
StepNum& : IF ProbeNum% = Np% THEN EXIT LOOP
        INCR ProbeNum% : FOR StepNum& = 0 TO LastStep& : PRINT #n9%,  USING$("######  ######.########",StepNum&,R(ProbeNum%,1,StepNum&)) : NEXT
StepNum& : IF ProbeNum% = Np% THEN EXIT LOOP
        INCR ProbeNum% : FOR StepNum& = 0 TO LastStep& : PRINT #n10%, USING$("######  ######.########",StepNum&,R(ProbeNum%,1,StepNum&)) : NEXT
StepNum& : IF ProbeNum% = Np% THEN EXIT LOOP
        INCR ProbeNum% : FOR StepNum& = 0 TO LastStep& : PRINT #n11,  USING$("######  ######.########",StepNum&,R(ProbeNum%,1,StepNum&)) : NEXT
StepNum& : IF ProbeNum% = Np% THEN EXIT LOOP
        INCR ProbeNum% : FOR StepNum& = 0 TO LastStep& : PRINT #n12%, USING$("######  ######.########",StepNum&,R(ProbeNum%,1,StepNum&)) : NEXT
StepNum& : IF ProbeNum% = Np% THEN EXIT LOOP
        INCR ProbeNum% : FOR StepNum& = 0 TO LastStep& : PRINT #n13%, USING$("######  ######.########",StepNum&,R(ProbeNum%,1,StepNum&)) : NEXT
StepNum& : IF ProbeNum% = Np% THEN EXIT LOOP
        INCR ProbeNum% : FOR StepNum& = 0 TO LastStep& : PRINT #n14%, USING$("######  ######.########",StepNum&,R(ProbeNum%,1,StepNum&)) : NEXT
StepNum& : IF ProbeNum% = Np% THEN EXIT LOOP
        INCR ProbeNum% : FOR StepNum& = 0 TO LastStep& : PRINT #n15%, USING$("######  ######.########",StepNum&,R(ProbeNum%,1,StepNum&)) : NEXT
StepNum& : IF ProbeNum% = Np% THEN EXIT LOOP
```

```
    LOOP

    ProbeNum% = 0

    DO 'close output data files
        INCR ProbeNum% : CLOSE #n1% : IF ProbeNum% = Np% THEN EXIT LOOP
        INCR ProbeNum% : CLOSE #n2% : IF ProbeNum% = Np% THEN EXIT LOOP
        INCR ProbeNum% : CLOSE #n3% : IF ProbeNum% = Np% THEN EXIT LOOP
        INCR ProbeNum% : CLOSE #n4% : IF ProbeNum% = Np% THEN EXIT LOOP
        INCR ProbeNum% : CLOSE #n5% : IF ProbeNum% = Np% THEN EXIT LOOP
        INCR ProbeNum% : CLOSE #n6% : IF ProbeNum% = Np% THEN EXIT LOOP
        INCR ProbeNum% : CLOSE #n7% : IF ProbeNum% = Np% THEN EXIT LOOP
        INCR ProbeNum% : CLOSE #n8% : IF ProbeNum% = Np% THEN EXIT LOOP
        INCR ProbeNum% : CLOSE #n9% : IF ProbeNum% = Np% THEN EXIT LOOP
        INCR ProbeNum% : CLOSE #n10% : IF ProbeNum% = Np% THEN EXIT LOOP
        INCR ProbeNum% : CLOSE #n11% : IF ProbeNum% = Np% THEN EXIT LOOP
        INCR ProbeNum% : CLOSE #n12% : IF ProbeNum% = Np% THEN EXIT LOOP
        INCR ProbeNum% : CLOSE #n13% : IF ProbeNum% = Np% THEN EXIT LOOP
        INCR ProbeNum% : CLOSE #n14% : IF ProbeNum% = Np% THEN EXIT LOOP
        INCR ProbeNum% : CLOSE #n15% : IF ProbeNum% = Np% THEN EXIT LOOP

    LOOP

    ProbeNum% = 0 : A$ = ""

    DO 'create file string for plot command file
        INCR ProbeNum% : A$ = A$ + Quote$ + "P"+REMOVE$(STR$(ProbeNum%),ANY" ") + Quote$ + " w l lw 2, " : IF ProbeNum% = Np% THEN EXIT LOOP
        INCR ProbeNum% : A$ = A$ + Quote$ + "P"+REMOVE$(STR$(ProbeNum%),ANY" ") + Quote$ + " w l lw 2, " : IF ProbeNum% = Np% THEN EXIT LOOP
        INCR ProbeNum% : A$ = A$ + Quote$ + "P"+REMOVE$(STR$(ProbeNum%),ANY" ") + Quote$ + " w l lw 2, " : IF ProbeNum% = Np% THEN EXIT LOOP
        INCR ProbeNum% : A$ = A$ + Quote$ + "P"+REMOVE$(STR$(ProbeNum%),ANY" ") + Quote$ + " w l lw 2, " : IF ProbeNum% = Np% THEN EXIT LOOP
        INCR ProbeNum% : A$ = A$ + Quote$ + "P"+REMOVE$(STR$(ProbeNum%),ANY" ") + Quote$ + " w l lw 2, " : IF ProbeNum% = Np% THEN EXIT LOOP
        INCR ProbeNum% : A$ = A$ + Quote$ + "P"+REMOVE$(STR$(ProbeNum%),ANY" ") + Quote$ + " w l lw 2, " : IF ProbeNum% = Np% THEN EXIT LOOP
        INCR ProbeNum% : A$ = A$ + Quote$ + "P"+REMOVE$(STR$(ProbeNum%),ANY" ") + Quote$ + " w l lw 2, " : IF ProbeNum% = Np% THEN EXIT LOOP
        INCR ProbeNum% : A$ = A$ + Quote$ + "P"+REMOVE$(STR$(ProbeNum%),ANY" ") + Quote$ + " w l lw 2, " : IF ProbeNum% = Np% THEN EXIT LOOP
        INCR ProbeNum% : A$ = A$ + Quote$ + "P"+REMOVE$(STR$(ProbeNum%),ANY" ") + Quote$ + " w l lw 2, " : IF ProbeNum% = Np% THEN EXIT LOOP
        INCR ProbeNum% : A$ = A$ + Quote$ + "P"+REMOVE$(STR$(ProbeNum%),ANY" ") + Quote$ + " w l lw 2, " : IF ProbeNum% = Np% THEN EXIT LOOP
        INCR ProbeNum% : A$ = A$ + Quote$ + "P"+REMOVE$(STR$(ProbeNum%),ANY" ") + Quote$ + " w l lw 2, " : IF ProbeNum% = Np% THEN EXIT LOOP
        INCR ProbeNum% : A$ = A$ + Quote$ + "P"+REMOVE$(STR$(ProbeNum%),ANY" ") + Quote$ + " w l lw 2, " : IF ProbeNum% = Np% THEN EXIT LOOP
        INCR ProbeNum% : A$ = A$ + Quote$ + "P"+REMOVE$(STR$(ProbeNum%),ANY" ") + Quote$ + " w l lw 2, " : IF ProbeNum% = Np% THEN EXIT LOOP
        INCR ProbeNum% : A$ = A$ + Quote$ + "P"+REMOVE$(STR$(ProbeNum%),ANY" ") + Quote$ + " w l lw 2, " : IF ProbeNum% = Np% THEN EXIT LOOP
        INCR ProbeNum% : A$ = A$ + Quote$ + "P"+REMOVE$(STR$(ProbeNum%),ANY" ") + Quote$ + " w l lw 2, " : IF ProbeNum% = Np% THEN EXIT LOOP

    LOOP

    A$ = LEFT$(A$,LEN(A$)-2)

    CALL
GetPlotAnnotation(PlotAnnotation$,Nd%,Np%,LastStep&,G,DeltaT,Alpha,Beta,Frep,M(),PlaceInitialProbes$,InitialAcceleration$,RepositionFactor$,FunctionN
ame$,Gamma)

    N% = FREEFILE

    OPEN "cmd2d.gp" FOR OUTPUT AS #N%

        PRINT #N%, "set label "       + Quote$ + PlotAnnotation$ + Quote$ + " at graph 0.5,0.95"
        PRINT #N%, "set grid xtics"
        PRINT #N%, "set grid ytics"
        PRINT #N%, "set title "  + Quote$ + "Evolution of "     + FunctionName$ + " Probe Positions"+ "\n" + RunID$ + Quote$
        PRINT #N%, "set xlabel " + Quote$ + "Time Step"      + Quote$
        PRINT #N%, "set ylabel " + Quote$ + "Probe Coordinate" + Quote$
        PRINT #N%, "plot "       + A$

    CLOSE #N%

    CALL CreateGNUplotINIfile(0.2##*Screenwidth&,0.2##*ScreenHeight&,0.6##*Screenwidth&,0.6##*ScreenHeight&) 'USAGE: CALL
CreateGNUplotINIfile(PlotwindowULC_X%,PlotwindowULC_Y%,PlotwindowWidth%,PlotwindowHeight%)

    ProcessID??? = SHELL("wgnuplot.exe"+" cmd2d.gp -") : CALL Delay(5##) 'before SUB Cleanup is called

END SUB

'------

SUB
DisplayRunParameters(FunctionName$,Nd%,Np%,Nt&,G,DeltaT,Alpha,Beta,Frep,R(),A(),M(),PlaceInitialProbes$,InitialAcceleration$,RepositionFactor$,RunCFO
$,ShrinkDS$,CheckForEarlyTermination$)

LOCAL A$, B$, YN&

    B$ = "" : IF PlaceInitialProbes$ = "UNIFORM ON-AXIS" AND Nd% > 1 THEN B$ = "  ["+REMOVE$(STR$(Np%/Nd%),ANY" ") + "/axis]"

    RunCFO$ = "NO"

A$ = "RUN CFO WITH THE" + CHR$(13) +_
        "FOLLOWING PARAMETERS?"                           + CHR$(13) + CHR$(13) +_
        "Function "       + FunctionName$               + " (" + REMOVE$(STR$(Nd%),ANY" ") + "-D)" + CHR$(13) +_
        "# time steps = " + REMOVE$(STR$(Nt&),ANY" ")    + CHR$(13) + _
        "Grav Const G = " + REMOVE$(STR$(G,2),ANY" ")    + CHR$(13) + _
        "Delta-T = "      + REMOVE$(STR$(DeltaT,3),ANY" ") + CHR$(13) + _
        "Exp Alpha = "    + REMOVE$(STR$(Alpha,3),ANY" ") + CHR$(13) + _
        "Exp Beta = "     + REMOVE$(STR$(Beta,3),ANY" ")  + CHR$(13) + _
        "Frep = "         + REMOVE$(STR$(Frep,4),ANY" ")  + "  ("+RepositionFactor$+ ")" + CHR$(13) + _
        "Initial Probes: " + PlaceInitialProbes$         + B$ + CHR$(13) + _
        "Initial Accel: " + InitialAcceleration$         + CHR$(13) + _
        "Check for Early Termination? " + CheckForEarlyTermination$ + CHR$(13) + _
        "Shrink Decision Space? "        + ShrinkDS$ + CHR$(13) +CHR$(13)

'   lResult& = MSGBOX(txt$ [, [style&], title$])

    A$ = "RUN CFO ON FUNCTION " + FunctionName$ + "?"

    YN& = MSGBOX(A$,%MB_YESNO,"CONFIRM RUN")

    IF YN& = %IDYES THEN RunCFO$ = "YES"

END SUB

'------

SUB StatusWindow(FunctionName$,StatusWindowHandle???)

    GRAPHIC WINDOW "Run Progress, "+FunctionName$,0.08##*Screenwidth&,0.08##*ScreenHeight&,0.25##*Screenwidth&,0.17##*ScreenHeight& TO
StatusWindowHandle???

    GRAPHIC ATTACH StatusWindowHandle???,0,REDRAW

    GRAPHIC FONT "Lucida Console",8,0 '"Courier New",8,0 'Fixed width fonts

    GRAPHIC SET PIXEL (35,15) : GRAPHIC PRINT "  Initializing...      " : GRAPHIC REDRAW
```

Page 73 of 82

```
END SUB

'------

SUB GetTestFunctionNumber(FunctionName$)

   LOCAL hDlg AS DWORD

   LOCAL N%, M%

   LOCAL FrameWidth&, FrameHeight&, BoxWidth&, BoxHeight&

' BoxWidth& = 276 : BoxHeight& = 300 : FrameWidth& = 82 : FrameHeight = BoxHeight&-5

   BoxWidth& = 276 : BoxHeight& = 300 : FrameWidth& = 90 : FrameHeight& = BoxHeight&-5

   DIALOG NEW 0, "CENTRAL FORCE OPTIMIZATION TEST FUNCTIONS",,, BoxWidth&, BoxHeight&, %WS_CAPTION OR %WS_SYSMENU, 0 TO hDlg

'----------------------------------------------------------------

   CONTROL ADD FRAME, hDlg, %IDC_FRAME1, "Test Functions",     5,  2, FrameWidth&, FrameHeight&
   CONTROL ADD FRAME, hDlg, %IDC_FRAME2, "GSO Test Functions", 105, 2, FrameWidth&, 255

   CONTROL ADD OPTION, hDlg, %IDC_Function_Number1, "Parrott F4",   10,  14, 70, 10, %WS_GROUP OR %WS_TABSTOP
   CONTROL ADD OPTION, hDlg, %IDC_Function_Number2, "SGO",          10,  24, 70, 10
   CONTROL ADD OPTION, hDlg, %IDC_Function_Number3, "Goldstein-Price",10, 34, 70, 10
   CONTROL ADD OPTION, hDlg, %IDC_Function_Number4, "Step",         10,  44, 70, 10
   CONTROL ADD OPTION, hDlg, %IDC_Function_Number5, "Schwefel 2.26",10,  54, 70, 10
   CONTROL ADD OPTION, hDlg, %IDC_Function_Number6, "Colville",     10,  64, 70, 10
   CONTROL ADD OPTION, hDlg, %IDC_Function_Number7, "Griewank",     10,  74, 70, 10

   CONTROL ADD OPTION, hDlg, %IDC_Function_Number31, "PBM #1",       10,  84, 70, 10
   CONTROL ADD OPTION, hDlg, %IDC_Function_Number32, "PBM #2",       10,  94, 70, 10
   CONTROL ADD OPTION, hDlg, %IDC_Function_Number33, "PBM #3",       10, 104, 70, 10
   CONTROL ADD OPTION, hDlg, %IDC_Function_Number34, "PBM #4",       10, 114, 70, 10
   CONTROL ADD OPTION, hDlg, %IDC_Function_Number35, "PBM #5",       10, 124, 70, 10
   CONTROL ADD OPTION, hDlg, %IDC_Function_Number36, "Himmelblau",   10, 134, 70, 10
   CONTROL ADD OPTION, hDlg, %IDC_Function_Number37, "Rosenbrock",   10, 144, 70, 10
   CONTROL ADD OPTION, hDlg, %IDC_Function_Number38, "Sphere",       10, 154, 70, 10
   CONTROL ADD OPTION, hDlg, %IDC_Function_Number39, "HimmelblauNLO",10, 164, 70, 10
   CONTROL ADD OPTION, hDlg, %IDC_Function_Number40, "Tripod",       10, 174, 70, 10
   CONTROL ADD OPTION, hDlg, %IDC_Function_Number41, "Rosenbrock F6",10, 184, 70, 10
   CONTROL ADD OPTION, hDlg, %IDC_Function_Number42, "Comp Spring",  10, 194, 70, 10
   CONTROL ADD OPTION, hDlg, %IDC_Function_Number43, "Gear Train",   10, 204, 70, 10
   CONTROL ADD OPTION, hDlg, %IDC_Function_Number44, "Loaded Monopole",10, 214, 70, 10
   CONTROL ADD OPTION, hDlg, %IDC_Function_Number45, "Reserved",     10, 224, 70, 10
   CONTROL ADD OPTION, hDlg, %IDC_Function_Number46, "Reserved",     10, 234, 70, 10
   CONTROL ADD OPTION, hDlg, %IDC_Function_Number47, "Reserved",     10, 244, 70, 10
   CONTROL ADD OPTION, hDlg, %IDC_Function_Number48, "Reserved",     10, 254, 70, 10
   CONTROL ADD OPTION, hDlg, %IDC_Function_Number49, "Reserved",     10, 264, 70, 10
   CONTROL ADD OPTION, hDlg, %IDC_Function_Number50, "Reserved",     10, 274, 70, 10

' --------------------- Test Functions from GSO Paper ---------------------
   CONTROL ADD OPTION, hDlg, %IDC_Function_Number8,  "F1" , 120,  14, 40, 10
   CONTROL ADD OPTION, hDlg, %IDC_Function_Number9,  "F2" , 120,  24, 40, 10
   CONTROL ADD OPTION, hDlg, %IDC_Function_Number10, "F3" , 120,  34, 40, 10
   CONTROL ADD OPTION, hDlg, %IDC_Function_Number11, "F4" , 120,  44, 40, 10
   CONTROL ADD OPTION, hDlg, %IDC_Function_Number12, "F5" , 120,  54, 40, 10
   CONTROL ADD OPTION, hDlg, %IDC_Function_Number13, "F6" , 120,  64, 40, 10
   CONTROL ADD OPTION, hDlg, %IDC_Function_Number14, "F7" , 120,  74, 40, 10
   CONTROL ADD OPTION, hDlg, %IDC_Function_Number15, "F8" , 120,  84, 40, 10
   CONTROL ADD OPTION, hDlg, %IDC_Function_Number16, "F9" , 120,  94, 40, 10
   CONTROL ADD OPTION, hDlg, %IDC_Function_Number17, "F10", 120, 104, 40, 10
   CONTROL ADD OPTION, hDlg, %IDC_Function_Number18, "F11", 120, 114, 40, 10
   CONTROL ADD OPTION, hDlg, %IDC_Function_Number19, "F12", 120, 124, 40, 10
   CONTROL ADD OPTION, hDlg, %IDC_Function_Number20, "F13", 120, 134, 40, 10
   CONTROL ADD OPTION, hDlg, %IDC_Function_Number21, "F14", 120, 144, 40, 10
   CONTROL ADD OPTION, hDlg, %IDC_Function_Number22, "F15", 120, 154, 40, 10
   CONTROL ADD OPTION, hDlg, %IDC_Function_Number23, "F16", 120, 164, 40, 10
   CONTROL ADD OPTION, hDlg, %IDC_Function_Number24, "F17", 120, 174, 40, 10
   CONTROL ADD OPTION, hDlg, %IDC_Function_Number25, "F18", 120, 184, 40, 10
   CONTROL ADD OPTION, hDlg, %IDC_Function_Number26, "F19", 120, 194, 40, 10
   CONTROL ADD OPTION, hDlg, %IDC_Function_Number27, "F20", 120, 204, 40, 10
   CONTROL ADD OPTION, hDlg, %IDC_Function_Number28, "F21", 120, 214, 40, 10
   CONTROL ADD OPTION, hDlg, %IDC_Function_Number29, "F22", 120, 224, 40, 10
   CONTROL ADD OPTION, hDlg, %IDC_Function_Number30, "F23", 120, 234, 40, 10

   CONTROL SET OPTION  hDlg, %IDC_Function_Number1, %IDC_Function_Number1, %IDC_Function_Number3 'default to Parrott F4

'----------------------------------------------------------------

   CONTROL ADD BUTTON, hDlg, %IDOK, "&OK", 200, 0.45##*BoxHeight&, 50, 14

'----------------------------------------------------------------

   DIALOG SHOW MODAL hDlg CALL DlgProc

   CALL Delay(1##)

   IF FunctionNumber% < 1 OR FunctionNumber% > 44 THEN

      FunctionNumber% = 1 : MSGBOX("Error in function number...")

   END IF

   SELECT CASE FunctionNumber%

      CASE 1 : FunctionName$ = "ParrottF4"
      CASE 2 : FunctionName$ = "SGO"
      CASE 3 : FunctionName$ = "GP"
      CASE 4 : FunctionName$ = "STEP"
      CASE 5 : FunctionName$ = "SCHWEFEL_226"
      CASE 6 : FunctionName$ = "COLVILLE"
      CASE 7 : FunctionName$ = "GRIEWANK"
      CASE 8 : FunctionName$ = "F1"
      CASE 9 : FunctionName$ = "F2"
      CASE 10 : FunctionName$ = "F3"
      CASE 11 : FunctionName$ = "F4"
      CASE 12 : FunctionName$ = "F5"
      CASE 13 : FunctionName$ = "F6"
      CASE 14 : FunctionName$ = "F7"
      CASE 15 : FunctionName$ = "F8"
      CASE 16 : FunctionName$ = "F9"
      CASE 17 : FunctionName$ = "F10"
      CASE 18 : FunctionName$ = "F11"
      CASE 19 : FunctionName$ = "F12"
```

```
            CASE 20: FunctionName$ = "F13"
            CASE 21: FunctionName$ = "F14"
            CASE 22: FunctionName$ = "F15"
            CASE 23: FunctionName$ = "F16"
            CASE 24: FunctionName$ = "F17"
            CASE 25: FunctionName$ = "F18"
            CASE 26: FunctionName$ = "F19"
            CASE 27: FunctionName$ = "F20"
            CASE 28: FunctionName$ = "F21"
            CASE 29: FunctionName$ = "F22"
            CASE 30: FunctionName$ = "F23"
            CASE 31: FunctionName$ = "PBM_1"
            CASE 32: FunctionName$ = "PBM_2"
            CASE 33: FunctionName$ = "PBM_3"
            CASE 34: FunctionName$ = "PBM_4"
            CASE 35: FunctionName$ = "PBM_5"
            CASE 36: FunctionName$ = "HIMMELBLAU"
            CASE 37: FunctionName$ = "ROSENBROCK"
            CASE 38: FunctionName$ = "SPHERE"
            CASE 39: FunctionName$ = "HIMMELBLAUNLO"
            CASE 40: FunctionName$ = "TRIPOD"
            CASE 41: FunctionName$ = "ROSENBROCKF6"
            CASE 42: FunctionName$ = "COMPRESSIONSPRING"
            CASE 43: FunctionName$ = "GEARTRAIN"
            CASE 44: FunctionName$ = "LD_MONO"

    END SELECT

END SUB

'----------

CALLBACK FUNCTION DlgProc() AS LONG

    '--------------------------------------------------------------
    ' Callback procedure for the main dialog
    '--------------------------------------------------------------
    LOCAL c, lRes AS LONG, sText AS STRING

    SELECT CASE AS LONG CBMSG

    CASE %WM_INITDIALOG' %WM_INITDIALOG is sent right before the dialog is shown.

    CASE %WM_COMMAND            ' <- a control is calling

        SELECT CASE AS LONG CBCTL  ' <- look at control's id

        CASE %IDOK                ' <- OK button or Enter key was pressed

            IF CBCTLMSG = %BN_CLICKED THEN
                '---------------------------------------
                ' Loop through the Function_Number controls
                ' to see which one is selected
                '---------------------------------------
                FOR c = %IDC_Function_Number1 TO %IDC_Function_Number50

                    CONTROL GET CHECK CBHNDL, c TO lRes

                    IF lRes THEN EXIT FOR

                NEXT 'c holds the id for selected test function.

                FunctionNumber% = c-120

                DIALOG END CBHNDL

            END IF
        END SELECT
    END SELECT

END FUNCTION

'------------------------------ PBM ANTENNA BENCHMARK FUNCTIONS ------------------------------

'Reference for benchmarks PBM_1 through PBM_5:

'Pantoja, M F., Bretones, A. R., Martin, R. G., "Benchmark Antenna Problems for Evolutionary
'Optimization Algorithms," IEEE Trans. Antennas & Propagation, vol. 55, no. 4, April 2007,
'pp. 1111-1121

FUNCTION PBM_1(R(),Nd%,p%,j&) 'PBM Benchmark #1: Max D for Variable-Length CF Dipole

    LOCAL Z, LengthWaves, ThetaRadians AS EXT

    LOCAL N%, Nsegs%, FeedSegNum%

    LOCAL NumSegs$, FeedSeg$, HalfLength$, Radius$, ThetaDeg$, Lyne$, GainDB$

    LengthWaves  = R(p%,1,j&)

    ThetaRadians = R(p%,2,j&)

    ThetaDeg$ = REMOVE$(STR$(ROUND(ThetaRadians*Rad2Deg,2)),ANY" ")

    IF TALLY(ThetaDeg$,".") = 0 THEN ThetaDeg$ = ThetaDeg$+"."

    Nsegs% = 2*(INT(100*LengthWaves)\2)+1 '100 segs per wavelength, must be an odd #, VOLTAGE SOURCE

    FeedSegNum% = Nsegs%\2 + 1 'center segment number, VOLTAGE SOURCE

    NumSegs$    = REMOVE$(STR$(Nsegs%),ANY" ")

    FeedSeg$    = REMOVE$(STR$(FeedSegNum%),ANY" ")

    HalfLength$ = REMOVE$(STR$(ROUND(LengthWaves/2##,6)),ANY" ")

    IF TALLY(HalfLength$,".") = 0 THEN HalfLength$ = HalfLength$+"."

    Radius$    = "0.00001" 'REMOVE$(STR$(ROUND(LengthWaves/1000##,6)),ANY" ")

    N% = FREEFILE

    OPEN "PBM1.NEC" FOR OUTPUT AS #N%

        PRINT #N%,"CM File: PBM1.NEC"
```

```
        PRINT #N%,"CM Run ID "+DATE$+" "+TIME$
        PRINT #N%,"CM Nd="+STR$(Nd%)+", p="STR$(p%)+", j="STR$(j&)
        PRINT #N%,"CM R(p,1,j)="+STR$(R(p%,1,j&))+", R(p,2,j)="+STR$(R(p%,2,j&))
        PRINT #N%,"CE"
        PRINT #N%,"GW 1,"+NumSegs$+",0.,0.,-"+HalfLength$+",0.,0.,"+HalfLength$+","+Radius$
        PRINT #N%,"GE"
        PRINT #N%,"EX 0,1,"+FeedSeg$+",0,1.,0." 'VOLTAGE SOURCE
        PRINT #N%,"FR 0,1,0,0,299.79564,0."
        PRINT #N%,"RP 0,1,1,1001,"+ThetaDeg$+",0.,0.,0.,1000.' 'gain at 1000 wavelengths range
        PRINT #N%,"XQ"
        PRINT #N%,"EN"

    CLOSE #N%

' - - - ANGLES - -        - - POWER GAINS -       - - - POLARIZATION - - -     - - E(THETA) - - -   - - - E(PHI) - - -
'  THETA     PHI      VERT.   HOR.    TOTAL     AXIAL    TILT  SENSE     MAGNITUDE    PHASE       MAGNITUDE      PHASE
' DEGREES  DEGREES     DB      DB       DB       RATIO    DEG.           VOLTS/M     DEGREES       VOLTS/M      DEGREES
'  90.00     0.00     3.91  -999.99    3.91    0.00000    0.00  LINEAR   1.29504E-04    5.37     0.00000E+00     -5.24
'123456789x123456789x123456789x123456789x123456789x123456789x123456789x123456789x123456789x123456789x123456789x123456789
'       10        20        30        40        50        60        70        80        90       100       110       120

    SHELL "n41_2k1.exe",0

    N% = FREEFILE

    OPEN "PBM1.OUT" FOR INPUT AS #N%

        WHILE NOT EOF(N%)

            LINE INPUT #N%, Lyne$

            IF INSTR(Lyne$,"DEGREES  DEGREES") > 0 THEN EXIT LOOP

        WEND 'position at next data line

        LINE INPUT #N%, Lyne$

    CLOSE #N%

    GainDB$ = REMOVE$(MID$(Lyne$,37,8),ANY" ")

    PBM_1 = 10^(VAL(GainDB$)/10##) 'Directivity

END FUNCTION 'PBM_1()

'----

FUNCTION PBM_2(R(),Nd%,p%,j&) 'PBM Benchmark #2: Max D for Variable-Separation Array of CF Dipoles

    LOCAL Z, DipoleSeparationWaves, ThetaRadians AS EXT

    LOCAL N%, i%

    LOCAL NumSegs$, FeedSeg$, Radius$, ThetaDeg$, Lyne$, GainDB$, Xcoord$, WireNum$

    DipoleSeparationWaves = R(p%,1,j&)

    ThetaRadians          = R(p%,2,j&)

    ThetaDeg$ = REMOVE$(STR$(ROUND(ThetaRadians*Rad2Deg,2)),ANY" ")

    IF TALLY(ThetaDeg$,".") = 0 THEN ThetaDeg$ = ThetaDeg$+"."

    NumSegs$ = "49"

    FeedSeg$ = "25"

    Radius$  = "0.00001"

    N% = FREEFILE

    OPEN "PBM2.NEC" FOR OUTPUT AS #N%

        PRINT #N%,"CM File: PBM2.NEC"
        PRINT #N%,"CM Run ID "+DATE$+" "+TIME$
        PRINT #N%,"CM Nd="+STR$(Nd%)+", p="STR$(p%)+", j="STR$(j&)
        PRINT #N%,"CM R(p,1,j)="+STR$(R(p%,1,j&))+", R(p,2,j)="+STR$(R(p%,2,j&))
        PRINT #N%,"CE"

        FOR i% = -9 TO 9 STEP 2
            WireNum$ = REMOVE$(STR$((i%+11)\2),ANY" ")
            Xcoord$  = REMOVE$(STR$(i%*DipoleSeparationWaves/2##),ANY" ")
            PRINT #N%,"GW "+WireNum$+","+NumSegs$+","+Xcoord$+",0.,-0.25,"+Xcoord$+",0.,0.25,"+Radius$
        NEXT i%

        PRINT #N%,"GE"

        FOR i% = 1 TO 10
            PRINT #N%,"EX 0,"+REMOVE$(STR$(i%),ANY" ")+",+FeedSeg$+",0,1.,0." 'VOLTAGE SOURCE
        NEXT i%
        PRINT #N%,"FR 0,1,0,0,299.79564,0."
        PRINT #N%,"RP 0,1,1,1001,"+ThetaDeg$+",90.,0.,0.,1000." 'gain at 1000 wavelengths range
        PRINT #N%,"XQ"
        PRINT #N%,"EN"

    CLOSE #N%

' - - - ANGLES - -        - - POWER GAINS -       - - - POLARIZATION - - -     - - E(THETA) - - -   - - - E(PHI) - - -
'  THETA     PHI      VERT.   HOR.    TOTAL     AXIAL    TILT  SENSE     MAGNITUDE    PHASE       MAGNITUDE      PHASE
' DEGREES  DEGREES     DB      DB       DB       RATIO    DEG.           VOLTS/M     DEGREES       VOLTS/M      DEGREES
'  90.00     0.00     3.91  -999.99    3.91    0.00000    0.00  LINEAR   1.29504E-04    5.37     0.00000E+00     -5.24
'123456789x123456789x123456789x123456789x123456789x123456789x123456789x123456789x123456789x123456789x123456789x123456789
'       10        20        30        40        50        60        70        80        90       100       110       120

    SHELL "n41_2k1.exe",0

    N% = FREEFILE

    OPEN "PBM2.OUT" FOR INPUT AS #N%

        WHILE NOT EOF(N%)

            LINE INPUT #N%, Lyne$

            IF INSTR(Lyne$,"DEGREES  DEGREES") > 0 THEN EXIT LOOP

        WEND 'position at next data line
```

```
        LINE INPUT #N%, Lyne$

    CLOSE #N%

    GainDB$ = REMOVE$(MID$(Lyne$,37,8),ANY" ")

    IF AddNoiseToPBM2$ = "YES" THEN

        Z = 10^(VAL(GainDB$)/10##) + GaussianDeviate(0##,0.4472##) 'Directivity with Gaussian noise (zero mean, 0.2 variance)

    ELSE

        Z = 10^(VAL(GainDB$)/10##) 'Directivity without noise

    END IF

    PBM_2 = Z

END FUNCTION 'PBM_2()

'----

FUNCTION PBM_3(R(),Nd%,p%,j&) 'PBM Benchmark #3: Max D for Circular Dipole Array

    LOCAL Beta, ThetaRadians, Alpha, ReV, ImV AS EXT

    LOCAL N%, i%

    LOCAL NumSegs$, FeedSeg$, Radius$, ThetaDeg$, Lyne$, GainDB$, Xcoord$, Ycoord$, WireNum$, ReEX$, ImEX$

    Beta         = R(p%,1,j&)

    ThetaRadians = R(p%,2,j&)

    ThetaDeg$ = REMOVE$(STR$(ROUND(ThetaRadians*Rad2Deg,2)),ANY" ")

    IF TALLY(ThetaDeg$,".") = 0 THEN ThetaDeg$ = ThetaDeg$+"."

    NumSegs$ = "49"

    FeedSeg$ = "25"

    Radius$  = "0.00001"

    N% = FREEFILE

    OPEN "PBM3.NEC" FOR OUTPUT AS #N%
        PRINT #N%,"CM File: PBM3.NEC"
        PRINT #N%,"CM Run ID "+DATE$+" "+TIME$
        PRINT #N%,"CM Nd="+STR$(Nd%)+", p="+STR$(p%)+", j="+STR$(j&)
        PRINT #N%,"CM R(p,1,j)="+STR$(R(p%,1,j&))+", R(p,2,j)="+STR$(R(p%,2,j&))
        PRINT #N%,"CE"

        FOR i% = 1 TO 8
            WireNum$ = REMOVE$(STR$(i%),ANY" ")

            SELECT CASE i%
                CASE 1 : Xcoord$ = "1"        : Ycoord$ = "0"
                CASE 2 : Xcoord$ = "0.70711"  : Ycoord$ = "0.70711"
                CASE 3 : Xcoord$ = "0"        : Ycoord$ = "1"
                CASE 4 : Xcoord$ = "-0.70711" : Ycoord$ = "0.70711"
                CASE 5 : Xcoord$ = "-1"       : Ycoord$ = "0"
                CASE 6 : Xcoord$ = "-0.70711" : Ycoord$ = "-0.70711"
                CASE 7 : Xcoord$ = "0"        : Ycoord$ = "-1"
                CASE 8 : Xcoord$ = "0.70711"  : Ycoord$ = "-0.70711"
            END SELECT

            PRINT #N%,"GW "+WireNum$+","+NumSegs$+","+Xcoord$+","+Ycoord$+",-0.25,"+Xcoord$+","+Ycoord$+",0.25,"+Radius$
        NEXT i%

        PRINT #N%,"GE"

        FOR i% = 1 TO 8
            Alpha = -COS(TwoPi*Beta*(i%-1))

            ReV   = COS(Alpha)
            ImV   = SIN(Alpha)

            ReEX$ = REMOVE$(STR$(ROUND(ReV,6)),ANY" ")
            ImEX$ = REMOVE$(STR$(ROUND(ImV,6)),ANY" ")

            IF TALLY(ReEX$,".") = 0 THEN ReEX$ = ReEX$+"."
            IF TALLY(ImEX$,".") = 0 THEN ImEX$ = ImEX$+"."

            PRINT #N%,"EX 0,"+REMOVE$(STR$(i%),ANY" ")+","+FeedSeg$+",0,"+ReEX$+","+ImEX$ 'VOLTAGE SOURCE
        NEXT i%

        PRINT #N%,"FR 0,1,0,0,299.79564,0."
        PRINT #N%,"RP 0,1,1,1001,"+ThetaDeg$+",0.,0.,0.,1000." 'gain at 1000 wavelengths range
        PRINT #N%,"XQ"
        PRINT #N%,"EN"

    CLOSE #N%

' - -  - ANGLES - -        - POWER GAINS -        AXIAL      - - -POLARIZATION - - -   - - E(THETA) - - -   - - - E(PHI) - - -
'  THETA      PHI     VERT.      HOR.     TOTAL      AXIAL    TILT   SENSE      MAGNITUDE   PHASE      MAGNITUDE    PHASE
' DEGREES  DEGREES      DB        DB        DB        RATIO    DEG.              VOLTS/M   DEGREES      VOLTS/M   DEGREES
'  90.00     0.00     3.91   -999.99     3.91     0.00000     0.00  LINEAR    1.29504E-04    5.37    0.00000E+00   -5.24
'123456789x123456789x123456789x123456789x123456789x123456789x123456789x123456789x123456789x123456789x123456789x123456789x
'     10        20        30        40        50        60        70        80        90       100       110       120

    SHELL "n41_2k1.exe",0

    N% = FREEFILE

    OPEN "PBM3.OUT" FOR INPUT AS #N%

        WHILE NOT EOF(N%)

            LINE INPUT #N%, Lyne$

            IF INSTR(Lyne$,"DEGREES  DEGREES") > 0 THEN EXIT LOOP

        WEND 'position at next data line
```

```
        LINE INPUT #N%, Lyne$

    CLOSE #N%

    GainDB$ = REMOVE$(MID$(Lyne$,37,8),ANY" ")
    PBM_3 = 10^(VAL(GainDB$)/10##) 'Directivity
END FUNCTION 'PBM_3()
'----
FUNCTION PBM_4(R(),Nd%,p%,j&) 'PBM Benchmark #4: Max D for Vee Dipole
    LOCAL TotalLengthWaves, AlphaRadians, ArmLength, Xlength, Zlength, Lfeed AS EXT

    LOCAL N%, i%, Nsegs%, Feedzcoord$

    LOCAL NumSegs$, Lyne$, GainDB$, Xcoord$, Zcoord$
    TotalLengthWaves = 2##*R(p%,1,j&)

    AlphaRadians     = R(p%,2,j&)

    Lfeed            = 0.01##

    Feedzcoord$      = REMOVE$(STR$(Lfeed),ANY" ")

    ArmLength = (TotalLengthWaves-2##*Lfeed)/2##

    Xlength  = ROUND(ArmLength*COS(AlphaRadians),6)

    Xcoord$  = REMOVE$(STR$(Xlength),ANY" ") : IF TALLY(Xcoord$,".") = 0 THEN Xcoord$ = Xcoord$+"."

    Zlength  = ROUND(ArmLength*SIN(AlphaRadians),6)

    Zcoord$  = REMOVE$(STR$(Zlength+Lfeed),ANY" ") : IF TALLY(Zcoord$,".") = 0 THEN Zcoord$ = Zcoord$+"."

    Nsegs%   = 2*(INT(TotalLengthWaves*100)\2) 'even number, total # segs

    NumSegs$ = REMOVE$(STR$(Nsegs%\2),ANY" ") '# segs per arm

    N% = FREEFILE

    OPEN "PBM4.NEC" FOR OUTPUT AS #N%
        PRINT #N%,"CM File: PBM4.NEC"
        PRINT #N%,"CM Run ID "+DATE$+" "+TIME$
        PRINT #N%,"CM Nd="+STR$(Nd%)+", p="+STR$(p%)+", j="+STR$(j&)
        PRINT #N%,"CM R(p,1,j)="+STR$(R(p%,1,j&))+",  R(p,2,j)="+STR$(R(p%,2,j&))
        PRINT #N%,"CE"

        PRINT #N%,"GW 1,5,0.,0.,-"+Feedzcoord$+",0.,0.,"+Feedzcoord$+",0.00001" 'feed wire, 1 segment, 0.01 wvln

        PRINT #N%,"GW 2,"+NumSegs$+",0.,0.,"+Feedzcoord$+","+Xcoord$+",0.,"+Zcoord$+",0.00001" 'upper arm

        PRINT #N%,"GW 3,"+NumSegs$+",0.,0.,-"+Feedzcoord$+","+Xcoord$+",0.,-"+Zcoord$+",0.00001" 'lower arm

        PRINT #N%,"GE"

        PRINT #N%,"EX 0,1,3,0,1.,0." 'VOLTAGE SOURCE

        PRINT #N%,"FR 0,1,0,0,299.79564,0."
        PRINT #N%,"RP 0,1,1,1001,90.,0.,0.,0.,1000." 'ENDFIRE gain at 1000 wavelengths range
        PRINT #N%,"XQ"
        PRINT #N%,"EN"

    CLOSE #N%
' - - ANGLES - -        - POWER GAINS -       - - -POLARIZATION - - -     - - E(THETA) - - -    - - - E(PHI) - - -
'   THETA      PHI       VERT.   HOR.   TOTAL    AXIAL   TILT  SENSE   MAGNITUDE    PHASE    MAGNITUDE    PHASE
' DEGREES   DEGREES       DB      DB      DB      RATIO   DEG.          VOLTS/M   DEGREES     VOLTS/M    DEGREES
'   90.00      0.00      3.91  -999.99   3.91   0.00000   0.00 LINEAR  1.29504E-04   5.37   0.00000E+00   -5.24
'123456789x123456789x123456789x123456789x123456789x123456789x123456789x123456789x123456789x123456789x123456789x
'     10        20        30        40        50        60        70        80        90       100       110       120

    SHELL "n41_2k1.exe",0

    N% = FREEFILE

    OPEN "PBM4.OUT" FOR INPUT AS #N%

        WHILE NOT EOF(N%)

            LINE INPUT #N%, Lyne$

            IF INSTR(Lyne$,"DEGREES  DEGREES") > 0 THEN EXIT LOOP

        WEND 'position at next data line

        LINE INPUT #N%, Lyne$

    CLOSE #N%

    GainDB$ = REMOVE$(MID$(Lyne$,37,8),ANY" ")
    PBM_4 = 10^(VAL(GainDB$)/10##) 'Directivity
END FUNCTION 'PBM_4()
'----
FUNCTION PBM_5(R(),Nd%,p%,j&) 'PBM Benchmark #5: N-element collinear array (Nd=N-1)
    LOCAL TotalLengthWaves, Di(), Ystart, Y1, Y2, SumDi AS EXT

    LOCAL N%, i%, q%

    LOCAL Lyne$, GainDB$
    REDIM Di(1 TO Nd%)

    FOR i% = 1 TO Nd%

        Di(i%) = R(p%,i%,j&) 'dipole separation, wavelengths

    NEXT i%
```

```
    TotalLengthWaves = 0##

    FOR i% = 1 TO Nd%

        TotalLengthWaves = TotalLengthWaves + Di(i%)

    NEXT i%

    TotalLengthWaves = TotalLengthWaves + 0.5## 'add half-wavelength of 1 meter at 299.8 MHz

    Ystart = -TotalLengthWaves/2##

    N% = FREEFILE

    OPEN "PBM5.NEC" FOR OUTPUT AS #N%

        PRINT #N%,"CM File: PBM5.NEC"
        PRINT #N%,"CM Run ID "+DATE$+" "+TIME$
        PRINT #N%,"CM Nd="+STR$(Nd%)+", p="+STR$(p%)+", j="+STR$(j&)
        PRINT #N%,"CM R(p,1,j)="+STR$(R(p%,1,j&))+", R(p,2,j)="+STR$(R(p%,2,j&))
        PRINT #N%,"CE"

        FOR i% = 1 TO Nd%+1

            SumDi = 0##

            FOR q% = 1 TO i%-1

                SumDi = SumDi + Di(q%)

            NEXT q%

            Y1 = ROUND(Ystart + SumDi,6)

            Y2 = ROUND(Y1+0.5##,6) 'add one-half wavelength for other end of dipole

            PRINT #N%,"GW "+REMOVE$(STR$(i%),ANY" ")+",49,0.,"+REMOVE$(STR$(Y1),ANY" ")+",0.,0.,"+REMOVE$(STR$(Y2),ANY" ")+",0.,0.00001"

        NEXT i%

        PRINT #N%,"GE"

        FOR i% = 1 TO Nd%+1
            PRINT #N%,"EX 0,"+REMOVE$(STR$(i%),ANY" ")+",25,0,1.,0." 'VOLTAGE SOURCES
        NEXT i%

        PRINT #N%,"FR 0,1,0,0,299.79564,0."
        PRINT #N%,"RP 0,1,1,1001,90.,0.,0.,0.,1000." 'gain at 1000 wavelengths range
        PRINT #N%,"XQ"
        PRINT #N%,"EN"

    CLOSE #N%

'      - -  ANGLES - -         - POWER GAINS -       - - - POLARIZATION - - -    - - E(THETA) - - -   - - - E(PHI) - - -
'       THETA      PHI     VERT.    HOR.    TOTAL    AXIAL   TILT  SENSE       MAGNITUDE    PHASE     MAGNITUDE     PHASE
'      DEGREES   DEGREES     DB      DB       DB      RATIO   DEG.             VOLTS/M    DEGREES     VOLTS/M     DEGREES
'       90.00      0.00    3.91  -999.99    3.91    0.00000   0.00  LINEAR    1.29504E-04    5.37    0.00000E+00    -5.24
'123456789x123456789x123456789x123456789x123456789x123456789x123456789x123456789x123456789x123456789x123456789x123456789x
'      10        20        30        40        50        60        70        80        90        100       110       120

    SHELL "n41_2k1.exe",0

    N% = FREEFILE

    OPEN "PBM5.OUT" FOR INPUT AS #N%

        WHILE NOT EOF(N%)

            LINE INPUT #N%, Lyne$

            IF INSTR(Lyne$,"DEGREES  DEGREES") > 0 THEN EXIT LOOP

        WEND 'position at next data line

        LINE INPUT #N%, Lyne$

    CLOSE #N%

    GainDB$ = REMOVE$(MID$(Lyne$,37,8),ANY" ")
    PBM_5 = 10^(VAL(GainDB$)/10##) 'Directivity

END FUNCTION 'PBM_5()
'----
FUNCTION LOADED_MONO(R(),Nd%,p%,j&) 'single resistor-loaded monopole

    LOCAL N%, i%, Nsegs%, NumFreqs%, LoadedSegNumber%

    LOCAL Lyne$, FR1$, FileStatus$

    LOCAL h_meters, a_meters, SegLength, LoadResistanceOhms, LoadHeightMeters AS EXT

    LOCAL Zo, FrequencyMHZ(), RadEfficiencyPCT(), MaxGainDBI(), MinGainDBI(), RinOhms(), XinOhms(), VSWR() AS EXT

    LOCAL FoM, MinimumRadiationEfficiency, MinimumMaxGain, MinVSWR, MaxVSWR, MinRin, MaxRin, MinXin, MaxXin AS EXT

    LOCAL StartFreqMHZ!, StopFreqMHZ!, FreqStepMHZ!, SegRes!, SegInd!, SegCap!

    REDIM FrequencyMHZ(1 TO 1), RadEfficiencyPCT(1 TO 1), MaxGainDBI(1 TO 1), MinGainDBI(1 TO 1), RinOhms(1 TO 1), XinOhms(1 TO 1), VSWR(1 TO 1)

    Zo = 50## '300 => 6:1 step-down ratio used by Rao & Debroux

    Nsegs% = 107 : h_meters = 10.7## : a_meters = 0.005##

    SegLength = h_meters/Nsegs%

    LoadResistanceOhms = R(p%,1,j&) : LoadHeightMeters = R(p%,2,j&) : LoadedSegNumber% = MIN(MAX(INT(0.5##+LoadHeightMeters/SegLength),1),Nsegs%)

    StartFreqMHZ! = 5! : StopFreqMHZ! = 30! : FreqStepMHZ! = 1! : NumFreqs% = 1 + (StopFreqMHZ!-StartFreqMHZ!)/FreqStepMHZ!

    FR1$ = "FR 0,"+Int2String$(NumFreqs%)+",0,0,"+FP2String2$(StartFreqMHZ!)+","+FP2String2$(FreqStepMHZ!)

    N% = FREEFILE
```

```
OPEN "LD_MONO.NEC" FOR OUTPUT AS #N%

    PRINT #N%,"CM File: LD_MONO.NEC"
    PRINT #N%,"CM DESIGN #1"
    PRINT #N%,"CM Fitness = [Min(e)+Min(Gmax)]/(MaxVSWR-MinVSWR)"
    PRINT #N%,"CM Run ID "+DATE$+" "+TIME$
    PRINT #N%,"CM Nd="+STR$(Nd%)+",  p="STR$(p%)+",  j="+STR$(j&)
    PRINT# N%,"CM Zo="+REMOVE$(STR$(Zo),ANY" ")+" ohms"
    PRINT #N%,"CE"
    PRINT #N%,"GW1,"+Int2String$(NSegs%)+",0.,0.,0.,0.,0.,"+FP2String$(h_meters)+","+FP2String$(a_meters)
    PRINT #N%,"GE1"  'required for gnd plane connection!

    SegInd! = 0! : SegCap! = 0! : SegRes! = LoadResistanceOhms

    PRINT
#N%,"LD0,1,"+Int2String$(LoadedSegNumber%)+","+Int2String$(LoadedSegNumber%)+","+FP2String$(SegRes!)+","+FP2String$(SegInd!)+","+FP2String$(SegCap
!)

    PRINT #N%,"GN1"

    PRINT #N%, FR1$ ' frequency card

    PRINT #N%, "EX 0,1,1,1,1.,0."

    PRINT #N%, "RP 0,10,1,1001,0.,0.,10.,0.,100000."

    PRINT #N%, "EN"

CLOSE #N%

SHELL "NEC2D_200_02-22-2011.EXE",0

CALL GetNECdata(NumFreqs%,Zo,FrequencyMHZ(),RadEfficiencyPCT(),MaxGainDBI(),MinGainDBI(),RinOhms(),XinOhms(),VSWR(),FileStatus$)

IF FileStatus$ = "OK" THEN

    MinimumRadiationEfficiency = RadEfficiencyPCT(1) : MinimumMaxGain = MaxGainDBI(1) : MinVSWR = VSWR(1) : MaxVSWR = VSWR(1)

    MinRin = RinOhms(1) : MaxRin = RinOhms(1) : Minxin = xinOhms(1) : Maxxin = xinOhms(1)

    FOR i% = 1 TO NumFreqs%
        IF RadEfficiencyPCT(i%) =< MinimumRadiationEfficiency THEN MinimumRadiationEfficiency = RadEfficiencyPCT(i%)
        IF MaxGainDBI(i%)        =< MinimumMaxGain           THEN MinimumMaxGain  = MaxGainDBI(i%)
        IF VSWR(i%)              =< MinVSWR                   THEN MinVSWR         = VSWR(i%)
        IF VSWR(i%)              >= MaxVSWR                   THEN MaxVSWR         = VSWR(i%)
        IF RinOhms(i%)           =< MinRin                   THEN MinRin          = RinOhms(i%)
        IF RinOhms(i%)           >= MaxRin                   THEN MaxRin          = RinOhms(i%)
        IF XinOhms(i%)           =< Minxin                   THEN Minxin          = xinOhms(i%)
        IF XinOhms(i%)           >= Maxxin                   THEN Maxxin          = xinOhms(i%)
    NEXT i%
'    FoM = MinimumRadiationEfficiency + 12##*MinimumMaxGain 'chosen to roughly balance these factors
'    FoM = 1##/(MaxVSWR-MinVSWR)  'chosen to minimize VSWR difference
'    FoM = 2##*MinimumRadiationEfficiency/sqr(((MaxRin-MinRin)^2+(Maxxin-Minxin)^2))
'    FoM = (MinimumRadiationEfficiency+MinimumMaxGain)/(MaxVSWR-MinVSWR)  'DESIGN #1

ELSE

    FoM = -98765##

END IF

N% = FREEFILE
OPEN "LD_MONO.NEC" FOR APPEND AS #N%
    PRINT #N%,""
    PRINT #N%,"FoM ="+STR$(ROUND(FoM,3))
    PRINT #N%,"VSWR Min/Max  ="+STR$(ROUND(MinVSWR,2))+"/"+STR$(ROUND(MaxVSWR,2))+"//"+REMOVE$(STR$(Zo),ANY" ")
    PRINT #N%,"Rin Min/Max   ="+STR$(ROUND(MinRin,2))+"/"+STR$(ROUND(MaxRin,2))
    PRINT #N%,"Xin Min/Max   ="+STR$(ROUND(Minxin,2))+"/"+STR$(ROUND(Maxxin,2))
    PRINT #N%,"MIN Gmax(dbi) ="+STR$(ROUND(MinimumMaxGain,2))
    PRINT #N%,"MIN Eff(%)    ="+STR$(ROUND(MinimumRadiationEfficiency,2))
    PRINT #N%,""
CLOSE #N%

N% = FREEFILE
OPEN "LD_MONO.OUT" FOR APPEND AS #N%
    PRINT #N%,""
    PRINT #N%,"FoM ="+STR$(ROUND(FoM,3))
    PRINT #N%,"VSWR Min/Max  ="+STR$(ROUND(MinVSWR,2))+"/"+STR$(ROUND(MaxVSWR,2))+"//"+REMOVE$(STR$(Zo),ANY" ")
    PRINT #N%,"Rin Min/Max   ="+STR$(ROUND(MinRin,2))+"/"+STR$(ROUND(MaxRin,2))
    PRINT #N%,"Xin Min/Max   ="+STR$(ROUND(Minxin,2))+"/"+STR$(ROUND(Maxxin,2))
    PRINT #N%,"MIN Gmax(dbi) ="+STR$(ROUND(MinimumMaxGain,2))
    PRINT #N%,"MIN Eff(%)    ="+STR$(ROUND(MinimumRadiationEfficiency,2))
    PRINT #N%,""
'   - - - ANGLES - -       - - POWER GAINS -        - - - POLARIZATION - - -  - - - E(THETA) - -    - - - E(PHI) - - -
'  THETA    PHI      VERT.   HOR.    TOTAL    AXIAL    TILT  SENSE      MAGNITUDE    PHASE      MAGNITUDE     PHASE
' DEGREES DEGREES     DB      DB      DB      RATIO    DEG.             VOLTS/M    DEGREES      VOLTS/M     DEGREES
'  90.00    0.00     3.91  -999.99   3.91   0.00000    0.00  LINEAR   1.29504E-04    5.37    0.00000E+00    -5.24
'123456789x123456789x123456789x123456789x123456789x123456789x123456789x123456789x123456789x123456789x123456789x
'       10        20        30        40        50        60        70        80        90       100       110       120

    LOADED_MONO = FoM

END FUNCTION 'LOADED_MONO()

'---------------------

    FUNCTION FP2String2$(X!)
    LOCAL A$
        A$=LTRIM$(RTRIM$(STR$(X!)))
        IF TALLY(A$,".") = 0! THEN A$ = A$ + "."
        FP2String2$ = A$
    END FUNCTION

'---

    FUNCTION FP2String$(X)
    LOCAL A$
        A$=LTRIM$(RTRIM$(STR$(X)))
        IF TALLY(A$,".") = 0## THEN A$ = A$ + "."
        FP2String$ = A$
    END FUNCTION

'---
```

```
    FUNCTION Int2String$(X%)
    LOCAL A$
        A$=LTRIM$(RTRIM$(STR$(X%)))
        Int2String$ = A$
    END FUNCTION

'---

    SUB GetNECdata(NumFreqs%,Zo,FrequencyMHZ(),RadEfficiencyPCT(),MaxGainDBI(),MinGainDBI(),RinOhms(),XinOhms(),VSWR(),FileStatus$)

    LOCAL N%, idx%, AngleNum%

    LOCAL Lyne$, Dum$

    LOCAL GmaxDBI, GminDBI AS EXT

    REDIM FrequencyMHZ(1 TO NumFreqs%),RadEfficiencyPCT(1 TO NumFreqs%),MaxGainDBI(1 TO NumFreqs%),MinGainDBI(1 TO NumFreqs%),RinOhms(1 TO
NumFreqs%),XinOhms(1 TO NumFreqs%),VSWR(1 TO NumFreqs%)

    FileStatus$ = "NOK"

    OPEN "LD_MONO.OUT" FOR INPUT AS #N%

        WHILE NOT EOF(N%)

            LINE INPUT #N%, Lyne$

                IF INSTR(Lyne$,"RUN TIME") > 0 THEN FileStatus$ = "OK"

        WEND

    CLOSE #N%

    IF FileStatus$ <> "OK" THEN EXIT SUB

    OPEN "LD_MONO.OUT" FOR INPUT AS #N%

        idx% = 1

        WHILE NOT EOF(N%)

            LINE INPUT #N%, Lyne$

            IF INSTR(Lyne$,"FREQUENCY=") > 0 THEN
                Lyne$ = REMOVE$(Lyne$,"MHZ") : Lyne$ = REMOVE$(Lyne$,"FREQUENCY= ") : FrequencyMHZ(idx%) = VAL(Lyne$)
'MSGBOX("idx="+STR$(idx%)+"   F="+STR$(FrequencyMHZ(idx%)))
            END IF

            IF INSTR(Lyne$,"INPUT PARAMETERS") > 0 THEN
                LINE INPUT #N%, Dum$ : LINE INPUT #N%, Dum$ : LINE INPUT #N%, Dum$ 'skip three lines
                LINE INPUT #N%, Lyne$ 'input next line with impedance data
                RinOhms(idx%) = VAL(MID$(Lyne$,61,12)) : XinOhms(idx%) = VAL(MID$(Lyne$,73,12)) : VSWR(idx%) =
StandingWaveRatio(Zo,RinOhms(idx%),XinOhms(idx%))
            END IF
            IF INSTR(Lyne$,"EFFICIENCY") > 0 THEN RadEfficiencyPCT(idx%) = VAL(REMOVE$(Lyne$,ANY Alphabet$+" ="))
            IF INSTR(Lyne$,"E(THETA)") > 0 THEN
                LINE INPUT #N%, Dum$ : LINE INPUT #N%, Dum$ 'skip two lines
                GmaxDBI = -9999## : GminDBI = -GmaxDBI
                FOR AngleNum% = 1 TO 10 'WARNING!! MUST USE 10 POLAR ANGLES! HARDWIRED!!
                    LINE INPUT #N%, Lyne$ 'input next TEN lines with pattern data
                    IF VAL(MID$(Lyne$,38,7)) >= GmaxDBI THEN GmaxDBI = VAL(MID$(Lyne$,38,7)) 'get max gain
                    IF (VAL(MID$(Lyne$,38,7)) =< GminDBI AND VAL(MID$(Lyne$,38,7)) > -999.99##) THEN GminDBI = VAL(MID$(Lyne$,38,7)) 'get min gain
                NEXT AngleNum%
                MaxGainDBI(idx%) = GmaxDBI
                MinGainDBI(idx%) = GminDBI
                INCR idx%
            END IF
'msgbox("idx="+STR$(idx%))

        WEND

    CLOSE #N%
```

```
'  TAG  SEG.    VOLTAGE (VOLTS)        CURRENT (AMPS)          IMPEDANCE (OHMS)         ADMITTANCE (MHOS)        POWER
'  NO.   NO.    REAL      IMAG.        REAL      IMAG.         REAL      IMAG.          REAL      IMAG.         (WATTS)
'   1     1  1.00000E+00 0.00000E+00 7.17910E-06 8.93193E-04 8.99811E+00-1.11951E+03 7.17910E-06 8.93193E-04 3.58955E-06
'123456789x123456789x123456789x123456789x123456789x123456789x123456789x123456789x123456789x123456789x123456789x123456789x
'   10       20        30        40        50        60        70        80        90        100       110       120      130
```

```
'  - - ANGLES - -        - - - POWER GAINS - - -    - - - POLARIZATION - - -    - - - E(THETA) - - -     - - - E(PHI) - - -
'  THETA   PHI     VERT.    HOR.    TOTAL    AXIAL   TILT  SENSE    MAGNITUDE   PHASE     MAGNITUDE   PHASE
' DEGREES DEGREES   DB      DB       DB      RATIO   DEG.          VOLTS/M    DEGREES     VOLTS/M    DEGREES
'   0.00    0.00  -999.99 -999.99 -999.99   0.00000  0.00 LINEAR  0.00000E+00  -240.17   0.00000E+00  -240.17
'  10.00    0.00   -18.97 -999.99  -18.97   0.00000  0.00 LINEAR  2.69380E-08   -61.44   0.00000E+00  -240.17
'123456789x123456789x123456789x123456789x123456789x123456789x123456789x123456789x123456789x123456789x123456789x123456789x
'   10       20        30        40        50        60        70        80        90        100       110       120      130
```

```
    END SUB 'GetNECdata()

'------------------------

    FUNCTION StandingWaveRatio(Zo,ReZ,ImZ)

    LOCAL ReRho, ImRho, MagRho, SWR AS EXT

        SWR = 9999##

        CALL ComplexDivide(ReZ-Zo,ImZ,ReZ+Zo,ImZ,ReRho,ImRho)

        MagRho = SQR(ReRho*ReRho+ImRho*ImRho)   'reflection coefficient

        IF MagRho <> 1## THEN SWR=(1##+MagRho)/(1##-MagRho)

        StandingWaveRatio = SWR

    END FUNCTION 'StandingWaveRatio()

'-----

    SUB ComplexMultiply(ReA,ImA,ReB,ImB,ReC,ImC)

'  Returns real and imaginary parts of product C=A*B

        ReC = ReA*ReB-ImA*ImB
        ImC = ImA*ReB+ReA*ImB
    END SUB
```

```
'-----

    SUB ComplexDivide(ReA,ImA,ReB,ImB,ReC,ImC)

'  Returns real and imaginary parts of quotient C=A/B

        LOCAL Denom AS EXT

        Denom = ReB*ReB+ImB*ImB
        ReC = (ReA*ReB+ImA*ImB)/Denom
        ImC = (ImA*ReB-ReA*ImB)/Denom
    END SUB

'************************************************* END PROGRAM *****************************************************
```